%% file: Main.tex
\documentclass[
	,12pt%
	,pagesize%
	,headings=small%
	,paper=a4%
	,parskip=false%
	,abstract=true
	,toc=bibliography
]{scrartcl}

\addtokomafont{sectioning}{\normalfont\bfseries}
\setkomafont{title}{\normalfont}
\setkomafont{subtitle}{\normalfont}
\setkomafont{author}{\normalfont}
\setkomafont{date}{\normalfont}
\addtokomafont{pageheadfoot}{\scshape\small}
\setkomafont{caption}{\footnotesize}
\setkomafont{captionlabel}{\usekomafont{caption}\bfseries}
\setcapindent{0pt}

\usepackage{etoolbox}
\usepackage{xspace}
\usepackage{xpatch}
\usepackage{ifthen}
\usepackage[normalem]{ulem}

\usepackage[utf8]{inputenc}
\usepackage[T1]{fontenc}
\usepackage{csquotes}
\usepackage[english]{babel}
\usepackage{graphicx}
\graphicspath{{Pictures/}{/}}
\usepackage[export]{adjustbox}

\usepackage[dvipsnames]{xcolor}
\usepackage{tikz-cd}

\usepackage[draft]{fixme}
\usepackage[nointlimits]{amsmath}
\usepackage{amssymb,mathrsfs,mathtools}
\usepackage{newtxtext}
\usepackage[varvw]{newtxmath}

\usepackage{upgreek}
\usepackage{cancel}
\usepackage{aliascnt} %alias counter
\usepackage[amsmath,thmmarks,hyperref]{ntheorem}
\usepackage{xkeyval}
\usepackage{multicol}
\usepackage{authblk}
\usepackage{enumitem}

\usepackage[final,%
	bookmarks,%
	bookmarksdepth=3,%
	breaklinks=true,%
	colorlinks=true,%
	urlcolor=NavyBlue,%
	linkcolor=NavyBlue,%
	citecolor=ForestGreen,%	
]{hyperref}%
\usepackage[all]{hypcap}
\usepackage[capitalise,nameinlink]{cleveref}

\crefname{equation}{}{}
\newlist{claims}{enumerate}{10}
\setlist[claims]{label*=\arabic*.,ref=\arabic*}
\crefname{claimsi}{Claim}{Claim}
\Crefname{claimsi}{Claim}{Claims}
\newlist{conditions}{enumerate}{10}
\setlist[conditions]{label*=\arabic*.,ref=\arabic*}
\crefname{conditionsi}{Condition}{Conditions}
\Crefname{conditionsi}{Condition}{Conditions}
\AtBeginEnvironment{abstract}{\footnotesize}
\makeatletter
\patchcmd{\@maketitle}{\huge}{\Large}{}{}
\makeatother
\usepackage[%
	backend = biber%
	,isbn=false
	,doi=true
	,url=false
	,backref=false
	,clearlang=true
	,maxbibnames=5
	,giveninits=true
]{biblatex}
\AtBeginBibliography{%
	\footnotesize%
}
\renewbibmacro*{doi+eprint+url}{%
    \printfield{doi}%
    \newunit\newblock%
    \iftoggle{bbx:eprint}{%
        \usebibmacro{eprint}%
    }{}%
    \newunit\newblock%
    \iffieldundef{doi}{%
        \usebibmacro{url+urldate}}%
        {}%
    }
\DeclareRedundantLanguages{en}{english,american,british,canadian,australian,newzealand,USenglish,UKenglish}
\DeclareRedundantLanguages{English}{english}
\numberwithin{equation}{section}

\input{Compatibility_HS}

\input{Macros_vdM}
\input{Macros_HS}

\input{Macros}
\input{NTheoremEnglish}

\input{Title}

\begin{document}

\maketitle

\begin{abstract}
\input{Abstract}
\end{abstract}
\input{Introduction}

\input{PalaisSmaleAbstractSetting}
\input{ManifoldA_new.tex}
\input{PalaisSmale}
\input{Smoothness}

\appendix
\crefalias{section}{appendix} 
\input{ReparametrizationToArcLength}

\input{ProductRule_BehzadanHolst}
\input{Acknowledgments}
\printbibliography
\end{document}

%% file: Compatibility_HS.tex
\providecommand{\rm}{}

\providecommand{\sc}{}

\renewcommand{\rm}{}

\renewcommand{\sc}{\scshape}

\newcommand{\bz}{\mathbf{z}}
\newcommand{\by}{\mathbf{y}}

\DeclareFontFamily{U}{mathx}{\hyphenchar\font45}
\DeclareFontShape{U}{mathx}{m}{n}{
<-6> mathx5 <6-7> mathx6 <7-8> matha7
<8-9> mathx8 <9-10> mathx9
<10-12> mathx10 <12-> mathx12
}{}
\DeclareSymbolFont{mathx}{U}{mathx}{m}{n}
\DeclareMathSymbol{\bigoasterisk}{\mathop}{mathx}{"C6}

\usepackage{ifthen}
\usepackage{enumitem}

\makeatletter
\newcommand{\Item}[1][]{%\refstepcounter{enumi}%
	\ifthenelse{\equal{#1}{}}{%
		\item%
	}{%
		\refstepcounter{Item}%
		\csname @item\endcsname[#1]%
		\protected@edef\@currentlabel{#1}%
		\protected@edef\cref@currentlabel{[enumi][][]{#1}}%
	}%
}
\makeatother

%% file: Macros_vdM.tex
\DeclareMathOperator{\pr}{pr}
\DeclareMathOperator{\PR}{Pr}

\DeclareMathOperator{\dist}{dist}
\DeclareMathOperator{\ARC}{arc}
\DeclareMathOperator{\CAPARC}{Arc}
\DeclareMathOperator{\intM}{intM}

\newcommand{\Mpq}[1][p,q]{\ensuremath{\intM^{\left(#1\right)}}}

\newcommand{\AL}{\mathcal{A}}
\newcommand{\BL}{\mathcal{B}}
\newcommand{\CL}{\mathcal{C}}

\newcommand{\EL}{\mathcal{E}}
\newcommand{\FL}{\mathcal{F}}
\newcommand{\GL}{\mathcal{G}}
\newcommand{\HL}{\mathcal{H}}

\newcommand{\KL}{\mathcal{K}}
\newcommand{\LL}{\mathcal{L}}
\newcommand{\ML}{\mathcal{M}}
\newcommand{\NL}{\mathcal{N}}
\newcommand{\OL}{\mathcal{O}}

\newcommand{\QL}{\mathcal{Q}}
\newcommand{\RL}{\mathcal{R}}
\newcommand{\SL}{\mathcal{S}}
\newcommand{\TL}{\mathcal{T}}
\newcommand{\UL}{\mathcal{U}}
\newcommand{\VL}{\mathcal{V}}

\newcommand{\FS}{\mathscr{F}}
\newcommand{\PS}{\mathscr{P}}

\newcommand{\Id}{{\mathrm{Id}}}

\DeclareMathOperator{\TP}{TP}
\newcommand{\Mpz}{\Mpq[p,2]}
\newcommand{\TPpq}[1][p,q]{\ensuremath{\TP^{\left(#1\right)}}}
\newcommand{\TPp}{\TPpq[p,2]}

\newcommand{\omitted}[1]{}

\newcommand{\Fo}{\,\,\,\text{for }\,\,}
\newcommand{\Foa}{\,\,\,\text{for all }\,\,}

\newcommand{\AND}{\,\,\,\text{and }\,\,}
\newcommand{\OR}{\,\,\,\text{or }\,\,}

\newcommand{\As}{\,\,\,\text{as }\,\,}

\newcommand{\heikodetail}[1]{}

%% file: Macros_HS.tex
\newcommand{\dd}{\mathop{}\!\mathrm{d}}
\renewcommand{\d}{\dd}
\newcommand{\ceq}{\coloneq}
\newcommand{\qec}{\eqqcolon}

\newcommand{\sdfrac}[2]{\mbox{\small$\displaystyle\frac{#1}{#2}$}}

\DeclarePairedDelimiterXPP{\pars}[1]{\mathop{}}{\lparen}{\rparen}{}{#1}
\DeclarePairedDelimiterXPP{\abs}[1]{\mathop{}}{\lvert}{\rvert}{}{#1}
\DeclarePairedDelimiterXPP{\norm}[1]{\mathop{}}{\lVert}{\rVert}{}{#1}
\DeclarePairedDelimiterXPP{\seminorm}[1]{\mathop{}}{\lbrack}{\rbrack}{}{#1}
\DeclarePairedDelimiterXPP{\inner}[1]{\mathop{}}{\langle}{\rangle}{}{#1}
\DeclarePairedDelimiterXPP{\iinner}[1]{\mathop{}}{\langle\!\langle}{\rangle\!\rangle}{}{#1}
\DeclarePairedDelimiterXPP{\brackets}[1]{\mathop{}}{\lbrack}{\rbrack}{}{#1}
\DeclarePairedDelimiterXPP{\braces}[1]{\mathop{}}{\lbrace}{\rbrace}{}{#1}

\DeclarePairedDelimiterXPP{\floor}[1]{\mathop{}}{\lfloor}{\rfloor}{}{#1}
\DeclarePairedDelimiterXPP{\ceil}[1]{\mathop{}}{\lceil}{\rceil}{}{#1}

\DeclarePairedDelimiterXPP{\intervalcc}[1]{\mathop{}}{\lbrack}{\rbrack}{}{#1}
\DeclarePairedDelimiterXPP{\intervalco}[1]{\mathop{}}{\lbrack}{\rparen}{}{#1}
\DeclarePairedDelimiterXPP{\intervaloc}[1]{\mathop{}}{\lparen}{\rbrack}{}{#1}
\DeclarePairedDelimiterXPP{\intervaloo}[1]{\mathop{}}{\lparen}{\rparen}{}{#1}

\DeclarePairedDelimiterXPP{\myset}[2]{\mathop{}}{\lbrace}{\rbrace}{}{#1\,\delimsize\vert\,\mathopen{}#2}

\newcommand{\Circle}{{\R\slash\Z}}

\newcommand{\AmbDim}{{n}}
\newcommand{\AmbSpace}{{\R^\AmbDim}}
\newcommand{\DomDim}{{m}}
\newcommand{\DomSpace}{\R^\DomDim}
\newcommand{\Domain}{{\varOmega}}

\newcommand{\BendingEnergy}{E_{\smash{\mathrm{b}}}}

\newcommand{\MixedEnergy}{\mathcal{T}}

\newcommand{\qand}{\quad \text{and} \quad}

\newcommand{\inj}{{\smash{\mathrm{i}}}}
\newcommand{\reg}{{\smash{\mathrm{r}}}}
\newcommand{\injreg}{{\smash{\mathrm{ir}}}}
\newcommand{\unitspeed}{{\smash{\mathrm{a}}}}
\newcommand{\injunitspeed}{{\smash{\mathrm{ia}}}}
\newcommand{\homeom}{{\smash{\mathrm{e}}}}
\newcommand{\periodic}{{\smash{\mathrm{p}}}}
\newcommand{\loc}{{\smash{\mathrm{loc}}}}

%% file: Macros.tex
\newcommand\D[1]{{#1(x)-#1(y)}}

\newcommand{\eps}{\ensuremath{\varepsilon}}

\newcommand{\F}[1][\eps]{F}%\ensuremath{F_{#1}}}

\newcommand{\g}{\ensuremath{\gamma}}
\newcommand{\G}{\ensuremath{\varGamma}}

\newcommand{\N}{\ensuremath{\mathbb{N}}}

\newcommand{\R}{\ensuremath{\mathbb{R}}}
\newcommand{\C}{\ensuremath{\mathbb{C}}}
\renewcommand{\rho}{\ensuremath{\varrho}}

\newcommand{\Z}{\ensuremath{\mathbb{Z}}}

\graphicspath{{./mathematica/}}
\newcommand{\intRZ}{\underset{\R/\Z}{\int}}

\newcommand{\iiintRZ}{\underset{(\R/\Z)^3}{\int}}
\renewcommand{\P}{\ensuremath{\,\mathit{P}}}
\renewcommand{\D}{\ensuremath{\,\mathrm{D}}}

\DeclareMathOperator{\T}{\mathrm{T}}
\newcommand{\embeds}{\hookrightarrow}

\newcommand{\Econtrolled}[1]{\hbox{$#1$-controlled}\xspace}
\newcommand{\ControlledSet}{\VL}
\newcommand{\BoundedLinOps}{L}

\newcommand{\PaSma}{(PS)\xspace}

%% file: NTheoremEnglish.tex
\newcommand{\mynewtheorem}[4] %{BEZEICHNER}{COUNTER}{TITEL}{PLURAL} - for correct enumeration with \cref from the cleveref package
{%
\newaliascnt{#1}{#2}%
\newtheorem{#1}[#1]{#3}%
\aliascntresetthe{#1}%
\crefname{#1}{#3}{#4}%
\Crefname{#1}{#3}{#4}%
}

\mynewtheorem{theorem}{basetheorem}{Theorem}{Theorems} 
\mynewtheorem{lemma}{basetheorem}{Lemma}{Lemmata} 
\mynewtheorem{proposition}{basetheorem}{Proposition}{Propositions} 
\mynewtheorem{corollary}{basetheorem}{Corollary}{Corollaries}
\mynewtheorem{prelemma}{basetheorem}{Pre-lemma}{Pre-lemmata} 
\mynewtheorem{problem}{basetheorem}{Problem}{Problems}

\theoremstyle{break}
\mynewtheorem{btheorem}{basetheorem}{Theorem}{Theorems} 
\mynewtheorem{blemma}{basetheorem}{Lemma}{Lemmata}
\mynewtheorem{bproposition}{basetheorem}{Proposition}{Propositions}
\mynewtheorem{bcorollary}{basetheorem}{Corollary}{Corollaries}
\mynewtheorem{bproblem}{basetheorem}{Problem}{Problems}

\theoremstyle{plain}
\theorembodyfont{\normalfont}
\mynewtheorem{definition}{basetheorem}{Definition}{Definitions}
\mynewtheorem{example}{basetheorem}{Example}{Examples}
\mynewtheorem{remark}{basetheorem}{Remark}{Remarks}

\theoremstyle{break}
\mynewtheorem{bdfn}{basetheorem}{Definition}{Definitions}
\mynewtheorem{bex}{basetheorem}{Example}{Examples}
\mynewtheorem{brem}{basetheorem}{Remark}{Remarks}

\theoremheaderfont{\itshape}
\theorembodyfont{\upshape}
\theoremstyle{nonumberplain}
\theoremseparator{.}
\theoremsymbol{\ensuremath{\Box}}	
\newtheorem{proof}{Proof}

\theoremstyle{empty}
\newtheorem{refproof}{Proof}

%% file: Title.tex
\title{On the Palais--Smale condition in geometric knot theory}

\author{Nicolas Freches}
\author{Henrik Schumacher}
\author{Daniel Steenebr\"ugge}
\author{Heiko von der Mosel}
\affil{Institute for Mathematics, RWTH Aachen University}

\date{\today}

%% file: Abstract.tex
We prove that various families of energies relevant in 
geometric knot theory satisfy the Palais--Smale condition (PS)
on submanifolds of arclength para\-metrized knots. 
These energies
include linear combinations of the Euler-Bernoulli
bending energy with a wide variety of non-local
knot energies, such as
O'Hara's self-repulsive potentials $E^{\alpha,p}$, generalized
tangent-point energies $\TP^{(p,q)}$, and generalized integral
Menger curvature functionals $\intM^{(p,q)}$. Even the
tangent-point energies $\TP^{(p,2)}$ for $p\in (4,5)$ alone
are shown to fulfill the (PS)-condition. For all energies mentioned
we can therefore prove existence of minimizing knots in any prescribed
ambient isotopy class, and we
provide long-time existence of their Hilbert-gradient flows, 
and subconvergence to critical knots as time goes to infinity.
In addition, we prove $C^\infty$-smoothness of all arclength-constrained 
critical knots, which shows in particular that these critical knots
are also critical for the energies on the larger open set of
regular knots under a fixed-length constraint.

%% file: Introduction.tex
\section{Introduction}\label{sec:intro}

\subsection{The Palais--Smale condition}\label{sec:palais-smale-cond}
A Palais--Smale sequence (or (PS)-sequence) for a continuously differentiable
energy function $\EL$ on a smooth Banach manifold $\ML$, is a sequence
$(x_k)_k\subset\ML$ such that $\sup_{k\in\N}\EL(x_k)<\infty$ and
$\lim_{k\to\infty}\|\D \EL_{x_k}\|_{(T_{x_k}\ML)^*}=0$, where 
$(T_x\ML)^*$ denotes the dual space of the tangent space $T_x\ML$ of 
the manifold $\ML$ at the point $x\in\ML$.

The energy $\EL$ is said to satisfy the \emph{Palais--Smale condition}
or (PS)-\emph{condition}, if every (PS)-sequence $(x_k)_k\subset\ML$
for $\EL$ contains a convergent subsequence, $x_{k_l}\to x_\infty\in\ML$
as $l\to\infty$.

The (PS)-condition has been central for the investigation of critical
points and of energy landscapes, with sophisticated tools such as
Mountain Pass theorems, 
Ljusternik-Schnirelman theory, or Morse-Bott results; 
see, e.g., the monographs \cite{rabinowitz_1986}, \cite{palais-terng_1988},
\cite{mawhin-willem_1989}, \cite{WillemMinimaxTheorems1996a},
\cite[Chapters II \& III]{struwe_2008a}, and
\cite{jabri_2003}. 

In \emph{geometric knot theory} one considers non-local energies that
penalize self-inter\-sec\-tions of geometric objects like curves or surfaces, 
so called \emph{self-avoidance energies} or \emph{knot energies}.
With these functionals one can treat variational problems in spaces of
embeddings under topological constraints, such as a prescribed knot type.
In this way, energy minimizing curves or Euclidean submanifolds have been
found in given ambient isotopy classes. In addition, bounds on the number of
isotopy types under given energy values were established, and classical knot
invariants such as the crossing number or the stick number can be controlled
in terms of knot energies; see the surveys 
\cite{ohara_2003,strzelecki-etal_2013a,
strzelecki-vdm_2013a,strzelecki-vdm_2014,blatt-reiter_2014,strzelecki-vdm_2018} 
and the many references 
therein. Moreover, critical points
were found by 
symmetric criticality, not only for self-avoidance energies   alone
\cite{cantarella-etal_2014b,gilsbach-vdm_2018,wings_2025}, 
but also for linear combinations of such
with the classic Euler-Bernoulli bending energy
\begin{equation}\label{eq:bending-energy}
        \BendingEnergy(\g) \ceq \frac{1}{2} \int_\g\kappa^2\d s
\end{equation}
for closed curves $\g:\R/\Z\to\R^3$, in search of \emph{elastic knots}
\cite{gerlach-etal_2017,gilsbach-etal_2023}. There are long-time existence
results for different kinds of  gradient flows for various knot energies, the 
$L^2$-gradient flow 
\cite{blatt_2012b,blatt_2018,blatt_2020a}, and Hilbert or Banach gradient flows
\cite{knappmann-etal_2023,matt-etal_2023}, as well as
$L^2$-gradient flows
for linear combinations with the bending energy \cref{eq:bending-energy}
in \cite{lin-schwetlick_2010,vonbrecht-blair_2017}. We should also point
out the substantial advances in the numerical treatment of such 
problems as in 
\cite{bartels-etal_2018,bartels-reiter_2020,bartels-reiter_2021,yu-etal_2021b,yu-etal_2021,sassen-etal_2024}.

However, a more sophisticated critical point theory for knot energies
has not been developed yet, partially because of the lack of the (PS)-condition
for these non-local energies. With this paper we want to close this gap
for a family of knot energies and for linear combinations of basically
all known differentiable knot energies
with the bending energy \cref{eq:bending-energy}.

\subsection{Three families of knot energies on curves}\label{sec:knot-energies}
O'Hara introduced and investigated
in \cite{ohara_1991a,ohara_1992a,ohara_1992b,ohara_1994}
\emph{repulsive potentials} of the form
\begin{equation}\label{eq:E-alpha-p}
E^{\alpha,p}\ceq\int_{\R/\Z}\int_{\R/\Z}\left(\frac1{|\g(u)-\g(v)|^\alpha}-\frac1{d_\g(s,t)^\alpha}\right)^p|\g'(u)||\g'(v)|\d u\d v,
\end{equation}
where $d_\g(u,v)$ denotes the intrinsic distance between the points $\g(u)$
and $\g(v)$ along the curve, regularizing the singularity of the first
summand. The most prominent example is the M\"obius energy
$E^{2,1}$, which  is not only scale-invariant but also invariant under
M\"obius transformations \cite[Theorem 2.1]{freedman-etal_1994}.

In contrast to O'Hara's energies \cref{eq:E-alpha-p} the \emph{(generalized)
tangent-point energies}
\begin{equation}\label{eq:TP}
 \TP^{(p,q)}(\gamma)= \int_{\R/\Z}\int_{-\frac{1}{2}}^\frac{1}{2}
\frac{\big|\P_{\gamma'(u)}^{\perp}(\gamma(u+w)-\gamma(u)) \big|^q}{|
\gamma(u+w)-\gamma(u)|^{p}}|\g'(u)||\g'(u+w)|\d w\d u
\end{equation}
have integrands that are finite for smooth embedded curves $\g \colon \R/\Z \to \R^n$
for a reasonable\footnote{The relevant ranges for the parameters of all
three energy families will be discussed below, in connection with
their respective energy spaces.}
parameter range for $p$ and $q$.
So, there is no need for a regularizing term. Here, $P_v^\perp$ means the
orthogonal projection onto the $(n-1)$-dimensional
linear subspace $v^\perp\subset\R^n$ for a vector $v\in\R^n\setminus\{0\}$.
The energies' name originates in the
geometric classic tangent-point energy $\TP_q\ceq2^q\TP^{(2q,q)}$ whose
integrand equals the $q$-th power of the inverse
tangent-point radius, i.e., the
radius of the unique circle through the points
$\g(u)$ and $\g(u+w)$ that is in addition tangent to the curve
$\g$ in $\g(u)$. 
This energy was suggested by Gonzalez and Maddocks in 
\cite[p. 4773]{gonzalez-maddocks_1999} as an alternative to self-repulsive
potentials, and was first investigated analytically in
\cite{strzelecki-vdm_2012}.

In addition to the tangent-point radius that defines $\TP_q$,
Gonzalez
and Maddocks suggested to also
consider the unique circle through three given curve points, which
leads to integral Menger curvature $\ML_p\ceq2^p\intM^{(p,p)}$
as a specific member of the family of
\emph{(generalized) integral Menger curvatures}
\begin{equation}\label{eq:intM}
\Mpq=\iiintRZ\frac{\left|\gamma'(u) \right|\left|\gamma'(v) \right|\left|\gamma'(w) \right|}{R^{(p,q)}(\gamma(u),\gamma(v),\gamma(w))}\d u \d w \d v,
\end{equation}
where
\begin{equation}\label{eq:integrand-intM}
R^{(p,q)}(x,y,z)\ceq\frac{\left( \left|y-z \right|\left|y-x \right|\left|z-x \right| \right)^p}{\left|(y-x)\wedge(z-x) \right|^q}\, \text{ for }x,y,z\in\R^n;
\end{equation}
see \cite{strzelecki-etal_2010} for first analytical
results on $\ML_p$ and \cite{strzelecki-etal_2013a} for its
impact on knot theory. 

Crucial for our investigations is the knowledge of the
exact \emph{energy spaces} for these families of knot energies.
By this we mean specific Banach spaces $\CL$ such that any arclength
parametrized closed curve has finite energy if and only if
it is embedded and contained in $\CL$.
For each of these energies the underlying energy space\footnote{For $\TP^{(p,q)}$
this was shown in \cite[Theorem 1.1]{blatt-reiter_2015a}, for $\intM^{(p,q)}$
see \cite[Theorem 1]{blatt-reiter_2015b}, and for $E^{\alpha,p}$ we 
refer to \cite[Theorem 1.1]{blatt_2012a}.}
 is a fractional
Sobolev space $W^{1+s,\rho}(\R/\Z,\R^n)$, where
\begin{equation}\label{eq:individual-energy-spaces}
\textstyle
(s,\rho)\ceq\begin{cases}
\big(\frac{\alpha p-1}{2p},2p\big) & \Fo \FL=E^{\alpha,p},\,
\alpha>0,p\ge 1,2<\alpha p < 2p+1.\\
\big(\frac{p-1-q}{q},q\big) & \textstyle
\Fo \FL=\TP^{(p,q)},\,p\in (q+2,2q+1),q>1,\\ 
\big(\frac{3p-2q-2}q,q\big) & \textstyle
\Fo \FL=\intM^{(p,q)},\,
p\in(\frac23 q+1,q+\frac23), q>1.
\end{cases}
\end{equation}
Note that in all these cases, $s\in (0,1)$.
If $\rho=2$ these fractional Sobolev spaces are Hilbert spaces, and we write
$W^{1+s,2}=H^{1+s}$, where we also admit
the case $s=1$ to include the bending energy \cref{eq:bending-energy},
whose energy space equals $H^2(\R/\Z,\R^n)$. In fact, in the present
paper we focus on the case when the energy space of the
leading order energy is a Hilbert space, allowing, e.g., for
the tangent-point energies $\TP^{(p,2)}$ for $p\in (4,5)$,
and integral Menger curvature $\intM^{(p,2)}$ for $p\in (\frac73,\frac83)$,
as well as total energies of the form $\TL \ceq \BendingEnergy + \theta \, \FL$ for some
$\theta>0$, where $\FL$ is any of the energies in
\cref{eq:individual-energy-spaces}.

Since the knot energies and the bending energy are invariant under the
action of the non-compact group of reparametrizations of the domain $\R/\Z$,
we restrict 
them to the set of injective, arclength 
parametrized curves,
\begin{equation}\label{eq:ALs}
        \AL^{1+s}
        \ceq
        \braces[\big]{
                \g\in H^{1+s}(\R/\Z,\R^n):\g_{[0,1)}\,\,\textnormal{injective}, \,|\g'|=1
        }
        ,
\end{equation}
and show in \Cref{thm:arclength-manifold} that $\AL^{1+s}$ is
a smooth submanifold of  $H^{1+s}(\R/\Z,\R^n)$ for any $s\in (\frac{1}{2},1]$.
It contains another smooth submanifold 
\begin{equation}\label{eq:ALsz}
        \textstyle
        \AL^{1+s}_\bz
        \ceq
        \braces[\big]{\g\in\AL^{1+s}:\g(0)=\bz},\quad s\in (\frac{1}{2},1],
\end{equation}
where $\bz\in\R^n$ is an arbitrary fixed point;
see \Cref{thm:arclength-manifold-point}. 
When restricting the energies
to these submanifolds we do not lose any topological information, since
we prove in \Cref{thm:strong-deformation-ALs} that the open subset of all
regular knots in $H^{1+s}(\R/\Z,\R^n)$ strongly deformation retracts
onto the submanifolds $\AL^{1+s}$ and $\AL^{1+s}_\bz$ 
for every $s\in (\frac{1}{2},1]$
and $\bz\in\R^n$. 
In addition, if two knots $\g_0$ and $\g_1$ lie in
the same connected component of $\AL^{1+s}$ then they are ambiently isotopic,
i.e., in the same knot class,
as shown in \Cref{thm:components-knot-classes}.

\subsection{Main results}
We start with the Palais--Smale condition for linear combinations of the
bending energy \cref{eq:bending-energy} with any energy contained in
one of the
three families \cref{eq:E-alpha-p}, \cref{eq:TP}, or \cref{eq:intM}. But we
have to shrink slightly\footnote{This does \emph{not} change the 
respective parameter range
for $E^{\alpha,1}$, or for $\TP^{(p,q)}$ and $\intM^{(p,q)}$ if $q\in (1,2]$.}
the range of parameters by requiring in addition to
\cref{eq:individual-energy-spaces} that
\begin{equation}\label{eq:shrink-range}
\textstyle
\alpha p<p+2\Fo E^{\alpha,p},\,\, p<\frac{3}{2} q +2\Fo
\TP^{(p,q)},\AND p<\frac{5}{6} q+1\Fo \intM^{(p,q)}
\end{equation}
to allow for the compact embedding $H^2\hookrightarrow W^{1+s,\rho}$ for 
each knot energy. 
\begin{theorem}[Bending energy plus knot energy above scale-invariance]
\label[theorem]{thm:bending-knot-energy}
Let $\theta>0$,  $\bz\in\R^n$, and $\FL\in\{E^{\alpha,p},\TP^{(p,q)},
\intM^{(p,q)}\}$, where the respective combinations of
parameters $\alpha$, $p$, and $q$ are chosen
as in \cref{eq:individual-energy-spaces} and \cref{eq:shrink-range}.
Then the total energy $\TL \ceq \BendingEnergy + \theta \, \FL$ restricted to
the arclength submanifold $\AL^2_\bz\subset H^2(\R/\Z,\R^n)$ 
defined in \cref{eq:ALsz}
satisfies the {\rm (PS)}-condition.
\end{theorem}
Notice that \cref{eq:individual-energy-spaces} excludes the respective
scale-invariant energies, whose energy spaces do \emph{not} embed compactly
into $C^1(\R/\Z,\R^n)$. For example, the energy space of the
M\"obius energy $E^{2,1}$ is the Hilbert space
$H^{\frac{3}{2}}(\R/\Z,\R^n)$. To include this energy as part of
the total energy, we
use in our proof
a slightly smaller Banach space that \emph{does} embed
into $C^1$ (similarly as in 
\cite{reiter-schumacher_2021} and \cite{matt-etal_2023}).
\begin{theorem}[Bending energy plus M\"obius energy]
\label[theorem]{thm:bending-moebius-energy}
The total energy $\TL \ceq \BendingEnergy+ \theta \, E^{2,1}$ restricted
to $\AL^2_\bz$ satisfies the {\rm (PS)}-condition for every $\theta>0$ and
$\bz\in\R^n$.
\end{theorem}
We strongly believe that the statement of \Cref{thm:bending-knot-energy}
carries over to all scale-invariant energies, i.e., to $E^{\alpha,p}$
for $\alpha p=2$, to $\TP^{(p,q)}$ for $p=q+2$, and to $\intM^{(p,q)}$
for $p=\frac23 q +1$, as well.

The case of a pure knot energy  is harder to handle, and we focus
here on the tangent-point energies $\TP^{(p,2)}$ for $p\in (4,5)$ whose
energy space $H^{1+\frac{p-3}2}(\R/\Z,\R^n)$ is a Hilbert space.
\begin{theorem}[Tangent-point energies $\TP^{(p,2)}$]
\label{thm:PS-TP}
Let $\bz\in\R^n$ and $p\in (4,5)$ be fixed and set
\begin{equation}\label{eq:s-TP}
        \textstyle
        s(p)\ceq\frac{p-3}{2}.
\end{equation}
Then the tangent-point energy $\TP^{(p,2)}$ restricted to
the arclength submanifold $\AL^{1+s(p)}_\bz\subset H^{1+s(p)}(\R/\Z,\R^n)$
satisfies the {\rm (PS)}-condition.
\end{theorem}
This result carries over to O'Hara's energies $E^{\alpha,1}$ for 
$\alpha\in (2,3)$ and to the integral Menger curvature functionals
$\intM^{(p,2)}$ for $p\in (\frac73,\frac83)$,   where
one has to replace \cref{eq:s-TP} by the respective differentiability
index $s=\frac{\alpha-1}2$ or 
$s=\frac{3}{2} p -3$ of the corresponding energy space;
see \cite{freches_2025}. At this point, it seems feasible to generalize our methods to the non-Hilbert case $ \TPpq $, with suitable restrictions on $ p $ and $ q $. This is currently work in progress (see \cite{frechesetal_2026}). 

As an immediate application of the (PS)-condition we can
prove the existence of minimizing (PS)-sequences and of
minimizing knots for the total energies
$\TL = \BendingEnergy + \theta \, \FL$ of \Cref{thm:bending-knot-energy}
in any connected component of the arclength
submanifold.

\begin{corollary}[Minimizing knots for total energies]
\label{cor:minimizing-knots-total-energy}
Fix a point \,$\bz\in\R^n$ and some $\theta>0$, and let $\TL = \BendingEnergy + \theta \, \FL$
be any of the total energies considered 
in \cref{thm:bending-knot-energy,thm:bending-moebius-energy}. 
Let $\KL$ be a connected component of the arclength submanifold $\AL^2_\bz$.
Then there are a minimizing {\rm (PS)}-sequence and
a minimizing knot for $\TL|_{\KL}$.
The same statements hold true for the whole arclength submanifold, i.e., for $\KL = \AL^2_\bz$.
\end{corollary}
This result contains the existence results in
\cite[Theorem 3.1]{vdm_1999} and \cite[Theorem 1.1]{gilsbach-etal_2023}\footnote{The additional symmetry constraint considered in \cite[Theorem 1.1]{gilsbach-etal_2023} is of lower order and can be omitted to obtain
global minimizers. This was carried
out  in \cite[Theorem 2.1]{gerlach-etal_2017} but there
for the non-smooth ropelength
functional, which  can be viewed 
as the $L^\infty$-limit of $\TP_q^{1/q}$, or of
$\ML_p^{1/p}$ as $q$ or $p$ tends to infinity.}, there
for specific
choices of knot energies
in pursuit of \emph{elastic knots} in any given knot class.  
As mentioned before,
connected components of $\AL^2_\bz$ correspond to
different knot classes as we show in \Cref{thm:components-knot-classes}.

For the pure  tangent-point energies $\TP^{(p,2)}$,
which are even
of class $C^{1,1}_\loc$ on 
 the underlying arclength submanifold of the larger Hilbert space $H^{1+s(p)}
 (\R/\Z,\R^n)$,
we obtain the analogous statement.

\begin{corollary}[Minimizing knots for $\TP^{(p,2)}$]
\label{cor:minimizing-knots-TP}
Let $\bz\in\R^n$, $p\in (4,5)$, and $s(p)$ as in \cref{eq:s-TP}.
Let $\KL$ be a connected component of the arclength submanifold $\AL^{1+s(p)}_\bz$.
Then there are a minimizing {\rm (PS)}-sequence and a minimizing knot for $\TP^{(p,2)}|_{\KL}$.
This is also true for $\KL = \AL^{1+s(p)}_\bz$.
\end{corollary}
The existence of minimizing knots in prescribed knot classes
under a fixed-length constraint
has been proven in \cite[Theorem 1.3]{blatt-reiter_2015a} for the 
more general energy family
$\TP^{(p,q)}$ with parameters $p,q$ as in \cref{eq:individual-energy-spaces}.

If the energies are of class $C^{1,1}_\loc$, then one can apply
Hilbert gradient-flow techniques on the arclength submanifold
to find critical points. This was carried out by Okabe and Schrader for
the pure bending energy $\BendingEnergy$ in \cite{okabe-schrader_2023}
 which inspired the following result.
Since at this point,  the necessary 
Lipschitz-estimates for the first
variation are not available
for \emph{every} knot energy considered above,
we restrict the following result
to the subfamilies of energy functionals $\FL$ with 
\begin{align}\label{eq:subfamilies}
\FL=E^{2,1},&\OR 
\textstyle
\FL\in 
\big\{\intM^{(p,2)} : p\in (\frac{7}{3} , \frac{8}{3}) \big\}
\cup
\big\{\intM^{(p,p)} : p\in (3,6)\big\},
\notag\\
&\OR\FL\in\big\{\TP^{(p,q)} : \textnormal{$p,q$ as in \cref{eq:individual-energy-spaces},
\cref{eq:shrink-range}, and $q\ge 2$}\big\}. 
\end{align}
\begin{theorem}[Gradient flow for total energies]
\label{thm:gradient-flow-total-energy}
        Fix $\theta >0$ and $\bz\in\R^n$, set $\ML\ceq\AL^2_\bz$,
        and consider the restriction
        $\widetilde{\TL}\ceq\TL|_{\ML}$ of the total energy
        $\TL \ceq \BendingEnergy + \theta \, \FL$ to the arclength submanifold $\ML$
        for any choice of $\FL$ in \cref{eq:subfamilies}.
        Then there is a function $\phi:[0,\infty)\times\ML\to\ML$ such
        that $\phi(\cdot,\g_0)\in C^1([0,\infty),\ML)$ for every knot
        $\g_0\in\ML$, that solves the Cauchy-Problem 
        \begin{equation}
                \label{eq:gradient-flow-total-energy}
                \begin{cases}
                        \frac{\d}{\d t}\phi(t,\g_0)& =-\nabla
                        \widetilde{\TL}(\phi(t,
                        \g_0))   
                                \Foa t\in [0,\infty),\\
                        \phi(0,\g_0)&=\g_0.  
                \end{cases}
        \end{equation}
        The function $t\mapsto\widetilde{\TL}(\phi(t,\g_0))$ 
        is strictly decreasing on $[0,\infty)$ unless
        $\g_0$ is 
        a critical knot for $\widetilde{\TL}$.
        For every $t\in [0,\infty)$ the knot $\phi(t,\g_0)$ is ambiently isotopic to $\g_0$, and as
        $t\to\infty$ the solution $\phi(t,\g_0)$ subconverges
        to a critical point $\g_\infty$ of  
        $\widetilde{\TL}$, which is ambiently isotopic
        to $\g_0$ as well. 
\end{theorem}  

This Hilbert gradient flow result is new, and it 
complements the long-time existence theorems
for the $L^2$-gradient flow of $\BendingEnergy+E^{2,1}$ in 
\cite[Theorem 1]{lin-schwetlick_2010},
as well as the arclength-constrained $L^2$-gradient flow in
\cite[Theorem 4.13]{vonbrecht-blair_2017} for total energies with a certain
family
of self-repulsive
potentials containing $E^{2,1}$.

Note, that for geometric energies, one usually investigates a natural gradient flow with respect to some geometric metric, often denoted $ L^2(\dd s) $, or $ H^2(\dd s) $. However, since we restrict the $ H^2 $-flow \cref{eq:gradient-flow-total-energy} to the arclength manifold, it coincides with the projected (i.e. intrinsic) $ H^2(\d s) $-gradient flow.

The subfamilies of energies in \cref{eq:subfamilies} contain 
the tangent-point energies $\TP^{(p,2)}$ for $p\in (4,5)$. For these
pure knot energies we can
also establish long-time existence for the Hilbert gradient flow
and its subconvergence
to critical points.
\begin{theorem}[Gradient flow for $\TP^{(p,2)}$]
\label{thm:gradient-flow-TP}
For fixed $\bz \in \R^n$, $p\in (4,5)$, and $s(p)$ as in \cref{eq:s-TP} set 
$\ML\ceq\AL^{1+s(p)}_\bz\subset H^{1+s(p)}(\R/\Z,\R^n)$. 
Then for  any given initial knot $\g_0\in\ML$
we have long-time existence for the gradient
flow \cref{eq:gradient-flow-total-energy} on $\ML$ 
within the same knot class $[\g_0]$ (with $\widetilde{\TL}$ 
replaced by the restriction $\widetilde{\TP}\ceq
\TP^{(p,2)}|_{\ML}$), along which $\widetilde{\TP}$ is strictly
decreasing unless the initial knot $\g_0$ is 
a $\widetilde{\TP}$-critical point. Moreover,
the flow subconverges to a $\widetilde{\TP}$-critical knot $\g_\infty$
of the  same
knot type as $\g_0$ as time
tends to infinity. 
\end{theorem}
This result carries over to the integral Menger curvatures
$\intM^{(p,2)}$ for $p\in (\frac{7}{3},\frac{8}{3})$ and to $ E^{\alpha,1} $ for $ \alpha \in (2,3) $, \cite{freches_2025}.
In contrast to the long-time existence result for a projected 
Hilbert-gradient
flow for $\intM^{(p,2)}$ in
\cite[Theorem 1.1]{knappmann-etal_2023} we do not need to project
the flow to control the parametrizations, and we do have subconvergence
of the flow
to a critical knot. Very recently, D\"ohrer in cooperation with the
first author proved a {\L}ojasiewicz--Simon inequality, which even leads
to full 
convergence of the flow for $\TP^{(p,2)}$ in \Cref{thm:gradient-flow-TP}; 
see \cite[Theorem 1.5]{doehrer-freches_2025}.
In addition, the gradient flows in \cref{thm:gradient-flow-total-energy,thm:gradient-flow-TP} preserve
symmetries, induced, for instance, by 
a group of isometries acting on the arclength submanifold; see 
\Cref{prop:symmetry-preserving-grad-flow} where we address such symmetries
in an abstract setting. 

In view of the existing regularity results\footnote{ 
\cite[Theorem 1.5]{blatt-reiter_2015a}, \cite[Theorem 4]{blatt-reiter_2015b}, \cite[Theorem 1.2]{blatt-reiter_2013}}
for critical points of knot energies
on the open subset of regular knots under a fixed-length
constraint 
one may ask whether also the more constrained arclength critical knots obtained in
\cref{cor:minimizing-knots-total-energy,cor:minimizing-knots-TP} and \cref{thm:gradient-flow-total-energy,thm:gradient-flow-TP} are smooth. 
We prove in \cref{thm:critical-points} in an abstract setting that
arclength constrained critical knots are actually also critical knots for the more
relaxed variational problem of fixed-length, but under the assumption of 
one more
order of differentiability of the critical knot. 
 For the concrete
energies treated above  we provide this additional regularity
by carefully analyzing the additional Lagrange parameter in the variational equation.
In fact, we prove full $C^\infty$-smoothness of critical knots
for all non-negative linear combinations
of bending energy $\BendingEnergy$ and tangent-point energies $\TP^{(p,2)}$ including
the cases when only one of these energies appears.
\begin{theorem}[Smoothness of critical points of total energies]
\label{thm:smoothness}
Let $\TL_{\nu\theta}\ceq\nu \, \BendingEnergy + \theta \, \TP^{(p,2)}$ for 
 $p\in (4,5)$, $\nu,\theta\ge 0$ with $\nu+\theta>0$,
and $\ML\in\{\AL^{1+s(\nu)},\AL^{1+s(\nu)}_\bz\}$ for a fixed point 
$\bz\in\R^n$, where $s(\nu)=1$ if $\nu >0$ and $s(\nu)=s(p)$ as in \cref{eq:s-TP} if
$\nu=0$. Then   
each critical knot $\g\in\ML$
of the restricted energy $\TL_{\nu\theta}|_{\ML}$ is of class $C^\infty(\R/\Z,\R^n)$.
\end{theorem}
If $ \nu>0 $ and $ \theta=0 $, then also $ \BendingEnergy $-critical immersions with possible self-intersections are smooth.\\
\Cref{thm:smoothness} also holds if one replaces the tangent-point energies
$\TP^{(p,2)}$ by the integral Menger curvature functionals $\intM^{(p,2)}$
for $p\in (\frac{7}{3},\frac{8}{3})$, or by O'Hara's energies
$E^{\alpha,1}$ for $\alpha\in (2,3)$; see \cite{freches_2025}. 
Critical points of these three energy families
are actually  also real-analytic, 
which follows from \Cref{thm:critical-points} and the regularity proofs
for the length-constrained critical points of  O'Hara's energies
in \cite{vorderobermeier_2020}, or of
integral Menger curvature in \cite{steenebruegge-vorderobermeier_2022}, the methods of which transfer to the tangent-point energies.

\subsection{Outline of the paper}
In \Cref{sec:PS-abstract} we develop sufficient criteria for the  validity of the (PS)-condition
in the abstract setting of general submanifolds of Hilbert spaces, not only
for a single energy as in \Cref{thm:PS-abstract} but also for total energies with a leading order term and
some additional energy of lower order; see \Cref{cor:PS-total-energy}. Moreover, with
\Cref{thm:gradflow_subconvergence} we
provide a long-time existence result for the gradient flow and subconvergence
to critical points  for energies of class $C^{1,1}_\loc$. 
That this flow preserves symmetry is shown in 
\Cref{prop:symmetry-preserving-grad-flow}. Then we
deduce as simple
consequences the existence of minimizing (PS)-sequences and energy minimizers;
see \cref{koro:critpoint_lowerenergy,koro:minimizer_existence}, as well as \cref{cor:critical+minimizers}
 for an alternative proof 
that works for energies that are only $C^1$. \Cref{sec:arclength} contains the proofs
that the sets $\AL^{1+s}$ and $\AL^{1+s}_\bz$ defined in \cref{eq:ALs}
and \cref{eq:ALsz} are actually smooth submanifolds of the Hilbert space
$H^{1+s}(\R/\Z,\R^n)$ for any $s\in (\frac{1}{2},1]$ and
$\bz\in\R^n$ 
(\Cref{thm:arclength-manifold} and \Cref{thm:arclength-manifold-point}).
We also define suitable projections onto the respective tangent spaces that are needed
in the abstract \Cref{thm:PS-abstract}. Then we prove in \Cref{sec:AL-toplogy} that the open set
of regular knots strongly deformation retracts onto the arclength submanifolds (see \Cref{thm:strong-deformation-ALs}), 
and that connected components
of the arclength submanifolds correspond to 
different knot classes (\Cref{thm:components-knot-classes}). 
In \Cref{sec:AL-critical} we prove that
arclength constrained critical knots are also critical under the more
relaxed fixed-length constraint, still in the abstract general setting; see \Cref{thm:critical-points} and \Cref{cor:critical-points}.

\Cref{sec:PS-bending+knot-energies} is dedicated to the proofs of
our main results on the Palais--Smale condition. We start in  \Cref{sec:PS-bending}, however,
 with the pure bending energy $\BendingEnergy$
as a prototype example, both for the verification of the (PS)-condition (\Cref{cor:PS_elastic}) and its consequences (arclength-constrained $\BendingEnergy$-critical 
points via minimization
or gradient flows, see \Cref{cor:elastic_crit_gradient-flow}). \Cref{sec:PS-bending} is meant as a
comprehensible blueprint for the more complicated proofs involving non-local terms
in the following \cref{sec:PS-elastic-knot} on total energies, and in
\cref{sec:tanpoint} on the pure tangent-point energies. 
\Cref{thm:bending-knot-energy,thm:bending-moebius-energy} are shown in \cref{sec:PS-elastic-knot}, whereas \Cref{thm:PS-TP} is proven in \cref{sec:tanpoint}.

Finally, \cref{sec:Smoothness} is devoted to the proof of \cref{thm:smoothness}. 
This regularity theorem follows from the combination of the smoothness result for critical points of the pure bending energy $ \BendingEnergy $ (\cref{thm:elastic-smoothness-new}) with those for critical knots of total energies where $ \BendingEnergy $ is the leading order term (\cref{thm:smoothness-T-critical-new}) and of the pure tangent-point energies $ \TPp $ (\cref{thm:smoothness-TP-new}). 

For notational convenience, we use combinations of the indices $\inj,\reg,\mathrm{a}$ to denote \underline{i}njective,
\underline{r}egular, or \underline{a}rclength parametrized
curves in a function space 
$\mathscr{F}$, e.g., $\FS_\inj,$ $\FS_\reg,$ $\FS_\injreg, $
or $\FS_\injunitspeed$.
Furthermore, we denote by $ L(V,W) $ the space of bounded linear maps between normed spaces, abbreviating $ L(V,V)=L(V) $. Unless otherwise specified, $ \left\langle \cdot,\cdot \right\rangle  $ and $ \left|\cdot \right| $, refer to the Euclidean inner product and its induced Euclidean norm, respectively.

%% file: PalaisSmaleAbstractSetting.tex
\section{The Palais--Smale condition in an abstract setting}\label{sec:PS-abstract}
\subsection{Energies on submanifolds of Hilbert spaces}
\label{sec:energies-hilbert}
In this section we provide sufficient criteria on the energy and on the underlying
spaces guaranteeing the validity of the Palais--Smale condition (PS). 
Throughout this section we use the following notation:
We let $\HL$ be a Hilbert space and $\BL$ be a Banach space such that there is a compact embedding $\iota_{\BL} \colon \HL \embeds \BL$.
We suppose that $\EL \colon \OL \to \R$ is an energy functional of class $C^1$
on an open subset $\OL \subset \HL$.
Moreover, we suppose that $\ML \subset \OL$ is a $C^1$-submanifold without boundary.
Summarized, we have the following setup:
\begin{enumerate}
	\Item[(E)]\label{condE}
	{\sc (Embedding)}
	$\ML\subset \OL\subset \HL \xhookrightarrow{\iota_\BL} \BL$ and $\iota_\BL$ is compact.
\end{enumerate}
Since $\ML \subset \HL$ is a $C^1$-submanifold, the tangent bundle $T \ML \subset \ML \times \HL$ is splitting, see \cite[Chapter III]{lang_2002}. More precisely, there is a continuous map $\pr^{\ML} \colon \ML \to \BoundedLinOps(\HL)$ such that for each $x \in \ML$ the operator $\pr^{\ML}_x \ceq \pr^{\ML}(x)$ is a continuous linear projector onto $T_x \ML$. Note that we do \emph{not} mean $\pr^{\ML}_x$ to be the orthonormal projector with respect to the inner product on $\HL$.

Our aim is to find sufficient criteria for the \PaSma-condition of $\widetilde{\EL} \ceq \EL|_{\ML}$.
To this end, we make further assumptions.
This will be eased by introducing the following notion:
We say that a subset $\ControlledSet \subset \ML$ is \emph{\Econtrolled{\EL}} if and only if $\ControlledSet \subset \HL$ is bounded and if $\EL$ is bounded on $\ControlledSet$.
Here are our assumptions:
\begin{enumerate}
	\Item[\rm (P)]\label{condP}
	{\sc (Projector)}
	The map $\pr^{\ML}$ is bounded on \Econtrolled{\EL} subsets of $\ML$.
	\Item[\rm (B)]\label{condB}
	{\sc (Boundedness)}
	Every  \PaSma-sequence for $\widetilde{\EL} = \EL|_{\ML}$ is bounded.
	\Item[\rm (L)]\label{condL}
	{\sc (Limit in $\ML$)}
	Every convergent \PaSma-sequence of $\widetilde{\EL}$ has its limit point in $\ML$.
	\Item[\rm (S)]\label{condS}
	{\sc (Splitting)}
	There are 
	a scalar function $\mu \colon \ML \to \intervaloo{0,\infty}$,
	a continuous bilinear form $\QL \colon \HL \times \HL \to \R$,
	and 
	a map $\RL \colon \ML \to \BL^*$
	such that the differential $\D\EL$ splits as follows:
	\[
		\D\EL_x h = \mu(x) \, \QL(x,h) + \RL_x \, \iota_\BL \, h 
		\quad 
		\text{for all $x \in \ML$ and all $h \in\HL$.}
	\]
	Moreover, we assume the following:
	\begin{itemize}
	\Item[\rm (S1)]\label{condS1}
	For every \Econtrolled{\EL} subset $\ControlledSet \subset \ML$ there are $0 < m \leq M < \infty$ such that 
	\[
		m \leq \mu(x) \leq M \quad \text{for all $x\in\ControlledSet$.}
	\]
	
	\Item[\rm (S2)]\label{condS2}
		$\QL$ is coercive, i.e., there is $c > 0$ such that 
	\[
		\QL(h,h) \ge c \norm{h}_\HL^2
		\quad 
		\text{for all $h \in \HL$.}
	\]
	
	\Item[\rm (S3)]\label{condS3}
	$\RL$ is bounded on \Econtrolled{\EL} subsets of $\ML$.
	\end{itemize}

	\Item[\rm (D)]\label{condD}
	{\sc (Differential in non-tangential directions)}
	There is a map $\SL \colon \ML \to \BL^*$ such that 
	\[
		\D{\EL}_x \pars{\Id_\HL-\pr_x^\ML} \, h = \SL_x \, \iota_\BL \, h
		\quad 
		\text{for all $x \in \ML$ and all $h \in \HL$}, 
	\]
	and $\SL$ is bounded on \Econtrolled{\EL} subsets of $ \ML $.
\end{enumerate}

Condition \ref{condS3} is an abstract way of saying $\RL_x$ is of ``lower order'' and thus of less importance compared to $\QL(x,\cdot)$.
Likewise, \ref{condD} states that $\D \EL_x h$ is of ``lower order'' if $h$ is transversal to $T_x \ML$ in the sense that $h$ lies in the range of the $\pars{\Id_\HL-\pr_x^\ML}$.

\begin{theorem}[\PaSma-condition for $\widetilde{\EL}=\EL|_\ML$]
\label{thm:PS-abstract}
Under the conditions \ref{condE}, \ref{condP}, \ref{condB}, \ref{condL}, \ref{condS}, and \ref{condD} above, 
the restricted
energy $\widetilde{\EL}=\EL|_\ML$ satisfies the \PaSma-condition.
\end{theorem}

\begin{proof}
	Let $(x_k)_k\subset\ML$ be a \PaSma-sequence for $\widetilde{\EL}$.
	By condition \ref{condB}, $\ControlledSet \ceq \braces{x_k}_{k \in \N}$ is a bounded subset of $ \HL $. So,
	passing over to a subsequence, we may therefore assume that $(x_{k})_{k \in \N}$ is a weakly convergent \PaSma-sequence. 
	We are going to prove that $(x_{k})_{k \in \N}$ must also be strongly convergent in $\HL$ to a point $x_\infty$;
	because condition \ref{condL} then implies that $x_\infty$ lies in $\ML$, which proves the theorem.

	To show that $(x_{k})_{k \in \N}$ is strongly convergent, we use an argument inspired by the proof of \cite[Lemma 2.1]{cerami-etal_1984}.
	Note that $\ControlledSet \ceq \braces{x_k : k \in \N}$ is \Econtrolled{\EL} by definition of a \PaSma-sequence. 

	Let $0 < m \leq M < 0$ be the constants from \ref{condS1}, and let $c > 0$ be the coercivity constant from \ref{condS2}.
	For $h_k \ceq x_k - x_\infty$ we get
	\begin{equation}
		c \, m \norm{h_k}_\HL^2
		\leq 
		\mu(x_k) \, Q(h_k,h_k)
		=
		\mu(x_k) \, Q(x_k,h_k) - \mu(x_k) \, Q(x_\infty,h_k)
		.
		\label{eq:ExploitingCoercivity}
	\end{equation}
	From \ref{condS} and \ref{condD}, we obtain
	\[
		\mu(x_k) \, Q(x_k,h_k)
		=
		\D \EL_{x_k} \pr^\ML_{x_k} h_k
		+
		\SL_{x_k} \, \iota_\BL \, h_k
		-
		\RL_{x_k} \, \iota_\BL \, h_k
		.
	\]
	The triangle inequality and condition \ref{condS1} lead to
	\begin{align*}
		c \, m \norm{h_{k}}_\HL^2
		&\leq
		\abs{\D \EL_{x_k} \pr^\ML_{x_k} h_k}
		+
		\pars[\big]{
			\norm{\SL_{x_k}}_{\BL^*}
			+
			\norm{\RL_{x_k}}_{\BL^*} 
		} \norm{\iota_{\BL}h_k}_{\BL}
		+
		M \abs{Q(x_\infty,h_k)}
		.
	\end{align*}
	Because $\VL$ is \Econtrolled{\EL},
	conditions \ref{condD}, and \ref{condS3} imply that 
	\[
		C_\SL \ceq \sup_{k \in \N} \norm{\SL_{x_k}}_{\BL^*} < \infty 
		,
		\qand 
		C_\RL \ceq \sup_{k \in \N} \norm{\RL_{x_k}}_{\BL^*} < \infty 
		.
	\]
	Thus, we obtain
	\begin{equation}
		\label{eq:StruweTrick1}
		c \, m \norm{h_{k}}_\HL^2
		\leq
		\abs{\D \EL_{x_k} \pr^\ML_{x_k} h_k}
		+
		\pars[\big]{ C_\SL + C_\RL }
		\norm{\iota_\BL \, h_k}_{\BL}
		+
		M \abs{Q(x_\infty,h_k)}
		.
	\end{equation}
	The bilinear form $Q$ is continuous, and we have weak convergence $h_k \rightharpoonup 0$ in~$\HL$.
	Hence, the third summand $M\abs{Q(x_\infty,h_k)}$ in \cref{eq:StruweTrick1} converges to $0$ as $k \to \infty$.
	The operator $\iota_\BL$ is compact and $h_k$ converges weakly to $0$ in~$\HL$. 
	Hence, $\iota_\BL \, h_k$ converges strongly to~$0$ in~$\BL$.
	Thus, also the second summand in \cref{eq:StruweTrick1} converges to~$0$.
	Finally, we have to discuss the first summand in \cref{eq:StruweTrick1}.
	By definition of the dual norm and the operator norm $\norm{\pr^\ML_{x_k}}_{\LL(\HL)}$, we have
	\[
		\abs{\D \EL_{x_k} \pr^\ML_{x_k} h_k}
		=
		\abs{\D \widetilde{\EL}_{x_k} \pr^\ML_{x_k} h_k}
		\leq
		\norm{\D \widetilde{\EL}_{x_k}}_{(T_{x_k}\ML)^*}
		\norm{\pr^\ML_{x_k}}_{\LL(\HL)}
		\norm{h_k}_{\HL}
		.
	\]	
	Because $(x_k)_{k \in \N}$ is a \PaSma-sequence for $\widetilde{\EL}$, we have $\norm{\D \widetilde{\EL}_{x_k}}_{(T_{x_k} \ML)^*} \to 0$ as $k \to \infty$. 
	Since $\VL$ is \Econtrolled{\EL} and because of \ref{condP},
	we have $\sup_{k \in \N} \norm{\pr^\ML_{x_k}}_{\LL(\HL)} < \infty$. In addition, the weakly convergent sequence $ (h_k)_{k} $ is bounded.
	This shows that also the first summand in \cref{eq:StruweTrick1} converges to $0$. Since $c \, m > 0$, we infer that $h_k \to 0$ strongly in $\HL$, thus $x_k \to x_\infty$ as $k \to \infty$.
\end{proof}

If $\ML= \OL \subset \HL$,
then assumptions \ref{condP}, and \ref{condD} are automatically satisfied, and we obtain a result similar in flavor to \cite[Chapter~II, Proposition 2.2]{struwe_2008a}. 
(Note that the Lax--Milgram theorem implies that coercivity of $\QL$ means that the operator $\HL \to \HL^*$, $x \mapsto \QL(x,\cdot)$ is continuously invertible.)

\begin{corollary}[\PaSma-condition for $\EL$ on $\OL$]
\label{cor:PS-open-set}
Let $\ML=\OL \subset \HL$ be an open set and let $\EL\in C^1(\OL)$. 
Suppose that \ref{condE}, \ref{condB}, \ref{condL},  and  \ref{condS} are satisfied.
Then $\EL$ satisfies the \PaSma-condition.
\end{corollary}

In view of applications in geometric knot theory we consider linear 
combinations of a non-negative energy functional with lower order
terms. 
For this, we extend the abstract setting \ref{condE} by another Banach space $\CL$, into which the Hilbert space $\HL$ embeds compactly, i.e., there is a compact linear embedding $\iota_\CL \colon \HL \to \CL$.
Moreover, we complement condition \ref{condP} on the projection $\pr_x^\ML \colon \HL\to T_x\ML$
by the following condition:

\begin{enumerate}
	\Item[\rm (P*)]\label{condP*}
	There is an operator-valued map $\Pr \colon \ML \to \BoundedLinOps(\CL)$
	such that 
	\[
		\iota_\CL \, \pr^{\ML}_x = \PR_x \, \iota_\CL
		\quad 
		\text{for all $x \in \ML$,}
	\]
	and the map $\Pr$ is bounded on \Econtrolled{\EL} subsets of $\ML$.
\end{enumerate}

In addition to the original energy $\EL\in C^1(\OL)$ with its restriction $\widetilde{\EL}=\EL|_\ML$ to the submanifold $\ML\subset\OL\subset\HL \embeds \BL$, 
we
consider another functional $\FL \in C^1(\UL)$ on an open subset $\UL\subset\CL$ with $\iota_\CL(\OL)\subset\UL$.
For some fixed $\theta>0$ we investigate the  \emph{total energy} $\TL \ceq \EL + \theta \pars{\FL \circ \iota_\CL} \in C^1(\OL)$. 
Note that, in general, a \PaSma-sequence  for the restriction $\widetilde{\TL} \ceq \TL|_\ML$, need not be a \PaSma-sequence for the individual
functionals $\widetilde{\EL}$ or $\widetilde{\FL} \ceq (\FL \circ \iota_\CL) |_\ML$, even if
both functionals are non-negative, since their differentials
may not tend to zero. To formulate a clear-cut result,
we replace assumptions \ref{condL} and
\ref{condB} by the following separate  
assumptions on the energies: 
\begin{enumerate}
	\Item[\rm (B*)]\label{condB*}
	For every $K \in \R$ the sublevel set $\braces[\big]{ x \in \ML : \EL(x) \leq K}$ of $\widetilde{\EL} = \EL_\ML$ is bounded in $\HL$.
	\Item[\rm (L*)]\label{condL*}
	Whenever a sequence $(x_k)_k\subset\ML$ is weakly convergent in $\HL$ to some point $x_\infty$ and $ \sup_{k\in\N}\FL(x_k)<\infty $,
	then $x_\infty \in \ML$.	
\end{enumerate}

\begin{corollary}[(PS)-condition for the total energy $\widetilde{\TL}=\TL|_\ML$]
\label{cor:PS-total-energy}
	Let $\iota_\CL \colon \HL \embeds \CL$ be a compact embedding.
	Let $\OL \subset \HL$ and $\UL \subset \CL$ be open subsets such that $\iota_\CL(\OL) \subset \UL$, and let $\EL \in C^1(\OL)$ and $\FL \in C^1(\UL)$ be bounded from below on $\ML$.
	Suppose that conditions \ref{condE},
	\ref{condP}, \ref{condP*},
	\ref{condB*}, \ref{condS}, \ref{condD},
	and \ref{condL*} are satisfied.
	Then for every $\theta>0$ the restriction  $\widetilde{\TL}=\TL|_\ML$ of the total energy $\TL = \EL + \theta \pars{\FL\circ\iota_\CL}$
	satisfies the \PaSma-condition.
\end{corollary}

\begin{proof}
	By adding suitable constants, we may assume that $\widetilde{\EL} \geq 0$ and $\widetilde{\FL} \geq 0$.
	Let $(x_k)_k\subset\ML$ be a \PaSma-sequence for $\widetilde{\TL}$.
	By definition, there is a $K \in \R$ such that $\widetilde{\TL}(x_k) \leq K$ for all $k \in \N$.
	Now $\widetilde{\FL} \geq 0$ implies
	\begin{equation}
		\label{eq:energy-bound}
		\EL(x_k)
		=
		\widetilde{\EL}(x_k)
		\le
		\widetilde{\EL}(x_k) + \theta \, \widetilde{\FL}(x_k)=\widetilde{\TL}(x_k)\le K\quad\Foa k\in\N
		.
	\end{equation}
	By \ref{condB*}, $(x_k)_{k \in \N}$ is bounded in $\HL$. Hence, it contains a weakly convergent subsequence with limit point $x_\infty \in \HL$. 
	As in the proof of \cref{thm:PS-abstract}, we pass over to this subsequence and denote it by $(x_k)_{k \in \N}$. 
	And as in that proof, we aim to show that the $x_k$ converge strongly to $x_\infty$ in $\HL$.
	Once this is established, condition \ref{condL*} implies $x_\infty \in \ML$
	since $(x_{k})_{k \in \N}$ is bounded under~$\FL$;
	this follows from $\widetilde{\EL} \geq 0$:
	\begin{equation*}
		\FL(x_k)
		=
		\widetilde{\FL}(x_k)
		\le
		\tfrac{1}{\theta}
		\pars[\big]{ \widetilde{\EL}(x_k) + \theta \, \widetilde{\FL}(x_k)}
		=
		\tfrac{1}{\theta} \, \widetilde{\TL}(x_k)
		\le 
		\tfrac{K}{\theta}
		\quad\Foa k\in\N.
	\end{equation*}
	To show strong convergence of $(x_k)_{k \in \N}$, we first recall from \cref{eq:energy-bound} that $\ControlledSet \ceq \braces{x_k}_{k \in \N}$  is \Econtrolled{\EL}.
	Precisely as in the proof of \cref{thm:PS-abstract}, we derive inequality \cref{eq:StruweTrick1} for $h_k \ceq x_k - x_\infty$ and show that the second and third summand in \cref{eq:StruweTrick1} converge to~$0$. Alas, we need a new argument for showing that the first summand $\abs{\D \EL_{x_k} \pr^\ML_{x_k} h_k}$ converges to~$0$ as $ k\to \infty $.
	This time, we decompose it and use condition \ref{condP*} as follows:
	\[
		\D \EL_{x_k} \pr^\ML_{x_k} h_k
		=
		\D \TL_{x_k} \pr^\ML_{x_k} h_k
		-
		\theta \, \D \FL_{x_k} \, \iota_\CL \, \pr^\ML_{x_k} h_k
		=
		\D \widetilde{\TL}_{x_k} \pr^\ML_{x_k} h_k
		-
		\theta \, \D \FL_{(\iota_\CL x_k)} \, \PR_{x_k}\, \iota_\CL \,h_k
		.
	\]
	This leads us to
	\[
		\abs{\D \EL_{x_k} \pr^\ML_{x_k} h_k}
		\leq 
		\norm{\D \widetilde{\TL}_{x_k}}_{(T_{x_k}\ML)^*}
		\norm{\pr_{x_k}^\ML}_{\LL(\HL)}
		\norm{h_k}_{\HL}
		+
		\theta \norm{\D \FL_{(\iota_\CL x_k)}}_{\CL^*}
		\norm{\PR_{x_k}}_{\LL(\CL)}
		\norm{\iota_\CL \,h_k}_{\CL}
		.
	\]
	Since $\iota_\CL$ is compact, we have $\iota_\CL \, x_k \to \iota_\CL \, x_\infty$ as $k \to \infty$. 
	By condition~\ref{condL*}, we have $x_\infty \in \ML \subset \OL$
	and $\iota_\CL x_\infty \in \iota_\CL(\OL) \subset \UL$.
	Since $x \mapsto \D \FL_x$ is continuous on $ \UL $, and by condition~\ref{condP}, and because of $\iota_\CL \, x_k \to \iota_\CL \, x_\infty$, we have
	\[
		\smash{\sup_{k \in \N}} \norm{\D \FL_{(\iota_\CL x_k)}} < \infty
		,
		\quad 
		\sup_{k \in \N} \norm{\pr_{x_k}^{\ML}}  < \infty
		,
		\qand 
		\smash{\lim_{k \to \N}} \norm{\iota_\CL \, h_k}_{\CL} = 0
		.
	\]
	The sequence $(x_k)_{k\in\N}$ is a \PaSma-sequence for $\smash{\widetilde{\TL}}$, and because of \ref{condP*}, and $h_k \rightharpoonup 0$, we have 
	\[
		\smash{\lim_{k \to \infty}}
		\norm{\D \widetilde{\TL}_{x_k}}_{(T_{x_k}\ML)^*}
		=
		0
		,
		\quad 
		\sup_{k \in \N} \norm{\PR_{x_k}}_{\LL(\CL)}  < \infty
		,
		\qand
		\smash{\sup_{k \in \N}} \norm{h_k}_{\HL}  < \infty
		.
	\]
	All this together implies that $\smash{\abs{\D \widetilde{\EL}_{x_k} \pr^\ML_{x_k} h_k}} \to 0$ as $k \to \infty$, which completes the proof.
\end{proof}

\subsection{Gradient flows on Riemannian manifolds}
\label{sec:grad_flows}
Let $(\ML,g)$ be a Riemannian manifold without boundary 
of class $C^2$ over a Hilbert 
space $\HL$ with  Riemannian metric $g$ given fibrewise as 
$g_x(\cdot,\cdot)$ 
inducing a norm $\|\cdot\|_{g_x} \ceq \sqrt{g_x(\cdot,\cdot)}$ on 
the tangent space $T_x\ML$. 
For two points $x,y$ 
in the same connected component of $\ML$ we consider
the  set of connecting $C^1$-paths,
\begin{equation}\label{eq:path-space}
	\PS(x,y) 
	\ceq \braces[\big]{c\in C^1([0,1],\ML) : c(0)=x,c(1)=y},
\end{equation}
and define the \emph{path metric} (or \emph{Palais distance})
\begin{equation}\label{eq:path-metric}
\dist_{\ML}(x,y) \ceq \inf_{c\,\in\PS(x,y)}\int_0^1
\|\dot{c}(t)\|_{g_{c(t)}}\d t,
\end{equation}
where $\dot{c}(t) \ceq \D c_t(1)\in T_{c(t)}\ML$
for the differential $\D c_t:T_t[0,1]\simeq\R\to
T_{c(t)}\ML$ of the $C^1$-map $c:[0,1]\to\ML$ (cf.{} \cite[Sections 6 \& 9]{palais_1963}).
Equipped with the path metric $\dist_{\ML}$, every connected
component
of 
$\ML$ is a metric space. This metric is consistent with the
topology of $\ML$ as shown in \cite[Section 9]{palais_1963}
and in the more general setting of Finsler manifolds in
\cite[Theorem 3.3]{palais_1966}. Therefore, every connected
component $\KL$ of $\ML$ is closed with respect to $\dist_{\ML}$;
see, e.g., \cite[Remark after Theorem 25.2]{munkres_2000}.

On this manifold we consider energy functionals $\EL:\ML\to\R$ 
of class $C^{1,1}_\loc$ whose gradient $\nabla \EL$ is
given by $\D\EL_x(h)=g_x(\nabla\EL(x),h)$ for $x\in\ML$ and $h\in
T_x\ML$. We assume that
\begin{enumerate}
	\Item[\rm (e)]\label{lower-bound} 
	({\sc Energy bounded below})
	There is a number $e\in\R$ such that
	$\EL(x)\ge e$ for all $x\in\ML$.
	\Item[\rm (L**)]\label{condL**} 
	({\sc limit in $\ML$}) 
	Every Cauchy sequence
	$(x_k)_{k \in \N}$ in $\ML$ that is bounded under $\EL$
	converges to some 
	$x_\infty\in\ML$ as $k\to\infty$.
\end{enumerate}
Long-time existence for the gradient flow
and subconvergence to critical points can
now be proved by slightly generalizing a classic argument 
that can be found, e.g., in \cite[Theorem 9.1.6]{palais-terng_1988} 
for Riemannian manifolds, in \cite[Theorem 1.1]{browder_1965} 
for Banach manifolds modelled over a uniformly convex Banach space, or 
in \cite[Theorem 5.4]{palais_1966} for general Finsler manifolds.

\begin{theorem}[Gradient flow]
	\label[theorem]{thm:gradflow_subconvergence}
	Let $\EL\in C^{1,1}_\loc(\ML)$ satisfy the Palais--Smale condition
	{\rm (PS)} and conditions \ref{lower-bound} and \ref{condL**}.
	There exists a function $\phi:[0,\infty)\times\ML\to\ML$, such
	that $\phi(\cdot,x)\in C^1([0,\infty),\ML)$ for every $x\in\ML$,
	and
	\begin{equation}\label{eq:gradient-flow}
	\begin{cases}
	\frac{\d}{\d t}\phi(t,x)& =-\nabla\EL(\phi(t,x)) \Foa t\in [0,\infty),\\
		\phi(0,x)&=x.  
	\end{cases}
	\end{equation}
The function $t\mapsto\EL(\phi(t,x))$ is strictly decreasing 
on $[0,\infty)$ unless
$x$ is
	$\EL$-critical.
For $x\in\ML$ and $t\to\infty$ the solution $\phi(t,x)$ subconverges  
to a critical point of $\EL$, which lies in the same
connected component of $\ML$ as the initial point $x$.
In particular, in each connected component
of $\ML$ the energy $\EL$ possesses at least one
critical point. 
\end{theorem}

\begin{proof}
	For $x\in\ML$ we obtain
	the function $\phi$ as the solution to the 
	initial value problem
	\begin{equation} 
		\begin{aligned}
			\label{eq:cauchy_problem}
		\textstyle
		\frac{\d}{\d t}\xi(t) &=-\nabla \EL(\xi(t)), &\hspace{10pt}\xi(0)&=x.
		\end{aligned}
	\end{equation}
	Such a solution exists and is unique on an 
	interval $(t^-,t^{+})$ containing $t=0$
	by the Picard--Lindelöf Theorem 
	(see, e.g., \cite[Chapter IV.2]{lang_2002}). 
	Here $t^{+}\in (0,\infty]$ 
	denotes the maximal positive existence time, 
	and for $t\in (t^-,t^+)$ we compute by means of 
	\cref{eq:cauchy_problem}
	\begin{align}
	\textstyle
	\frac{\d}{\d t}\EL(\xi(t)) &
=\D\EL_{\xi(t)}(\dot{\xi}(t))= g_{\xi(t)}\big(\nabla \EL(\xi(t)),
	\dot{\xi}(t)\big)\label{eq:decreasing-energy}\\
		&\overset{\cref{eq:cauchy_problem}}{=}  
- g_{\xi(t)}\big(\nabla \EL(\xi(t)),\nabla \EL(\xi(t))\big)
=-\|\nabla\EL(\xi(t))\|^2_{g_{\xi(t)}}\le 0.\notag
	\end{align}
Hence, $\phi(\cdot,x) \ceq \xi(\cdot)$ 
solves 
\cref{eq:gradient-flow} on the positive maximal time interval $[0,t^+)$, 
and the energy is strictly decreasing along the flow unless $\D\EL_x=0$,
or equivalently, $\nabla\EL(x)=0$, in which case $\xi(t)=x$ and therefore
$\EL(\xi(t))=\EL(x)$ for all $t\in
[0,t^+)$ by uniqeness of the solution to \cref{eq:cauchy_problem}.

To show that $t^+=\infty$ we assume to the contrary that $t^+<\infty$,
and show below
that in that case the one-sided limit
$\xi(t^+) \ceq \lim_{t\nearrow t^+}\xi(t)$ exists, so that 
we could find a solution of \cref{eq:cauchy_problem} with $x$ replaced
by $\xi(t^+)$ on an
open
 time interval containing $t^+$, which contradicts the 
maximality of $t^+$ by uniqueness of the solution.

Towards the existence of the limit point
$\xi(t^+)\in\ML$ for some finite maximal time $t^+\in
	(0,\infty)$
consider 
a strictly increasing
sequence $(t_k)_k\subset[0,t^+)$ with $t_k\to t^+$ as $k\to\infty$,
and define 
the sequence of points
$x_k=\xi(t_k)\in\ML$. We first show that $x_k$ is a Cauchy 
sequence in $\ML$, so that by condition
\ref{condL**} we find $x_\infty\in\ML$ such that $x_k\to x_\infty$ as
$k\to\infty$, since we know  that
$\EL(x_k)\le\EL(x)$ for all $k\in\N$ by \cref{eq:decreasing-energy}. Going back
to an arbitrary sequence $s_m\nearrow t^+$ we will then show that also 
$y_m \ceq \xi(s_m)\to x_\infty$ as $m\to\infty$, so that $\xi(t^+)$ indeed exists.

We have by \cref{eq:decreasing-energy} 
\begin{align}
	\EL(\xi(t))-\EL(\xi(s))\overset{\cref{eq:decreasing-energy}}{=}
	\textstyle
	-\int_s^t \norm{ \nabla \EL(\xi(\tau))}^2_{g_{\xi(\tau)}}\d \tau
	\Fo 0\le s< t<t^+,\label{eq:energy-identity}
\end{align}
which yields the following by means of \ref{lower-bound}:
\begin{align*}
\textstyle
\int_{s}^{t} \norm{\nabla \EL(\xi(\tau))}_{g_{\xi(\tau)}}
\d \tau 
&\textstyle
	\le 
	\pars{t-s}^{1/2} 
	\pars[\big]{ \int_s^t
		\norm{\nabla \EL(\xi(\tau))}^2_{g_{\xi(\tau)}}\d \tau 
	}^{1/2}\\
&\overset{\ref{lower-bound}}{\le}(t-s)^{\frac{1}{2}}(\EL(x)-e)^{\frac{1}{2}}
\quad\Foa 0\le s<t<t^+.
	\end{align*}
	Hence, by definition \cref{eq:path-metric} of the path metric,
	and by the differential equation \cref{eq:cauchy_problem},
	\begin{align}
		\dist_\ML(x_j,x_k)&\overset{\cref{eq:path-metric}}{\le}
		\textstyle
		\int_{t_j}^{t_k}\|\dot{\xi}(\tau)\|_{g_{\xi(\tau)}}\d
		\tau \overset{\cref{eq:cauchy_problem}}{=}
		\int_{t_j}^{t_k}\|\nabla \EL(\xi(\tau)) \|_{g_{\xi(\tau)}}
		\d 
		\tau\notag\\
		&\le (t_k-t_j)^{\frac{1}{2}}(\EL(x)-e)
		\longrightarrow 0 \As j,k\to\infty;\label{eq:cauchy-estimate}
	\end{align}
hence $(x_k)_k$ is indeed a Cauchy sequence so that $x_k\to x_\infty$
for some $x_\infty\in\ML$. For any other sequence $s_m\nearrow t^+$ repeat
the argument in \cref{eq:cauchy-estimate} with $y_m \ceq \xi(s_m)$ instead of
$x_j$ to find
\begin{align*}
	\dist_\ML(y_m,x_\infty)
	=
	\lim_{k\to\infty} \dist_\ML(y_m,x_k)
	&\le
	\lim_{k\to\infty} \pars{t_k-s_m}^{1/2} \pars{\EL(x)-e}^{1/2}
	\\
	& = 
	\pars{t^+-s_m}^{1/2} \pars{\EL(x)-e}^{1/2}
	\to 
	0
	\As m\to\infty,
\end{align*}
so that the one-sided limit $\xi(t^+)$ indeed exists.
Therefore, 
	the solution $\phi(\cdot,x)$ exists on $[0,\infty)$.

Combining the energy identity \cref{eq:energy-identity} for 
$t \ceq t_k\nearrow\infty$ and 
$s \ceq 0$ with condition \ref{lower-bound}
we obtain
$$
\textstyle
\int_0^\infty\|\nabla\EL(\xi(\tau))\|^2_{g_{\xi(\tau)}}\d \tau=
\lim_{k\to\infty}\int_0^{t_k}\|\nabla\EL(\xi(\tau))\|^2_{g_{\xi(\tau)}}\d \tau
\le \EL(x)-e.
$$
Therefore, there exists a sequence $s_l\nearrow\infty$ such that
for the points $z_l \ceq \xi(s_l)\in\ML$ one has
$\|\nabla\EL(z_l)\|_{g_{z_l}}\to 0$ 
as $l\to\infty$. Together with 
$\EL(z_l)\le\EL(x)$ for all $l\in\N$ this shows that
$(z_l)_l\subset\ML$ is a (PS)-sequence, which by means of the (PS)-condition
implies
that there is a subsequence $(z_{l_i})_i\subset (z_l)_l$ 
converging to some point
$z_\infty\in\ML$ as $i\to\infty$. Since $\EL$ is of class $C^1$ we find
$$
\|\D\EL_{z_\infty}\|_{(T_{z_\infty}\ML)^*}=\lim_{i\to\infty}
\|\D\EL_{z_{l_i}}\|_{(T_{z_{l_i}}\ML)^*}=\lim_{i\to\infty}
\|\nabla\EL(z_{l_i})\|_{g_{z_{l_i}}}=0,
$$
that is, $z_\infty$ is $\EL$-critical. 
All points $z_{l_i}$ are connected to the initial point $x$ by
the continuous path $\xi|_{[0,s_{l_i}]}$, so they lie in the same
connected component $\KL$ of $\ML$ as $x$. Since $\KL$ is closed,
one also has $z_\infty\in\KL$.
\end{proof}
Palais considers in \cite[Section 2]{palais_1979} a group $G$ of
isometries of a Riemannian manifold and a $G$-invariant energy $\EL:\ML\to
\R$ to provide an exemplary simple setting where his
principle of symmetric criticality holds true. It turns out that in this
situation the gradient flow \cref{eq:cauchy_problem} is symmetry
preserving.
\begin{proposition}[Symmetry-preserving gradient flow]
\label{prop:symmetry-preserving-grad-flow}
Let $G$ be a group of $C^1$-isometries of the Riemannian manifold $(\ML,g)$. 
Suppose
that $\EL\in C^{1,1}_\textnormal{loc}(\ML)$ is invariant under $G$, i.e.,
$\EL=\EL\circ\G$ for all $\G\in G$, and satisfies the assumptions of
\Cref{thm:gradflow_subconvergence}. Then the gradient flow 
\cref{eq:cauchy_problem} preserves the $G$-symmetry, that is,
for every initial point $x\in\ML$ that satisfies $\G x=x$ for all $\G\in G$, the solution
$\xi$ of \cref{eq:cauchy_problem}  satisfies $\G \xi(t)=\xi(t)$ for all
$\G \in G$ and $t\in [0,\infty)$. In addition, $\xi$ subconverges to a 
$G$-symmetric critical point $x_\infty$, i.e., with $\G x_\infty = x_\infty$
for all $\G\in G$.
In particular, in any connected component $\KL$ of $\ML$ that contains
a $G$-symmetric point $x$ there is a $G$-symmetric $\EL$-critical point.
\end{proposition}
\begin{proof}
Palais showed in \cite[Section 2]{palais_1979} that 
$\nabla\EL(\G x)=
\D\G_x(\nabla\EL(x))$ for all $\G\in G$ and $x\in\ML$. 
This implies for any $G$-symmetric
initial point $x\in\ML$ and any fixed group element $\G\in G$ that
the solution $\xi$ of \cref{eq:cauchy_problem}
generates the
function $\varXi(t) \ceq \G\xi(t)$, which satisfies $\varXi(0)=\G\xi(0)=
\G x=x$ and
$$
\dot{\varXi}(t)=\D\G_{\xi(t)}(\dot\xi(t))\overset{\cref{eq:cauchy_problem}}{=}
-\D\G_{\xi(t)}\nabla\EL(\xi(t))=-\nabla\EL\G\xi(t)=-\nabla\EL\varXi(t)
$$
for all $t\in [0,\infty)$. 
In other words, $\varXi$ solves the same initial
value problem as $\xi$ does, so that $\varXi=\xi$ by uniqueness of the solution.
Since $\G\in G$ was chosen arbitrarily, we have shown that \cref{eq:cauchy_problem} is 
symmetry-preserving. If $x_\infty\in\ML$ is an $\EL$-critical limit point of 
$\xi$ then there is a sequence 
$t_i\to\infty$ such that 
$$
\G x_\infty=\G \lim_{i\to\infty}\xi(t_i)=\lim_{i\to\infty}\G\xi(t_i)=
\lim_{i\to\infty}\xi(t_i)=x_\infty,
$$
so $x_\infty $ is $G$-symmetric as well. Finally, to find a $G$-symmetric
$\EL$-critical point in a component $\KL$ of $\ML$ containing a $G$-symmetric
point $x\in\KL$, we simply start the gradient flow \cref{eq:cauchy_problem}
in that point $x$. This flow subconverges to a 
$G$-symmetric $\EL$-critical point
in that same component $\KL$.
\end{proof}
\begin{corollary}[Critical points of lower energy]
{\cite[Corollary 9.1.8]{palais-terng_1988}}
	\label[corollary]{koro:critpoint_lowerenergy}
Let $\EL$ satisfy the assumptions of \Cref{thm:gradflow_subconvergence}.
Then, for any point $x\in\ML$ with $\D\EL_x\not=0$ there
exists a critical point $x_{\text{c}}$ of $\EL$ in the connected component
of $\ML$ that contains $x$, such that $\EL(x_{\text{c}})<\EL(x)$.
If $x$ is $G$-symmetric and $\EL$ is $G$-invariant under a group
$G$ of $C^1$-isometries of $(\ML,g)$ then $x_c$ is $G$-symmetric as well.
\end{corollary}
\begin{proof}
	Choose $x_{\text{c}}\in\ML$ as any limit point of $\phi(t,x)$
	as $t\to\infty$, 
	which is contained in the same connected
	component $\KL\subset\ML$ as the initial point $x$ according
	to \Cref{thm:gradflow_subconvergence}. 
	Since $x$ is not a critical point we have $\EL(x_\text{c})<\EL(x)$
	because of the strict energy decay along the flow. The last
	statement follows from \Cref{prop:symmetry-preserving-grad-flow}.
\end{proof}
\begin{corollary}[Energy minimizers]{\cite[Theorem 9.1.9]{palais-terng_1988}}
	\label[corollary]{koro:minimizer_existence}
	Under the assumptions of \Cref{thm:gradflow_subconvergence}
	 the energy $\EL$ attains its infimum on $\ML$, and additionally on 
 every connected component of $\ML$. Moreover, on $\ML$, and also
 in each of its connected components
 there exists a minimizing {\rm (PS)}-sequence. 
 If a connected component
 $\KL$ of $\ML$ contains a $G$-symmetric point $x$ and $\EL$ is $G$-invariant
 under a group $G$ of $C^1$-isometries 
 of $(\ML,g)$, then there is a $G$-symmetric
 element $x_s\in\KL$ with $\EL(x_s)=\inf\{\EL(y):y\in\KL,\,\textnormal{$y$
 is $G$-symmetric}\}$, and a corresponding {\rm (PS)}-sequence of
 $G$-symmetric elements in $\KL$.
\end{corollary}
\begin{proof}
Let $x_k$ be a minimizing sequence for $\EL$ on $\ML$, i.e.,
$\lim_{k\to\infty}\EL(x_k)=\inf(\EL)$. 
Without loss of generality we can assume, that each $x_k$ is a critical point, since otherwise we just replace $x_k$ 
by a critical point of lower energy by \cref{koro:critpoint_lowerenergy}. Then, since $\EL$ satisfies the (PS)-condition, 
there exists a subsequence $(x_{k_i})_i\subset (x_k)_k$ 
converging to a critical point 
$x_{\text{c}}$ as $i\to\infty$, and by continuity of $\EL$ one has
$
\EL(x_{\text{c}})=\lim_{i\to\infty}\EL(x_{k_i})=\inf(\EL).
$
We can repeat the same argument on each connected component of $\ML$ instead
of $\ML$, since each such component is closed.
Finally, if we consider a minimizing sequence of $G$-symmetric
elements $x_k\in\KL$, then the same reasoning is valid due to the last
statement of \cref{koro:critpoint_lowerenergy} and because of the fact
that $G$-symmetry is a closed condition since each $\G\in G$ is a
continuous mapping from $\ML$ onto $\ML$.
\end{proof}

Even if the energy is only in $C^1$, we can use Ekeland's variational
principle to obtain energy-minimizing
Palais--Smale sequences subconverging to
energy minimizers. 

\begin{corollary}[Minimizers]
\label{cor:critical+minimizers}
Let $\EL\in C^1(\ML)$ satisfy the Palais--Smale condition {\rm (PS)} and conditions \ref{lower-bound} and \ref{condL**}.
Let $\KL$ be a connected component of $\ML$.
Then there are a minimizing {\rm (PS)}-sequence and a minimizer of $\EL |_{\KL}$.

If  $\KL$ contains a $G$-symmetric element and if
 $\EL$ is $G$-invariant, then
there exists a (not necessarily $G$-symmetric) {\rm (PS)}-sequence for $\EL|_\KL$ that converges to
a $G$-symmetric element in $\KL$ that realizes $\inf\{\EL(y):y\in\KL,\,\textnormal{$y$
is $G$-symmetric}\}$.
These statements hold true also for $\KL = \ML$.
\end{corollary}
\begin{proof}
We  choose a minimizing sequence $(x_k)_{k \in \N}$ of $\EL |_{\KL}$ such that $\EL(x_k) \le \inf(\EL|_{\KL}) + 1/k^2$.
The function $\EL$ is continuous and the set $\KL$ is closed.
Hence, condition \ref{condL**} implies that the sublevel set 
$
	\mathcal{U} \ceq \braces[\big]{x \in \KL : \EL(x) \leq \inf (\EL|_\KL) + 2}
$ 
equipped with the restriction of $\dist_\ML$ is a complete metric space.
Now we apply Ekeland's variational principle (see, e.g., \cite[Chapter~1,Theorem~5.1]{struwe_2008a}) to find $y_k \in \mathcal{U}$
such that 
$
	\EL(y_k) \le \EL(x_k)
$,
$
	\dist_\ML(x_k,y_k) \le 1/k
$,
and
\begin{equation}
	\label{eq:ekeland-ineq}
	\EL(y_k) \leq \EL(z) + \tfrac{1}{k} \dist_\ML(y_k,z) \quad \Foa z\in\mathcal{U}
	.
\end{equation}
We have $\EL(y_k) \leq \EL(x_k) \leq \inf(\EL|_{\KL}) + 1 < \inf(\EL|_\KL) + 2$.
Because $\EL$ is continuous and because $\KL$ is open, $ y_k $ is contained in the interior of $ \UL $.
So, for every tangent vector $h\in T_{y_k}\ML$ 
there is a $\sigma > 0$ and a $C^1$-path $c \colon \intervaloo{-\sigma,\sigma} \to \mathcal{U}$ satisfying $c(0)=y_k$ and $\dot c(0) = h$.
By continuity, we can take $\sigma$ so small that
\begin{equation}\label{eq:tangent-small}
	\textstyle
	\norm{\dot c(t)}_{g_{c(t)}}\le 2\norm{h}_{g_{c(0)}} = 2 \norm{h}_{T_{y_k}\ML}
\quad\Foa t\in (-\sigma,\sigma).
\end{equation}
Applying \cref{eq:ekeland-ineq} to the point
$z \ceq c(t)$ for any $t\in (0,\sigma)$
we obtain
\begin{equation}\label{eq:PS-weg}
\textstyle
	\EL(y_k)-\EL(c(t))
	\overset{\cref{eq:ekeland-ineq}}{\le}
	\frac1{k}\dist_\ML(y_k,c(t)) \le \frac1{k}\int_0^t
	\norm{c'(\tau)}_{g_{c(\tau)}}\d \tau\le \frac2{k}\norm{h}_{T_{y_k}\ML}\cdot t
	.
\end{equation}
Dividing by $t>0$ and taking the limit $t\searrow 0$,
we arrive at
$\D \EL_{y_k} h \le \frac{2}{k} \norm{h}_{T_{y_k}\ML}$, hence
\[ 
	\norm{\D\EL_{y_k}}_{(T_{y_k}\ML)^*} \le \tfrac{2}{k}\to 0
	\quad  \text{as} \quad  k \to\infty
	.
\] 
Thus, $(y_k)_k$ is a minimizing (PS)-sequence for $\EL|_{\KL}$.
In particular, it is a (PS)-sequence for $\EL$.
Because $\EL$ satisfies the (PS)-condition, the sequence subconverges to some $y_\infty \in \ML$: $y_{k_i}\to y_\infty$ as $i\to\infty$.
Since $\mathcal{U}$ is closed, we have $y_\infty \in \mathcal{U} \subset \KL$, hence
\[
	\inf(\EL|_{\KL})
	\leq 
	\EL(y_\infty)
	=
	\lim_{i\to\infty} \EL(y_{k_i})
	\le
	\lim_{i\to\infty} \EL(x_{k_i})
	=
	\inf(\EL|_{\KL})
	. 
\]
Thus, $y_\infty$ is the desired minimizer.

If $\EL$ is $G$-invariant,
we can start with a minimizing sequence $(x_k)_k\subset\KL$ of $G$-symmetric
elements, and then follow the same line of arguments to show
that a subsequence $(x_{k_i})_i\subset (x_k)_k$ 
converges to some $y_\infty\in \KL$ as $i\to\infty$. Because $G$-symmetry is a closed condition, we see that the limit point $y_\infty$ is $G$-symmetric, too. Finally, we notice that the above argument goes through, too, if we replace $\KL$ by $\ML$.
\end{proof}

%% file: ManifoldA_new.tex
 
\section{The submanifolds of arclength parametrized curves and knots}\label{sec:arclength}
Throughout this section we fix any parameter $s\in (\frac{1}{2},1]$ and
consider the (for $s\in (\frac{1}{2},1)$ fractional) Sobolev space 
\begin{equation}\label{eq:frac-space}
H^{1+s}(\R/\Z,\R^n) \ceq W^{1+s,2}(\R/\Z,\R^n),
\end{equation}
which embeds into the classic H\"older space
$C^{1,s-1/2}(\R/\Z,\R^n)$ and 
hence compactly into $C^1(\R/\Z,\R^n)$; 
see \cite[Theorem A.2]{knappmann-etal_2023}
for a short self-contained proof in the one-dimensional periodic setting
for $s\in (\frac{1}{2},1).$
Therefore, 
the subsets $H^{1+s}_\reg(\R/\Z,\R^n)$ or $H^{1+s}_\inj(\R/\Z,\R^n)$ 
of regular, respectively injective
curves are open subsets of the Hilbert space
$H^{1+s}(\R/\Z,\R^n)$, as well as their
intersection $H^{1+s}_\injreg(\R/\Z,\R^n)$.

%Using a hack to suppress warnings that PDF bookmarks don't like TeX code.
\subsection{The arclength manifolds \texorpdfstring{$\NL^{1+s}$}{N1+s} and \texorpdfstring{$\AL^{1+s}$}{A1+s}}
\label{sec:AL}

To identify the subset 
\begin{equation}\label{eq:arclenth-set}
H^{1+s}_\unitspeed(\R/\Z,\R^n)=\{\g\in H^{1+s}(\R/\Z,\R^n):|\g'|=1\,
\textnormal{
on $\R/\Z$}\}
\end{equation}
of arclength parametrized curves as a submanifold of the open set
$H^{1+s}_\reg(\R/\Z,\R^n)$ we follow the idea of 
\cite[Section 2]{scholtes-etal_2022} and
\cite{knappmann-etal_2023} and introduce
the scalar-valued \emph{logarithmic strain}
$\varSigma:H^{1+s}_\reg(\R/\Z,\R^n)\to H^s(\R/\Z)$ defined as
\begin{equation}\label{eq:logarithmic-strain}
	\varSigma(\g) \ceq \log\big( |\g'|\big)  \in H^s(\R/\Z) \quad\Fo \g\in H^{1+s}_\reg(\R/\Z,\R^n),
\end{equation}
which is a smooth function on the set $H^{1+s}_\reg(\R/\Z,\R^n)$ of
regular curves. It follows from \cite[Proposition 5.1]{knappmann-etal_2023}
for $s\in (\frac{1}{2},1)$
that the differential $\D\varSigma_\g:T_\g H^{1+s}_\reg(\R/\Z,\R^n)
\simeq H^{1+s}(\R/\Z,\R^n)\to T_{\varSigma(\g)}H^s(\R/\Z,\R^n)\simeq
H^s(\R/\Z,\R^n)$ has the form
\begin{equation}\label{eq:log-strain-differential}
	\D\varSigma_\g h=\langle\g',h'\rangle \quad\Fo
	\g\in H^{1+s}_\unitspeed(\R/\Z,\R^n),\, h\in H^{1+s}(\R/\Z,\R^n),
\end{equation}
and is bounded, i.e. has operator norm
\begin{equation}\label{eq:log-strain-bounded-differential}
\|\D\varSigma_\g\|_{\LL(H^{1+s},H^s)}\le C_\varSigma,
\end{equation}
where the constant  $C_\varSigma=C_\varSigma(n,s,\|\g'\|_{H^s})$
depends 
non-decreasingly on the norm $\|\g'\|_{H^s}$
of the unit tangent $\g'$. In view of the explicit
formula \cref{eq:log-strain-differential}
for $\D\varSigma_\g$ it can easily be checked that
the same holds true for $s=1$.
Moreover, \cite[Lemma 5.2]{knappmann-etal_2023}
yields for $s\in (\frac{1}{2},1)$
the existence of a bounded right inverse $Y_\g\in
\LL(H^s(\R/\Z),H^{1+s}(\R/\Z,\R^n))$ of the differential $\D\varSigma_\g$ 
satisfying
$\D\varSigma_\g(Y_\g(f))=f$ for all $f\in H^s(\R/\Z)$,
so that $\D\varSigma_\g$ is surjective for any
$\g\in H^{1+s}_\unitspeed(\R/\Z,\R^n)$. Again, as before, 
the same formula\footnote{Explicitly given in \cref{eq:right-inverse-Y}
in our proof
of \cref{prop:projection-ALs} below.}
for $Y_\g$ given in the proof of \cite[Lemma 5.2]{knappmann-etal_2023}
yields a right inverse of $\D\varSigma_\g$
also in the case $s=1$.
Therefore, we can apply the
preimage theorem\footnote{Notice that the nullspace $\ker(\D\varSigma_\g)$
splits $H^{1+s}(\R/\Z,\R^n)$ because the latter is a Hilbert space.}
(see, e.g., \cite[Theorem 4.J]{zeidler_1993}) 
to
\begin{equation}
\begin{aligned}
	\label{eq:arclength-manifold}
&\NL^{1+s}\ceq H^{1+s}_\unitspeed(\R/\Z,\R^n)=\big(\varSigma|_{H^{1+s}_\reg(\R/\Z,\R^n)}
\big)^{-1}(0),\hspace{5pt}\text{ and }\\
&\AL^{1+s} \ceq H^{1+s}_\injunitspeed(\R/\Z,\R^n)=\big(\varSigma|_{H^{1+s}_\injreg(\R/\Z,\R^n)}\big)^{-1}(0)=\NL^{1+s}\cap H^{1+s}_\inj\left( \R/\Z,\R^n \right).
\end{aligned}
\end{equation}

\begin{theorem}[Arclength manifolds]
\label{thm:arclength-manifold}
Let $s\in (\frac{1}{2},1]$.
The sets $\NL^{1+s}$ of arclength para\-met\-rized curves, and $\AL^{1+s}$ of arclength para\-met\-rized knots are smooth submanifolds
of the open subsets $H^{1+s}_\reg(\R/\Z,\R^n)$ and $ H^{1+s}_\injreg(\R/\Z,\R^n) $, respectively. 
\end{theorem}
To secure compactness of (PS)-sequences  later,
we can fix an arbitrary 
point $\bz\in\R^n$ and
consider the subsets
\begin{equation}\label{eq:arclenth-point}
\NL^{1+s}_\bz\ceq\{\g\in\NL^{1+s}:\g(0)=\bz\},\, \text{ and }\,\AL^{1+s}_\bz\ceq\{\g\in\AL^{1+s}:\g(0)=\bz\}.
\end{equation}
The corresponding ambient space $H^{1+s}_\bz(\R/\Z,\R^n)\ceq\{\g\in
H^{1+s}(\R/\Z,\R^n):\g(0)=\bz\}$ is an affine subspace
and therefore also a smooth submanifold of $H^{1+s}(\R/\Z,\R^n)$, but in order
to come up with a suitable projection later, we express $\NL^{1+s}_\bz$ as
preimage of $(0,\bz)\in H^s(\R/\Z)\times\R^n$, i.e., 
\begin{align*}\label{eq:NL_z-preimage}
	\NL^{1+s}_\bz=\big(
\varSigma|_{H^{1+s}_\reg(\R/\Z,\R^n)}\times P_0\big)^{-1}(0,\bz),
\end{align*}
where the point evaluation $P_0:H^{1+s}_\reg(\R/\Z,\R^n)\to\R^n$ is 
defined by
\begin{equation}\label{eq:point-evaluation}
P_0(\g)\ceq\g(0).
\end{equation}
This linear map is a smooth submersion at every $\g\in H^{1+s}_\reg(\R/\Z,\R^n)$
with differential $\D(P_0)_\g:T_\g H^{1+s}_\reg(\R/\Z,\R^n)\simeq
H^{1+s}(\R/\Z,\R^n)\to T_{P_0(\g)}\R ^n\simeq \R^n$  given by
\begin{equation}\label{eq:P0-differential}
\D(P_0)_\g h=P_0 h=h(0)\quad\Foa h\in H^{1+s}(\R/\Z,\R^n).
\end{equation}
The differential of the product map $\varSigma\times P_0:
H^{1+s}_\reg(\R/\Z,\R^n)\to H^s(\R/\Z)\times\R^n$ is
$$
\D(\varSigma,P_0)_\g=(\D\varSigma_\g,\D(P_0)_\g)
\overset{\cref{eq:P0-differential}}{=}(\D\varSigma_\g,P_0),
$$
 and it possesses as bounded right
inverse the map $Z_\g:H^s(\R/\Z)\times\R^n\to H^{1+s}(\R/\Z,\R^n)$
defined via the right inverse $Y_\g$ of $\D\varSigma_\g$ (see the arguments that lead to
\cref{thm:arclength-manifold}) as
\begin{align*}
	Z_\g(f,\by) \ceq Y_\g(f)+\by-Y_\g(f)(0)\quad\Fo (f,\by)\in
H^s(\R/\Z)\times\R^n.
\end{align*}
Indeed, one computes with the help of \cref{eq:P0-differential}, 
\cref{eq:point-evaluation} and
\cref{eq:log-strain-differential}
\begin{align*}
	\D(\varSigma,P_0)_\g & \big(Z_\g(f,\by)\big)
	\overset{\cref{eq:P0-differential}}{=}
	\big(\D\varSigma_\g
	\big(Y_\g(f)+\by-Y_\g(f)(0)\big),P_0\big(Y_\g(f)+\by-
	Y_\g(f)(0)\big)\big)\\
	&\overset{\cref{eq:log-strain-differential},\cref{eq:point-evaluation}}{=}
	\big(\D\varSigma_\g(Y_\g(f)),\by\big)=(f,\by)\quad\Foa (f,\by)
	\in H^s(\R/\Z)\times\R^n.
\end{align*}
We have thus shown that the differential $\D(\varSigma,P_0)_\g$ is surjective
for any curve $\g\in H^{1+s}_\reg(\R/\Z,\R^n)$, so that we may apply the
preimage theorem as in the proof of \cref{thm:arclength-manifold}.
Including injectivity, we get 
\begin{equation}\label{eq:AL_z-preimage}
	\AL^{1+s}_\bz=\big(
	\varSigma|_{H^{1+s}_\injreg(\R/\Z,\R^n)}\times P_0\big)^{-1}(0,\bz)=\NL^{1+s}_{\bz}\cap H^{1+s}_{\inj}.
\end{equation}

\begin{theorem}[Arclength manifolds with fixed initial point]
\label{thm:arclength-manifold-point}
Let $s\in (\frac{1}{2},1]$ and $\bz\in\R^n$, then the sets $\NL^{1+s}_\bz$ and 
$\AL^{1+s}_\bz$
of arclength parametrized curves and knots $\g$ with fixed initial point $\g(0)=\bz$ are
smooth submanifolds of the open subsets $H^{1+s}_\reg(\R/\Z,\R^n)$ and $ H^{1+s}_{\injreg}(\R/\Z,\R^n) $, respectively. 
\end{theorem}

The tangent space of the arclength manifold $\NL^{1+s}$ at a curve $\g\in
\NL^{1+s}$
is given by
\begin{equation}\label{eq:tangent-space-AL}
	T_\g\NL^{1+s} 
	=
	\ker (\D\varSigma_\g)
	=
	\braces[\big]{ 
		h\in H^{1+s}(\R/\Z,\R^n) : \inner{\g',h'}=0
	},
\end{equation}
whereas for $\NL_\bz^{1+s}$ we have
\begin{equation}\label{eq:tangent-space-ALz}
	T_\g\NL^{1+s}_\bz 
	=
	\ker(\D(\varSigma,P_0)_\g)
	=
	\braces[\big]{ 
		h\in H^{1+s}(\R/\Z,\R^n) : \inner{\g',h'}=0,\,h(0)=0
	}.
\end{equation}
Analogously, one finds $ T_{\g}\AL^{1+s}=T_\g\NL^{1+s} $ for $ \gamma\in \AL^{1+s} $ and $ T_\g\AL^{1+s}_{\bz}=T_\g\NL^{1+s}_{\bz} $ if $ \gamma\in\AL^{1+s}_{\bz} $. 
\begin{proposition}[Bounded linear projection onto the tangent space]
\label{prop:projection-ALs}
Let $s\in (\frac{1}{2},1]$, $ \ML\in \left\{\AL^{1+s},\NL^{1+s}\right\}$, and $\g\in\ML$. Then the mapping
$\pr_\g^{\ML}:H^{1+s}(\R/\Z,\R^n)\to H^{1+s}(\R/\Z,\R^n)$ given by
\begin{equation}\label{eq:projection-ALs}
\pr_\g^{\ML}h\ceq(\Id-Y_\g \D\varSigma_\g)h
\end{equation}
is a linear and bounded projection onto the tangent space $T_\g\ML$
of $\ML$ at $\g$.
\end{proposition}
\begin{proof}
We need the explicit formula for 
the right inverse $Y_\g$ of $\D\varSigma_\g$
as provided in \cite[Lemma 5.2]{knappmann-etal_2023}:
\begin{equation}\label{eq:right-inverse-Y}
	\textstyle
	Y_\g f(t)
	=
	\int_0^t \big( f(u) \, \g'(u) + \P_{\g'(u)}^\perp V_{\gamma}f\big) \d u,
	\;\;
	\text{with}
	\;\;
	V_{\gamma}f=-\varTheta_\g^{-1}\int_0^1 f(u) \, \g'(u)\d u\in\R^n
	,
\end{equation}
where $\varTheta_\g:\R^n\to\R^n$ is 
the bounded invertible linear operator
given by
$\varTheta_\g\ceq\int_0^1\P_{\g'(u)}^\perp \d u$. Note that $ \varTheta_{\gamma}^{-1} $ is indeed bounded because of \cite[Lemma 5.3]{knappmann-etal_2023}, since $ \gamma'\in C^{0,s-\frac{1}{2}}\left( \R/\Z,\R^n \right) $ by Morrey's embedding. 
The vector $V_{\gamma}f$ guarantees that $Y_\g(f)$ is $1$-periodic. 
It depends linearly  on $f$ and satisfies the estimate
\begin{equation}\label{eq:V-h-est}
 |V_{\gamma}f|\le C\|f\|_{C^0}\le C_EC\|f\|_{H^s}\quad\Foa f\in H^s(\R/\Z),
\end{equation}
where the constant $C=C(s,\|\g'\|_{H^s})$ depends non-de\-creas\-ingly on
the norm $\|\g'\|_{H^s}$ and $C_E=C_E(s)$ is the embedding constant from
the compact Morrey embedding $H^s(\R/\Z)\hookrightarrow
 C^0(\R/\Z)$, which is standard for $s=1$; for
 $s\in (\frac{1}{2},1)$ 
 see, e.g., \cite[Theorem A.2]{knappmann-etal_2023}. 
So, we compute
\begin{equation}\label{eq:derivative-right-inverse}
	\pars{Y_\g \, f}'(t) = f(t) \, \g'(t)+\P_{\g'(t)}^\perp V_{\gamma}f;
\end{equation}
$
	\inner{ \pars{Y_\g f}'(t),\g'(t) }=f(t)$ for all $t\in\R/\Z$ and 
$f\in H^s(\R/\Z),
$
and therefore, by \cref{eq:log-strain-differential}
\begin{align}
	\inner{ (\pr_\g^{\ML}h)'(t),\g'(t) }&
	\overset{\cref{eq:log-strain-differential}}{=}
	\inner[\big]{ 
		h'(t) - \pars[\big]{Y_\g \inner{\g', h'}}'(t)
		,
		\g'(t)
	}=0\label{eq:tangent-orthogonal}
\end{align}
for all $h\in H^{1+s}(\R/\Z,\R^n).$
This 
shows that $\pr_\g^{\ML}h\in T_\g\ML$ for all $h\in H^{1+s}(\R/\Z,\R^n).$
Obviously, $\pr_\g^{\ML}$ is linear, so that it is a projection as desired if
we can show that its restriction to $T_\g\ML$ equals the identity (see, e.g., \cite[9.13(1)]{alt_2016}). But linearity of $f\mapsto Y_\g f$ 
indeed yields for $h\in
T_\g\ML$  by means of \cref{eq:tangent-space-AL}
\begin{equation}\label{eq:proj-id}
	\pr_\g^{\ML}h
	=
	h-Y_\g \inner{ \g',h'} 
	\overset{\cref{eq:tangent-space-AL}}{=}h-Y_\g(0)
	=
	h.
\end{equation}
Combining the explicit formulas \cref{eq:right-inverse-Y}
and \cref{eq:derivative-right-inverse} with the estimates
\cref{eq:V-h-est}, \cref{eq:log-strain-bounded-differential} 
and the product rule 
(see, e.g., \cite[Proposition A.5]{knappmann-etal_2023}),     we obtain
\begin{equation}\label{eq:projection-bound}
\|\pr_\g^{\ML}h\|_{H^{1+s}}\le (1+\tilde C C_\varSigma)\|h\|_{H^{1+s}}=:C_\g\|h\|_{H^{1+s}}\quad\Foa 
h\in H^{1+s}(\R/\Z,\R^n),
\end{equation}
where the constant $\tilde C=\tilde C(n,s,\|\g\|_{H^{1+s}})$ and therefore also $C_\g$
depends non-de\-creas\-ingly on
the norm $\|\g\|_{H^{1+s}}.$
\end{proof}
\begin{corollary}[Bounded linear projection onto $T_\g\NL^{1+s}_\bz$]
\label{cor:projection-ALsz}
For $s\in (\frac{1}{2},1]$, $\bz\in\R^n$, $ \ML\in \left\{\NL^{1+s},\AL^{1+s}\right\} $, and $\g\in\ML_\bz\ceq\{\eta\in \ML: \eta(0)=\bz\}$, the
mapping $\pr_\g^{\ML_\bz}:H^{1+s}(\R/\Z,\R^n)\to 
H^{1+s}(\R/\Z,\R^n)$
defined by
\begin{equation}\label{eq:projection-ALsz}
\pr_\g^{\ML_\bz}h\ceq\pr_\g^{\ML}h-
\pr_\g^{\ML}h(0)
\end{equation}
is a linear and bounded projection onto the tangent space
$T_\g\ML_\bz$ of $\ML_\bz$ at $\g$. 
Here, $\pr_\g^{\ML}$ is the projection provided by \cref{prop:projection-ALs}.
\end{corollary}
\begin{proof}
Notice that the
linearity and boundedness of $\pr_\g^{\ML}$ 
transfers to $\pr_\g^{\ML_\bz}$, and with 
$(\pr_\g^{\ML_\bz}h)'
=
(\pr_\g^{\ML}h)'$ and $\pr_\g^{\ML_\bz}h(0)=0$ we find
by \cref{eq:tangent-orthogonal} and
\cref{eq:tangent-space-ALz} that $\pr_\g^{\ML_\bz}h
\in T_\g\ML_\bz$. Finally, for $h\in T_\g\ML_\bz
\subset T_\g\ML$ we have $h(0)=0$ and therefore, due to \cref{eq:proj-id}, 
\begin{align*}
	\pr_\g^{\ML_\bz}h=
	\pr_\g^{\ML}h\overset{\cref{eq:proj-id}}{=}h.
\end{align*}
\end{proof}

Now \cref{thm:arclength-manifold,thm:arclength-manifold-point},
\cref{prop:projection-ALs}, and \cref{cor:projection-ALsz} 
together with the estimate \cref{eq:projection-bound}
verify the basic abstract conditions
\ref{condE} and \ref{condP} for the arclength manifolds $\AL^{1+s}$ 
and
$\AL^{1+s}_\bz$.
\begin{corollary}[\ref{condE} and \ref{condP} for arclength manifolds]
\label{cor:basic-assumptions}
Let $s\in (\frac{1}{2},1]$ and $\bz\in\R^n$ be fixed. Then condition \ref{condE}
holds true for $\ML\in \left\{\AL^{1+s},\AL_{\bz}^{1+s},\NL^{1+s},\NL_{\bz}^{1+s}\right\}$, if we choose 
$\OL \ceq H^{1+s}_\reg(\R/\Z,\R^n)$, or $ \OL\ceq H^{1+s}_\injreg(\R/\Z,\R^n)$ respectively, as an open subset of the Hilbert space
$\HL \ceq H^{1+s}(\R/\Z,\R^n)$, which embeds compactly into
the Banach space $\BL \ceq C^1(\R/\Z,\R^n)$. 
Moreover, the linear projections $\pr_\g^{\ML}:
T_\g\OL\to T_\g\ML$ for $\g\in\ML$ have property \ref{condP}.
\end{corollary}

\begin{remark}[Explicit Lagrange multiplier]\label{rem:explicit-lagrange-multiplier}
	The concrete formulas for the right inverse $ Y_{\gamma} $ in \cref{eq:right-inverse-Y} and for the projection $ \pr^{\ML}_{\gamma} $ in \cref{eq:projection-ALs} for $ \ML\in \{\AL^{1+s},\NL^{1+s}\} $ allow us to derive an explicit expression for the Lagrange multiplier arising from the arclength constraint. This turns out to be extremely useful in the regularity theory in \cref{sec:Smoothness}. 

	Suppose $ \gamma\in \ML $ is critical for the restricted energy $ \widetilde{\EL}\ceq\EL\vert_{\ML} $, where $ \EL \in C^1(\OL)$ for an open subset of the Hilbert space $ \HL=H^{1+s}\left( \R/\Z,\R^n \right) $ with $ \ML\subset\OL $. Then $ \D \EL_{\gamma}h=\D \widetilde{\EL}_{\gamma}h =0$ for all $h\in \T_{\gamma}\ML $, which is equivalent to 
	\begin{equation}\label{eq:Lagra_Mult}
	\begin{aligned}
		0=\D\EL_{\gamma} \pr^{\ML}_{\gamma}\overset{\cref{eq:projection-ALs}}{=} \D\EL_{\gamma}\left( \Id_{\HL}-Y_{\gamma}\D\Sigma_{\gamma} \right).
	\end{aligned}
	\end{equation}
	By the Lagrange multiplier theorem, see \cite[Corollary 3.5.29]{abraham-etal_1988}, there exists $ \Lambda\in \left( H^{s}(\R/\Z) \right)^* $ such that $ \D\EL_{\gamma}-\Lambda\D\Sigma_{\gamma}=0 $. This in combination with \cref{eq:Lagra_Mult} yields the formula 
	\[\Lambda=\D\EL_{\gamma}Y_{\gamma}.\]
\end{remark}

\subsection{The topology of arclength manifolds}
\label{sec:AL-toplogy}
Fix a point $\bz\in\R^n$ and $s\in \intervaloc{0,1}$. The Hilbert space $ \HL=H^{1+s}\left( \R/\Z,\R^n \right) $ with its inner product\footnote{For $ s=1 $ simply replace the double integral in \cref{eq:fractionalscalarproduct} by $ \left\langle \gamma'',h'' \right\rangle_{L^2}  $.}
\begin{equation}\label{eq:fractionalscalarproduct}
\begin{aligned}
	\inner{f,h}_\HL&\ceq\inner{f,h}_{H^{1+s}}=\inner{f,g}_{L^2}+\inner{f',g'}_{L^2}\\
	&+\int_{\R/\Z}\int_{\R/\Z}\frac{\inner{f(u)-f(v),h(u)-h(v)}}{
	\abs{u-v}_{\R/\Z}^{1+2s}}\d u \d v,\quad 
	\text{for $f,h \in \HL$},
\end{aligned}
\end{equation}

induces a Riemannian metric $g$ on  the arclength submanifolds
\begin{align*}
	&\NL^{1+s}_\bz\subset\NL^{1+s}\subset H^{1+s}_\reg(\R/\Z,\R^n)\subset
H^{1+s}(\R/\Z,\R^n)\,\text{ and }\\
&\AL^{1+s}_\bz\subset\AL^{1+s}\subset H^{1+s}_\injreg(\R/\Z,\R^n)\subset
H^{1+s}(\R/\Z,\R^n).
\end{align*}
More precisely, let $ \ML\in \left\{\NL^{1+s},\NL^{1+s}_{\bz},\AL^{1+s},\AL^{1+s}_{\bz}\right\} $,
$g_\g(f,h)\ceq\langle f,h\rangle_\HL$
for $\g\in\ML$ and $f,h\in T_\g\ML$. So, $\ML$ is a
smooth Riemannian submanifold of the open set $H^{1+s}_\reg(\R/\Z,\R^n) $, or $H^{1+s}_\injreg(\R/\Z,\R^n)$,
respectively; see
\cite[Theorem 3.6]{palais_1966}.
As in the abstract setting of \cref{sec:grad_flows} we have
the path metric (cf.{} \cref{eq:path-metric})
\begin{equation}\label{eq:path-metricALsz}
\dist_{\ML}(\g,\eta)\ceq\inf_{c\in \PS(\g,\eta)}\int_0^1 
\sqrt{g_{c(t)}(\dot{c}(t),\dot{c}(t))}\d t=
\inf_{c\in \PS(\g,\eta)}\int_0^1 \|\dot{c}(t)\|_\HL \d t 
\end{equation}
measuring the distance of
two curves $\g,\eta\in\ML$
in the same connected component of $\ML$, where 
\begin{equation}\label{eq:path-spaceAL}
	\PS(\g,\eta)
	\ceq
	\braces[\big]{
		c\in C^1([0,1],\ML)
		:
		c(0)=\g,c(1)=\eta
	}.
\end{equation}

To investigate the topology of $\ML$ 
more closely,
we look at a mapping that transforms an arbitrary regular closed curve $\g
\in H^{1+s}_\reg(\R/\Z,\R^n)$
into an arclength parametrized curve. The arclength $\psi_\g:
[0,1]\to [0,1]$
(rescaled to the unit interval and then extended to all of $\R$ such that 
$\psi_\g(t+1)=1+\psi_\g(t)$ for all $t\in\R$)  
is given by the strictly monotonic function 
$\psi_\g(t)\ceq\frac1{\LL(\g)}\LL(\g|_{[0,t]})$ for $t\in [0,1]$.
It has a strictly monotonic inverse
$\psi_\g^{-1}:[0,1]
\to [0,1]$.
We prove in \Cref{thm:arclength} 
in \cref{sec:arclength-reparam}
that the reparametrization and rescaling mapping
$\ARC:H^{1+s}_\reg(\R/\Z,\R^n)\to H^{1+s}_\unitspeed(\R/\Z,\R^n)$
given by
\begin{equation}\label{eq:ARC}
	\ARC(\g)
	\ceq
	\sdfrac{1}{\LL(\g)}\,\g\circ \psi_\g^{-1}\quad\Fo \g\in H^{1+s}_\reg(\R/\Z,\R^n)
\end{equation}
is continuous with respect to the $H^{1+s}$-topology for any $s\in (\frac{1}{2},1)$.
For $s=1$ this was proven by Reiter in \cite[Theorem 1.13]{reiter_2009}.
Notice that rescaling and reparametrization does not destroy the injectivity
of a curve so that we have the inclusion
$\ARC(H^{1+s}_\injreg(\R/\Z,\R^n))\subset
\AL^{1+s}$ for all $s\in (\frac{1}{2},1]$.

\begin{theorem}[Arclength manifolds are strong deformation retracts]
\label{thm:strong-deformation-ALs}
The open set\newline $H_\reg^{1+s}(\R/\Z,\R^n)\subset H^{1+s}(\R/\Z,\R^n)$ strongly
deformation retracts onto the submanifolds $\NL^{1+s}$ and 
$\NL^{1+s}_\bz$
for any $s\in (\frac{1}{2},1]$ and $\bz\in\R^n$, whereas $ H^{1+s}_{\injreg}\left( \R/\Z,\R^n \right) $ strongly deformation retracts onto $ \AL^{1+s} $ and $ \AL^{1+s}_{\bz} $. 
\end{theorem}

\begin{proof}
If $\g\in\NL^{1+s}$ then $\LL(\g)=1$ and $\psi_\g(t)=t$ for any $t\in\R$ so
that $\ARC(\g)=\g$ by definition \cref{eq:ARC}. 
Therefore, the mapping $\ARC \colon H^{1+s}_\reg(\R/\Z,\R^n)\to
H^{1+s}_\unitspeed(\R/\Z,\R^n)=\NL^{1+s}$ 
is continuous with $\ARC|_{\NL^{1+s}}=
\Id_{\NL^{1+s}}$, whence a retraction onto $\NL^{1+s}$.
For a given point $\bz\in\R^n$ we simply compose
the map $\ARC$ with the smooth
shift $\tau_\bz:H^{1+s}(\R/\Z,\R^n)\to
H^{1+s}(\R/\Z,\R^n)$ defined as $\tau_\bz(h) \ceq h+\bz-h(0)$
satisfying $\tau_\bz(\g)=\g$ for all $\g\in\AL^{1+s}_\bz$ 
so that we
obtain the retraction $\tau_\bz\circ\ARC:H^{1+s}_\injreg(\R/\Z,\R^n)
\to\NL^{1+s}_\bz$.

Now, we need to construct a homotopy between the identity map on the set
$\OL \ceq H^{1+s}_\reg(\R/\Z,\R^n)$ and the retraction onto $\NL^{1+s}$, and $\NL^{1+s}_\bz$, respectively.
Consider for any curve $\eta\in \OL$
the family of parameter transformations 
\begin{equation}\label{eq:familiy-param}
\phi_\eta(\sigma) \ceq (1-\sigma)\Id_\R+
\sigma\psi_\eta^{-1}:[0,1]\to [0,1]\quad\Fo \sigma\in [0,1],
\end{equation}
where, as before, $\psi_\eta^{-1}$ is the strictly monotonic inverse
of the rescaled arclength $\psi_\eta(t)=\frac1{\LL(\eta)}\LL(\eta|_{[0,t]})$.
As before, we may extend $\phi_\eta$ to all of $\R$ with
$\phi_\eta(t+1)=1+\phi_\eta(t)$ for all $t\in\R$.
Clearly, 
$\phi_\eta(\sigma)(0)=0$ and $\phi_\eta(\sigma)(1)=1$, and $\phi_\eta(\sigma)
(\cdot)$
is strictly
increasing for all $\sigma\in [0,1]$, since
$\phi_{\eta}(\sigma)'(t)=(1-\sigma)+\sigma (\psi_\g^{-1})'(t)>0$ for
any $t\in [0,1]$. We can interpret this family $\phi_\eta(\sigma)$ as 
a mapping $\phi_{\cdot}(\cdot):\OL\times [0,1]\to
H^{1+s}_\reg([0,1])$, which clearly is continuous on
$\OL\times [0,1]$ 
due to \Cref{cor:restriction-inverse-operator} and 
\Cref{lem:parametrization_continuous} in \cref{sec:arclength-reparam}.
Next, we modify the mapping $\ARC$ in \cref{eq:ARC} to obtain
a family of mappings $\CAPARC(\eta,\sigma)$ for $\sigma\in [0,1]$ defined as
\begin{equation}\label{eq:ARC-sigma}
\textstyle
\CAPARC(\eta,\sigma)\ceq\frac1{1+\sigma(\LL(\eta)-1)}\,\eta\circ\phi_\eta(\sigma)
\quad\Fo\eta\in \OL.
\end{equation}
\Cref{composition_continuous} in \cref{sec:arclength-reparam}
implies that the map 
$\CAPARC(\cdot,\cdot):\OL\times [0,1]\to \OL
$ is continuous. With $\phi_\eta(0)(t)=t$ for all $t\in\R$ one has $\CAPARC(\eta,0)=\eta$
for all $\eta\in \OL$, and $\CAPARC(\eta,1)\in\NL^{1+s}$ since $\phi_\eta(1)=
\psi_\eta^{-1}$ so that $\CAPARC(\eta,1)=\ARC(\eta)$ by \cref{eq:ARC}, 
and the original map
$\ARC$ was shown 
to be a retraction onto $\NL^{1+s}$.
For $\g\in\NL^{1+s}$ one has $\LL(\g)=1$ and $\psi_\g^{-1}(t)=t$ whence 
$\phi_\g(\sigma)(t)=t$ for all $t\in\R$ and $\sigma\in [0,1]$, 
and therefore 
$
\CAPARC(\g,\sigma)=\frac1{1+\sigma(1 -1)}\,\g=\g.$
Hence, $\CAPARC(\cdot,\cdot)$ is a strong deformation 
retraction of $\OL$ onto $\NL^{1+s}$.

If we fix a point $\bz\in\R^n$ to treat the arclength submanifold
$\NL^{1+s}_{\bz}$ instead of $\NL^{1+s}$, then we take 
the family of shifts $T_\bz(h,\sigma)\ceq
h+\sigma(\bz-h(0))$ for $h\in\OL$. 
Then the composition
$T_\bz(\CAPARC(\cdot,\sigma),\sigma)$ for $\sigma\in [0,1]$
yields a strong deformation retraction of $\OL$ onto $\NL^{1+s}_\bz$.
Indeed, we have
$T_\bz(h,0)=h$, 
$T_\bz(h,1)=h+\bz-h(0)$ hence $T_\bz(h,1)'=
h'$, and
$T_\bz(h,1)(0)=\bz$ for all $h\in \OL$.
So, $T_\bz(\CAPARC(\eta,0),0)=\CAPARC(\eta,0)=\eta$, and
$T_\bz(\CAPARC(\eta,1),1)=\CAPARC(\eta,1)+\bz-\CAPARC(\eta,1)(0)
\in\NL^{1+s}_\bz$ for every $\eta\in\OL$. In addition,
for $\g\in\NL^{1+s}_\bz$ we have shown above
that $T_\bz(\CAPARC(\g,\sigma),\sigma)=
T_\bz(\g,\sigma),$
which equals $\g$
since $\g(0)=\bz$.

Replacing $ \NL^{1+s} $ and $ H^{1+s}_{\reg} $ by $ \AL^{1+s} $ and $ H^{1+s}_{\injreg} $ respectively, yields the desired statement for the submanifolds of arclength parametrized knots. 
\end{proof}
\begin{remark}\label{rem:topological-justification-AL}
\cref{thm:strong-deformation-ALs} 
justifies restricting knot energies to the submanifold of arclength parametrized
curves topologically: As the knot space $H^{1+s}_\injreg(\R/\Z,\R^n)$ strongly
deformation retracts onto $\AL^{1+s}$ and also onto $\AL^{1+s}_\bz$, all
three topological spaces are homotopy equivalent and hence contain
the same topological information. 
\end{remark}
The last result of this section gives a justification from the perspective of knot theory. 
\begin{theorem}[Connected components are knot classes]
\label{thm:components-knot-classes}
Let $s\in (\frac{1}{2},1]$. If two knots $\g_0$ and $\g_1$ lie in the same
connected component of $\AL^{1+s}$ then they are $C^1$-ambiently isotopic.
If, on the other hand,
two knots $\g_0,\g_1\in\AL^{1+s}$ are $C^1$-isotopic in 
$\R^n$,
then they lie in the same connected component
of $\AL^{1+s}$.
All statements remain 
true if we replace $\AL^{1+s}$ by $\AL^{1+s}_\bz$ for some fixed point
$\bz\in\R^n$.
\end{theorem}
\begin{proof}
Fix $s\in (\frac{1}{2},1]$.
It suffices to treat $\AL \ceq \AL^{1+s}$.
If $\g_0,\g_1\in\AL$ are contained in the same connected component
of the smooth submanifold $\AL\subset\OL\ceq
H^{1+s}_\injreg(\R/\Z,\R^n)$ then there exists a path $c\in  
C^\infty([0,1],\AL)$ such that
$c(0)=\g_0$ and $c(1)=\g_1$; see, e.g., \cite[Lemma 3.1]{palais_1966}.
Set $J(t,\sigma) \ceq c(\sigma)(t)$ for $t\in\R/\Z$ and $\sigma\in [0,1]$, to
find that the curve
$J(\cdot,\sigma)\in\AL$ is a $C^1$-embedding for any $\sigma\in [0,1]$
by means of the Morrey embedding $H^{1+s}(\R/\Z,\R^n)\hookrightarrow
C^1(\R/\Z,\R^n)$; see \cite[Theorem A.2]{knappmann-etal_2023}. 
Moreover, one can show with this embedding 
that $J$ is also of class $C^1(\R/\Z\times [0,1],
\R^n)$ so that $J$ constitutes a $C^1$-isotopy with $J(\cdot,0)=\g_0$ and
$J(\cdot,1)=\g_1$. According to \cite[Theorem 1.2]{blatt_2009a} the two knots
$\g_0$ and $\g_1$ are therefore ambiently isotopic, i.e., in the same knot
class.

Let $J \in C^1(\Circle\times[0,1] ,\AmbSpace)$ be an isotopy satisfying $J(\cdot,0) = \gamma_0$ and $J(\cdot,1) = \gamma_1$. 
We consider the metric space $X \ceq C^1(\Circle,\AmbSpace)$ with the $C^1$-norm.
Define $\varGamma \in C^0([0,1],X)$ by $\varGamma(t) \ceq J(\cdot,t)$.
The set $U \ceq C^1_{\mathrm{ir}}(\Circle,\AmbSpace)$ is open in $X$,
and $\varGamma$ maps into it.
Observe that the distance function
\[
    \dist(\gamma, X \setminus U) \ceq \smash{\inf_{\eta \in X \setminus U} \|\gamma - \eta\|_{C^1}} \quad \text{ for }\gamma\in U
\]
is Lipschitz continuous.
Thus, the function $F \colon \intervalcc{0,1} \to \R$, $F(t) \ceq \dist(\varGamma(t), X \setminus U)$ is continuous. 
The interval $\intervalcc{0,1}$ is compact, so $F$ must have a minimizer $t_0 \in \intervalcc{0,1}$.
Because of $\varGamma(t_0) \in U$ and because $U$ is open, the minimum value must be positive: $2\varepsilon\ceq F(t_0) > 0$. 
Since $C^\infty(\Circle\times[0,1],\AmbSpace)$ is dense in $C^1(\Circle\times[0,1],\AmbSpace)$, there is an  $H \in C^\infty(\Circle\times[0,1],\AmbSpace)$ so that
$
    \norm{J - H}_{C^1(\Circle\times[0,1],\R^n)}
    <
    \varepsilon
$.
For the path $G$ defined by $G(t) \ceq H(t,\cdot)$ this means that
\[
\norm[\big]{ 
		t \mapsto \norm{ \varGamma(t) - G(t) }_{C^1} 
	}_{L^\infty} 
    \leq 
    \norm{J - H}_{C^1([0,1] \times \Circle,\R^n)}
    <
    \varepsilon
	.
\]
By construction, we have for every $t \in [0,1]$ that
\[
    G(t) 
	\in 
	B_\varepsilon(\varGamma(t)) \cap C^2(\Circle,\AmbSpace)
	\subset
	U \cap C^2(\Circle,\AmbSpace) 
    = C^2_{\mathrm{ir}}(\Circle,\AmbSpace)
	.
\]
Since the embedding $C^2(\Circle,\AmbSpace) \hookrightarrow H^{1+s}(\Circle,\AmbSpace)$ is linear and continuous, the path $G$ is an element of
$C^\infty([0,1],H^{1+s}_{\mathrm{ir}}(\Circle,\AmbSpace))$.
Now we define the linear homotopies
\[
    \varGamma_1(t) \ceq (1-t) \, \gamma_0 + t \, G(0)
    \quad \text{and} \quad 
    \varGamma_2(t) \ceq (1-t) \, G(1) + t \, \gamma_1
    ,
    \quad 
    \text{for $t \in [0,1]$}
    .
\]
Notice that $G(0) \in B_\varepsilon(\gamma_0)$, $G(1) \in B_\varepsilon(\gamma_1)$, that these balls contained in $ C^1_{\injreg}(\Circle,\AmbSpace) $ are convex,  and that knot classes are stable under small $ C^1 $-pertubations, see \cite{blatt_2009a}. 
Thus, we have $\varGamma_1, \varGamma_2 \in C^\infty([0,1],H^{1+s}_{\mathrm{ir}}(\Circle,\AmbSpace))$.
Concatenation of the paths $\varGamma_1$, $G$, and $\varGamma_2$ yields a continuous path from $\gamma_0$ to $\gamma_1$ in $H^{1+s}_{\mathrm{ir}}(\Circle,\AmbSpace)$.
Finally, we retract this path to a continuous path in $\mathcal{A}$, by means of \cref{thm:strong-deformation-ALs}.
\end{proof}

\subsection{Critical knots in arclength submanifolds}\label{sec:AL-critical}
In this abstract setting we consider for $s\in (\frac{1}{2},1]$  
an arbitrary  continuously differentiable parameter-invariant
energy 
$$
\EL \colon \OL \ceq H^{1+s}_\reg(\R/\Z,\R^n)\to\R.
$$
The latter means that 
$\EL(\eta)=\EL(\eta\circ\psi)$ for any $\eta\in\OL$,
where $\psi\in H^{1+s}([0,1])$ is a
strictly increasing parameter transformation with $\psi(0)=0$ and $\psi(1)=1$,
which can be extended to all of $\R$ to satisfy $\psi(t+1)=1+\psi(t)$
for all $t\in\R$.  

One of the simplest examples of such a parameter-invariant functional is the
length $\LL \colon \OL\to\R$, which is actually a smooth submersion, since its differential
$\D\LL_\eta\colon T_\eta\OL\simeq H^{1+s}(\R/\Z,\R^n)\to T_{\LL(\eta)}\R\simeq\R$
satisfies $\D\LL_\eta(\eta)=\LL(\eta)$ for all $\eta\in\OL$ because of the 
positive
$1$-homogeneity of $\LL$ and Euler's homogeneous function 
theorem. 
Hence, by the preimage theorem \cite[Theorem 4.J]{zeidler_1993}, the set
\begin{equation}\label{eq:Omega-s}
\varOmega^{1+s}\ceq\{\g\in\OL:\LL(\g)=1\}=\LL^{-1}(1),
\end{equation}
is a smooth submanifold of $\OL$
with its tangent space given by
\begin{equation}\label{eq:TgOmega-s}
	T_\g\varOmega^{1+s}
	=
	\bigg\{
		\nu\in H^{1+s}(\R/\Z,\R^n) 
		: 
		\D\LL_\g(\nu)=\int_0^1\sdfrac{\langle\g'(u),\nu'(u)\rangle}{|\g'(u)|}\d u=0
	\bigg\}.
\end{equation}
We also consider the smooth submanifolds with a fixed point $\bz\in\R^n$:
\begin{equation}\label{eq:Omega-sz}
\varOmega^{1+s}_\bz\ceq\{\g\in\varOmega^{1+s}:\g(0)=\bz\}
\end{equation}
with tangent space
\begin{equation}\label{eq:TgOmega-sz}
T_\g\varOmega^{1+s}_\bz=\{\nu\in T_\g\varOmega^{1+s}:\nu(0)=0\}.
\end{equation}
We now relate the arclength manifolds $\NL^{1+s}$ and $\NL^{1+s}_{\bz}$ 
to $\varOmega^{1+s}$ and
$\varOmega^{1+s}_{\bz}$ via general parameter-invariant energies $\EL\in C^1(\OL)$
and their critical points under a fixed-length constraint.

Note that as before, since $ \AL^{1+s} $ and $ \AL^{1+s}_{\bz} $ are open subsets of $ \NL^{1+s} $ and $ \NL_{\bz}^{1+s} $ respectively, the results of this section hold as well for arclength parametrized knots if we choose $ \OL$ to be the open set $ H^{1+s}_{\injreg}(\R/\Z,\R^n) $. This is due to the fact that the corresponding tangent spaces coincide; see \cref{eq:tangent-space-AL,eq:tangent-space-ALz}.
Only $|\g'|=1$ will be used in the arguments of this subsection, injectivity of $\g$ does not play a r{\^o}le here.

\begin{theorem}[Arclength criticality vs. criticality under length constraint]
\label{thm:critical-points}
Let $s\in (\frac{1}{2},1]$ and $ \ML\in \{\NL^{1+s},\AL^{1+s}\} $.
If $\g\in\ML\cap H^{2+s}(\R/\Z,\R^n)$ is a critical point of the restricted
parameter-invariant energy $\widetilde{\EL}\ceq\EL|_{\ML}$, 
then $\g$ is a critical
point of $\EL-\lambda\LL$ on $\OL$ for some $\lambda\in\R$. If  
$\EL$ is 
positively $\alpha$-homogeneous for some 
$\alpha\in\R$, then
$\lambda = \alpha \, \EL(\g) = \alpha \, \widetilde{\EL}(\g)$. In particular, if $\alpha\not=0$ and
$\widetilde{\EL}(\g)\not= 0$ then $\lambda\not= 0$, whereas
if $\EL$ is scale-invariant or $\widetilde{\EL}(\g)=0$, then $\lambda=0$.
\end{theorem}

Notice that we require higher regularity of the critical point, 
which we will verify for the energies in
our applications in geometric knot theory in \cref{sec:Smoothness}. 
This higher regularity  is necessary in
the following decomposition of tangent vectors $\nu\in T_\g\varOmega^{1+s}$
into a tangential and a non-tangential part  with respect to $\g$. 
We lose one derivative of $\g$ when we split off the tangential part.
\begin{lemma}[Decomposition of tangent vectors]
\label{lem:decomposition}
For any $\g\in \NL^{1+s}\cap H^{2+s}(\R/\Z,\R^n)\subset\varOmega^{1+s}$ 
and $\nu\in T_\g\varOmega^{1+s}$ there exists
$h\in T_\g\NL^{1+s}$ and a function $\varrho\in H^{1+s}(\R/\Z)$ with $\varrho(0)=\varrho(1)=0$, such 
that 
\begin{equation}\label{eq:decomposition}
	\nu(t)=h(t) + \varrho(t) \, \g'(t)\quad\Foa t\in\R/\Z.
\end{equation}
Likewise, if $\g\in\NL^{1+s}_{\bz}\cap H^{2+s}(\R/\Z,\R^n)$ for some fixed $\bz\in\R^n$
then every $\nu\in T_\g\varOmega^{1+s}_\bz$ can be decomposed as in \cref{eq:decomposition}
where now $h\in T_\g\NL^{1+s}_\bz$.
\end{lemma}
\begin{proof}
Set $\varrho(t)\ceq\int_0^t \inner{\nu'(u),\g'(u)}_{\R^n} \dd u$, which is of class
$H^{1+s}(\R/\Z)$ due to the product rule (see, e.g., \cite[Proposition A.5]{knappmann-etal_2023}), and satisfies $\varrho(0)=0=\varrho(1)$ since $\nu\in T_\g\varOmega^{1+s}$ and
$\abs{\g'}=1$ on $\R/\Z$. The function $h\ceq\nu-\varrho\g'$ is in $H^{1+s}(\R/\Z,\R^n)$
because of our regularity assumption on $\g$ and by the product rule, and one
computes (abbreviating $\inner{\cdot,\cdot} \ceq \inner{\cdot,\cdot}_{\R^n}$)
\begin{align*}
	\inner{h',\g'} 
	&=
	\inner{\nu',\g'} -\varrho' \inner{\g',\g'}-\rho \inner{\g'',\g'}
	=
	\inner{\nu',\g'} - \inner{\nu',\g'} \inner{\g',\g'}
	=
	0,
\end{align*}
where we used the definition of $\varrho$ and the fact that $\g$ is arclength parametrized.
Therefore, $h\in T_\g\NL^{1+s}$ (see \cref{eq:tangent-space-AL}), which establishes the decomposition \cref{eq:decomposition}.

With the same Ansatz for $\varrho$ one obtains for $\g\in\NL^{1+s}_\bz$ 
and $\nu
\in T_\g\varOmega_\bz^{1+s}$ that $h(0)=\nu(0)-\varrho(0)\g'(0)=0$ so that 
$h\in T_\g\NL^{1+s}_\bz$ by means of \cref{eq:tangent-space-ALz}.
\end{proof}
With this decomposition of tangent vectors of $\varOmega^{1+s}$ and 
$\varOmega^{1+s}_\bz$ we can
now prove criticality under the fixed-length constraint.

\begin{refproof}[Proof of \cref{thm:critical-points}]
As mentioned earlier, it suffices to prove the theorem for $ \ML=\NL^{1+s} $. For arbitrary $\nu\in T_\g\varOmega^{1+s}$ we find 
according to \cref{lem:decomposition}
some $h\in T_\g\ML$ and a function $\varrho\in H^{1+s}(\R/\Z)$ with $\varrho(0)=\varrho(1)=0$ such that
\cref{eq:decomposition} holds true. Consider the family $\varphi_\tau\ceq
\Id_\R + \tau\varrho\in H^{1+s}([0,1])$ satisfying $\varphi_\tau(0)=0$ and
$\varphi_\tau(1)=1$, and 
$$
\varphi_\tau'(t)=1+\tau\varrho'(t)\ge 1 - \abs{\tau} \norm{\varrho'}_{C^0}\ge
1-\tau_0 \,C_E \norm{\varrho}_{H^{1+s}}>0\quad\Foa t\in [0,1],
$$
if $|\tau|\le\tau_0\ll 1$. Here, $C_E$ denotes the constant of the Morrey embedding
theorem \cite[Theorem A.2]{knappmann-etal_2023}. So, $\psi\ceq\varphi_\tau$ is an
admissible strictly increasing parameter transformation under which $\EL$ is assumed to be invariant, i.e.,
the function $\tau\mapsto \EL(\g\circ\varphi_\tau)$ is constant on $(-\tau_0,\tau_0)$.
With 
the additional properties 
$$
	\textstyle
	\frac{\partial}{\partial \tau}\varphi_\tau(t) |_{\tau=0}
	=
	\varrho(t)
	\quad\AND\quad
	\varphi_0(t) = t \quad\Foa t\in [0,1],
$$
we conclude
\begin{align*}
	0 &
	\textstyle
	= 
	\frac{\partial}{\partial \tau}\EL(\g\circ\varphi_\tau) \big|_{\tau=0}
	= 
	\D\EL_\g \!\pars[\big]{ \pars{\g'\circ\varphi_0} \cdot\frac{\partial}{\partial \tau}\varphi_\tau |_{\tau=0} }
	=
	\D\EL_\g(\varrho\g')
	\\
	&\overset{\cref{eq:decomposition}}{=}
	\D\EL_\g(\nu-h)
	=
	\D\EL_\g(\nu)-\D\EL_\g(h)
	=
	\D\EL_\g(\nu)-\D\widetilde{\EL}_\g(h)
	=
	\D\EL_\g(\nu)
	.
\end{align*}
Since $\nu\in T_\g\varOmega^{1+s}$ 
was arbitrary we infer from the Lagrange multiplier 
theorem (see \cite[Corollary 3.5.29]{abraham-etal_1988}) 
that there is some number $\lambda\in\R$ such that
$(\D\EL-\lambda \D\LL)_\g=0$.

For the evaluation of $\lambda$ we use that $\LL(\g)=1$,
the positive
$1$-homogeneity of $\LL$, and the positive $\alpha$-homogeneity of $\EL$ to
compute with Euler's homogeneous function theorem
$
\lambda = \lambda \, \LL(\g) = \lambda \D\LL_\g(\g) = \D\EL_\g(\g) = \alpha \, \EL(\g) =
 \alpha \, \widetilde{\EL}(\g).
$
\end{refproof}

Fixing a point $\bz\in\R^n$ adds a vectorial constant Lagrange multiplier.
\begin{corollary}[Criticality with fixed point]
\label{cor:critical-points}
Let $s\in (\frac{1}{2},1]$, $\bz\in\R^n$ and $ \ML_{\bz}\in \{\NL^{1+s}_{\bz},\AL^{1+s}_{\bz}\} $. If $\g\in\ML_{\bz}\cap
H^{2+s}(\R/\Z,\R^n)$ is a critical knot of the restricted 
parameter-invariant energy
$\widetilde{\EL}\ceq\EL|_{\ML_{\bz}}$, then $\g$ is critical for
$\EL-\lambda \, \LL- \inner{\Upsilon,P_0(\cdot)}$ for some
$\lambda\in\R$ and $\Upsilon\in\R^n$. 
If $\EL$ is 
positively $\alpha$-homogeneous for some $\alpha\in\R$,
then $\lambda + \inner{\Upsilon,\bz}=\alpha \, \widetilde{\EL}(\g).$
\end{corollary}
\begin{proof}Again we may focus on $ \ML_{\bz}=\NL^{1+s}_{\bz} $.
We can use \cref{lem:decomposition} for an arbitrary tangent vector
$\nu\in T_\g\varOmega^{1+s}_\bz$ to obtain \cref{eq:decomposition}
with a function $\varrho\in H^{1+s}(\R/\Z)$ with $\varrho(0)=\varrho(1)=0$
and tangent vector $h\in T_\g\ML_{\bz}$. Constructing a family of
admissible parameter transformations as in the previous proof we obtain
$\D\EL_\g(\nu)=0$. Observing that we can represent 
the smooth submanifold $\varOmega_\bz^{1+s}$ as
preimage as $\varOmega_\bz^{1+s}=(\LL\times P_0)^{-1}(1,\bz)$,
where $P_0$ is the point evaluation defined in \cref{eq:point-evaluation},
we can apply the
Lagrange multiplier rule
\cite[Corollary 3.5.29]{abraham-etal_1988}
to obtain some number $\lambda\in\R$ and a vector $\Upsilon\in\R^n$ such
that $\g$ is a critical point of the functional $\EL-\lambda\LL -\langle
\Upsilon,P_0(\cdot)\rangle$. 
For the relation between $\lambda$ and $\Upsilon$ we compute similarly
as in the previous proof $ \lambda + \inner{\Upsilon,\bz}
	 =
	\lambda \, \LL(\g) + \inner{\Upsilon,P_0(\g)}
	=
	\lambda \D\LL_\g(\g) + \inner{ \Upsilon,(\D P_0)_\g(\g) }
	=
	\D\EL_\g(\g)
	= 
	\alpha \, \widetilde{\EL}(\g). $
\end{proof}

If we assume  
that the energy $\EL$ is also invariant under translation, 
we do not have to deal with the  additional
Lagrange parameter $\Upsilon\in\R^n$.

\begin{lemma}\label{cor:critical-points_wo_fixed_point}
Let $ \ML=\{\NL^{1+s},\AL^{1+s}\} $ and $ \ML_{\bz}=\{\eta\in \ML:\eta(0)=\bz\} $. If the energy $\EL\in C^1(\OL)$ is invariant under translations, then 
a critical point $\gamma \in \ML_\bz$ of $\EL|_{\ML_\bz}$ is also a 
critical point of $\EL|_{\ML}$.
\end{lemma}
\begin{proof}
	Since the energy is invariant under translations we have 
\[
	\textstyle
	\D\EL_\g p
	=
	\frac{\partial}{\partial\varepsilon}
	\EL(\gamma+\varepsilon p)
	\big|_{\varepsilon=0}
	=
	0\quad\Foa p\in\R^n
	.
\]
In particular, this implies
\[
	\textstyle
	0=\D\EL_\g h(0)
	=
	\D\EL_\g h - \D\EL_\g (h-h(0))
	=
	\D\EL_\g h
	\quad
	\Foa
	h\in T_\g\ML
	,
\]
since $h-h(0)\in T_\g\ML_\bz$;  see \cref{eq:tangent-space-AL} and
\cref{eq:tangent-space-ALz}.
\end{proof}

%% file: PalaisSmale.tex
\section{The Palais--Smale condition for bending energy and knot energies}
\label{sec:PS-bending+knot-energies}

\subsection{Bending energy}
\label{sec:PS-bending}
We consider the Euler-Bernoulli bending energy $\BendingEnergy$ defined in
\cref{eq:bending-energy}
on the open subset
$\OL \ceq H^2_\reg(\R/\Z,\R^n)$ of immersions in the Hilbert space 
$\HL \ceq H^2(\R/\Z,\R^n)$, which itself embeds compactly into 
$\BL\ceq C^1(\R/\Z,\R^n)$. Its differential 
$\D(\BendingEnergy)_\g:
T_\g \OL\simeq\HL\to T_{\BendingEnergy(\g)}\R\simeq\R $
may be written in a parameter-invariant representation (see, e.g., 
\cite[Proof of Lemma A.1]{dallacqua-pozzi_2014})
as
\begin{align*}
        \D(\BendingEnergy)_{\gamma}h\textstyle
= \intRZ\left\langle \kappa,\partial_s\nabla_s h \right\rangle \d s 
- \frac{1}{2} \intRZ \left|\kappa \right|^2 \left\langle \partial_s 
\gamma,\partial_s h \right\rangle \d s\quad\Fo \g\in\OL,\,h\in\HL,
\end{align*}
where $ \kappa=\partial^2_s \gamma $ is the curvature vector field of 
$ \gamma $,
 $ \partial_s =\frac{1}{|\gamma'|}\partial_u $ the 
 arclength derivative, $ \nabla_s h = \P^{\perp}_{\gamma'}(\partial_s h)=
 \partial_s  h-\langle\partial_s h,\frac{\g'}{|\g'|}\rangle
 \frac{\g'}{|\g'|}$
 the normal derivative,
 and $\d s=|\g'(u)|\d u$ denotes
integration by arclength. 
If $\g\in \NL^{2}\subset\OL\subset\HL$ 
then the differential
$\D(\BendingEnergy)_\g$ 
simplifies (because of $|\g'|=1$ and $\g''\perp \g'$) to
\begin{align}
\D(\BendingEnergy)_{\gamma}h&\textstyle
= \intRZ \big\langle \gamma'',( \P_{\gamma'}^\perp(h') )' 
\big\rangle \d u 
-\frac{1}{2} \intRZ\left|\gamma'' \right|^2 
\left\langle \gamma',h' \right\rangle \d u\notag\\
                &
		\textstyle
		=\intRZ 
        \pars[\big]{
            \inner{ \gamma'',h'' } 
            - 
            \frac{3}{2}\abs{\gamma'' }^2
            \inner{ \gamma',h' }
        }
        \d u
		\quad\Foa h\in\HL.
        \label{eq:Diff_elastic_formula}
\end{align}

\begin{theorem}
\label{thm:PS_elastic}
Let $\bz\in\R^n$ be fixed. Then 
the restriction $\widetilde{E}_b\ceq\BendingEnergy|_{\ML}$ of the elastic energy $\BendingEnergy$ 
to any one of the submanifolds $\ML_{\bz}\in\{\AL^2_\bz,\NL^2_\bz\}$ contained in $\OL\subset\HL\hookrightarrow\BL$
satisfies conditions \ref{lower-bound}, \ref{condB*},
\ref{condS} and \ref{condD}. 
\end{theorem}
\begin{proof}
Condition \ref{lower-bound} is trivially satisfied for $e\ceq0$, and $\BendingEnergy\in C^1(\OL)$, which follows from \cite[Proposition 4.8]{okabe-schrader_2023}.
With 
$$\textstyle
\sup_{k\in\N}\widetilde{E}_b(\g_k)=\sup_{k\in\N}\|\g_k''\|_{L^2}^2<\infty
\quad\AND\quad\|\g_k\|_{C^1}\le |\bz|+2\quad\Foa k\in\N
$$ 
for any sequence
$(\g_k)_k\subset\ML_{\bz}$
we obtain $\sup_{k\in\N}
\norm{\g_k}_\HL<\infty$ proving condition \ref{condB*}.

Concerning the splitting \ref{condS}, we add the inner product
$\pm\langle\g,h\rangle_{H^1}$ of the Hilbert space $H^1(\R/\Z,\R^n)$
to \cref{eq:Diff_elastic_formula}
to arrive at
\begin{align}\label{eq:elastic-splitting}
        \D(\BendingEnergy)_\g h&\textstyle
        =
        \langle\g,h\rangle_{\HL} - \int_{\R/\Z}\big(\langle\g',h'\rangle
        (1+\frac{3}{2}|\g''|^2) +\langle\g,h\rangle\big) \d u
        \notag\\
        &\qec  \textstyle
        \QL(\g,h)+\RL_\g \iota_{\BL}h \quad\Foa h\in\HL,\gamma\in \ML_{\bz}.
\end{align}
Clearly,  the inner product $\QL$ is a continuous bilinear form on
$\HL\times\HL$ satisfying the requirement \ref{condS2} in
the splitting \ref{condS} with coercivity
constant $c=1$, here. We use the constant function $\mu \equiv 1$, so condition \ref{condS1} is satisfied, as well. The remainder term $\RL$ can be bounded as follows:
\begin{equation}
    \label{eq:remainder-estimate}
    \abs{ \RL_\g \, h } 
    \le\textstyle
    \norm{h}_{C^1}\int_{\R/\Z}\big(1+\frac{3}{2}|\g''|^2+|\g|\big)\d u\le
    \norm{h}_{C^1}\big(2+|\bz|+\frac{3}{2} \norm{\g}^2_\HL\big)\,\,\Foa h\in \BL.
\end{equation}
So, $\RL$ is bounded on all bounded sets and hence, on all \Econtrolled{\BendingEnergy} sets. 
This verifies \ref{condS3}.

We observe from \cref{eq:Diff_elastic_formula} that $\D(\BendingEnergy)_\g \, h$ contains only derivatives of the tangent vector
$h\in\HL$, which in condition \ref{condD} is replaced by
the term $(\Id_\HL-\pr_\g^\ML) \, h$.
Therefore,
by means of \cref{eq:projection-ALsz,eq:Diff_elastic_formula,eq:projection-ALs} 
\[\D(\BendingEnergy)_{\gamma}(\Id_{\HL}-\pr^{\ML}_{\gamma})h=\D(\BendingEnergy)_{\gamma}Y_{\gamma}\left\langle \gamma',h' \right\rangle.\]
To prove \ref{condD} we need to analyze the mapping
\begin{equation}\label{eq:conD_mapping_EB}
\begin{aligned}\textstyle
    f\mapsto\D(\BendingEnergy)_{\gamma}Y_\gamma f \overset{\cref{eq:Diff_elastic_formula}}{=}\int_{\R/\Z}\left\langle \gamma'',(Y_\gamma f)'' \right\rangle \d u-\frac{3}{2}\int_{\R/\Z}\left|\gamma'' \right|^2 \left\langle \gamma',(Y_\gamma f)'  \right\rangle\d u,
\end{aligned}
\end{equation}
where it suffices to compute these terms for $ f\in C^{\infty}(\R/\Z) $. By means of the explicit formulas in \cref{eq:derivative-right-inverse} we obtain
\begin{equation}
\begin{aligned}\label{eq:first_derivative_rightInverse}
    \left\langle \gamma',(Y_\gamma f)' \right\rangle =f
\end{aligned}
\end{equation}
and 
\begin{equation}\label{eq:second_derivative_rightInverse}
\begin{aligned}
    \left( Y_\gamma f \right)''&=f'\gamma'+f\gamma''+(\P^{\perp}_{\gamma'})'(V_{\gamma}f)\\
    &=f'\gamma'+f\gamma''-\left\langle V_{\gamma}f,\gamma'' \right\rangle \gamma'-\left\langle V_{\gamma}f,\gamma' \right\rangle\gamma''.
\end{aligned}
\end{equation}
Because of $ \left|\gamma' \right|=1 $ one has $ \gamma'\perp\gamma'' $; hence 
\begin{equation}\label{eq:product_second_derivative}
\begin{aligned}
    \left\langle \gamma'',(Y_\gamma f)'' \right\rangle =\left|\gamma'' \right|^2\left(f-\left\langle V_{\gamma}f,\gamma' \right\rangle \right).
\end{aligned}
\end{equation}
Substituting \cref{eq:first_derivative_rightInverse,eq:second_derivative_rightInverse} into \cref{eq:conD_mapping_EB} yields 
\begin{equation}\label{eq:Eb_lowerorderrightinverse}
\begin{aligned}\textstyle
    \D(\BendingEnergy)_{\gamma}Y_{\gamma}f=-\frac{1}{2}\int_{\R/\Z}\left|\gamma'' \right|^2f\d u-\int_{\R/\Z}\left|\gamma'' \right|^2 \left\langle V_{\gamma}f,\gamma' \right\rangle\d u.
\end{aligned}
\end{equation}
The estimate \cref{eq:V-h-est} leads to $ \left|\D(\BendingEnergy)_{\gamma}Y_\gamma f \right|\le C\left\|f \right\|_{C^0} $ for all $ f\in  C^\infty(\R/\Z) $, where $ C=C(s,\left\|\gamma \right\|_{\HL}) $ depends non-decreasingly on the norm $ \left\|\gamma \right\|_{\HL} $. 
Therefore, by density of $ C^\infty(\R/\Z) $ in $ C^0(\R/\Z) $, the map in \cref{eq:conD_mapping_EB} extends to a continuous linear map on $ C^0(\R/\Z) $, which implies that 
\[S_{\gamma}\ceq \D(\BendingEnergy)_{\gamma}Y_\gamma \left\langle \gamma',\cdot \right\rangle\in \BL^* \]
is bounded on \Econtrolled{\BendingEnergy} sets as required in \ref{condD}. 
\end{proof}

From the
compact embedding $\HL\hookrightarrow\BL=C^1(\R/\Z,\R^n)$ we obtain
for any sequence $(\g_k)_k\subset\NL^2_{\bz}$ with $\g_k\rightharpoonup\g\in\HL$ that
$|\g'|=\lim_{k\to\infty}|\g_k'|=1$ and $\g(0)=\lim_{k\to\infty}\g_k(0)=\bz$; 
hence $\g\in\NL^2_{\bz}$, which establishes condition \ref{condL*} for 
$\widetilde{\FL}\equiv\widetilde{\EL} \ceq \BendingEnergy|_{\NL^2_{\bz}}$; and hence \ref{condL}.
Moreover, \ref{condE} and \ref{condP}
hold true according to \Cref{cor:basic-assumptions} for $s=1$, and
condition \ref{condB*} established in \Cref{thm:PS_elastic}
implies \ref{condB}, which allows us
to apply \Cref{thm:PS-abstract} to $\widetilde{\EL}=\BendingEnergy|_{\NL^2_{\bz}}$
to reprove the results of \cite[Theorem 1.7]{langer-singer_1985},  
\cite[Section 5 Theorem 5.6]{schrader_2016}, and
\cite[Proposition 5.4]{okabe-schrader_2023}, also cf.
to \cite[Proposition 3.2]{linner_1998}.

\begin{bcorollary}[Bending energy satisfies (PS) on $\NL^2_{\bz}$]
\label{cor:PS_elastic}
The restricted bending energy $\BendingEnergy|_{\NL^2_{\bz}}$ satisfies the
{\rm (PS)}-condition. 
\end{bcorollary}

As a consequence we can apply \Cref{thm:gradflow_subconvergence} and
\cref{koro:critpoint_lowerenergy,koro:minimizer_existence}   
to conclude\footnote{In \cite[Lemma 3.1]{okabe-schrader_2023}
$\BendingEnergy$ was shown
 to be of class
$C^{1,1}_\textnormal{loc}$ on $H^2_\reg(\R/\Z,\R^n)$.}
for $\widetilde{E}_b \ceq \BendingEnergy|_{\NL^2_{\bz}}$
with the following.

\begin{bcorollary}[Critical points and gradient flow for $\BendingEnergy$ on $\NL^2_{\bz}$]
\label{cor:elastic_crit_gradient-flow}
For any curve $\g\in\NL^2_{\bz}$ such that $\D(\widetilde{E}_b)_\g\not=0$ there is
some $\g_c\in\NL^2_{\bz}$ with $\D(\widetilde{E}_b)_{\g_c}=0$ and 
$\widetilde{E}_b(\g_c)<\widetilde{E}_b(\g)$. There is an 
$\widetilde{E}_b$-minimizing immersion
$\g_\textnormal{min}\in\NL^2_{\bz}$, i.e., $\widetilde{E}_b(\g_\textnormal{min})=
\inf_{\NL^2_{\bz}}\widetilde{E}_b(\cdot).$
Moreover, we have long-time existence for 
the gradient flow of $\widetilde{E}_b$ 
on $\NL^2_{\bz}$ along which $\widetilde{E}_b$ is strictly decreasing 
unless the initial curve is an
$\widetilde{E}_b$-critical point.
Furthermore, the flow subconverges to an $\widetilde{E}_b$-critical
point as time tends to infinity.
\end{bcorollary}
Before proving the Palais-Smale condition for linear combinations of bending and self-avoidance energies let us point out here that the arguments to establish \ref{condD} for $ \BendingEnergy $ in \cref{thm:PS_elastic} can be generalized considerably. The following lemma will be applied later to verify the (PS)-condition for the tangent-point energy and also in \cref{sec:Smoothness} to prove regularity of critical points. 
\begin{lemma}\label{lem:DEbYisLowerOrderIIII}
    Let $ \gamma\in H^{t}(\R/\Z,\R^n) $. 
    \begin{enumerate}
        \item[\textrm{(i)}] If $ t>\frac{3}{2} $, then the map $ f\mapsto V_{\gamma}f$ defined by \cref{eq:right-inverse-Y} extends continuously to a linear map in $ L\left( H^{1-t}(\R/\Z),\R^n \right)= L\left( (H^{t-1}(\R/\Z))^*,\R^n \right)$. 
        \item[\textrm{(ii)}] If $ \gamma\in \NL^{t} $ and $ t\ge 2 $, then 
        \begin{align*}
            \D(\BendingEnergy)_{\gamma}Y_\gamma\in H^{\sigma(t-2,\varepsilon)}(\R/\Z,\R^n)\quad \Foa \varepsilon>0,
        \end{align*}
        where 
        \begin{align*}
            \sigma(\tau,\varepsilon)\ceq 
            \begin{cases}
                2\tau-\frac{1}{2}-\varepsilon &\quad\Fo \tau\in \left[0,\frac{1}{2}\right],\\
                \tau &\quad \Fo\tau>\frac{1}{2}.
            \end{cases}
        \end{align*}
    \end{enumerate}
\end{lemma}
\begin{proof}
    \textrm{(i)}  Let $ \{ e_1,\dots,e_n\} $ be an orthonormal basis of $ \R^n $. In view of 
    \begin{align*}\textstyle
        \left|V_{\gamma}f \right|\overset{\cref{eq:right-inverse-Y}}{\le}\left\|\Theta^{-1}_{\gamma} \right\|_{L(\R^n)}\left|\int_{\R/\Z}f\gamma'\d x \right|=\left\|\Theta^{-1}_{\gamma} \right\|_{L(\R^n)}\left|\sum_{l=1}^ne_l \left\langle f,\gamma'_l \right\rangle_{L^2}  \right|
    \end{align*}
    for $ \gamma'_l\ceq \left\langle \gamma',e_l \right\rangle  $, it suffices to express the $ L^2 $-inner product for each $ l=1,\dots,n $ in terms of Fourier coefficients. Indeed, 
    \begin{align*}\textstyle
        \sum_{k\in\Z}\widehat{(\gamma'_{l})}_{k}\bar{\hat{f}}_{k}=\sum_{k\in\Z}\left|k \right|^{t-1}\widehat{(\gamma'_{l})}_{k}\bar{\hat{f}}_{k}\left|k \right|^{1-t}
    \end{align*}
    is finite if $ (\hat{f}_{k}\left|k \right|^{1-t} )_{k}\in \ell^2$, i.e., if $ f\in H^{1-t}(\R/\Z) $, because $ (\widehat{(\gamma'_{l})}_{k}\left|k \right|^{t-1} )_{k}\in \ell^2 $ by virtue of our regularity assumptions on $ \gamma $. 

    \textrm{(ii)}  According to the product rule, \cref{cor:gamma''gamma''}, the mapping $ \gamma\mapsto\left|\gamma''  \right|^2\qec g $ is of class $ C^0\left( H^{t}(\R/\Z,\R^n),H^{\sigma(t-2,\varepsilon)}(\R/\Z) \right) $ for all $ \varepsilon>0 $. 
    Therefore,  
    \begin{align*}\textstyle
        \left\langle g,f \right\rangle_{L^2} =\sum_{k\in\Z}\hat{g}_{k}\hat{f}_{k}=\sum_{k\in\Z}\hat{g}_{k}\left|k \right|^{\sigma(t-2,\varepsilon)}\left|k \right|^{-\sigma(t-2,\varepsilon)}\hat{f}_{k}
    \end{align*}
    is finite if $ (\left|k \right|^{-\sigma(t-2,\varepsilon)}\hat{f}_{k})_{k}\in \ell^2 $, or $ f\in H^{-\sigma(t-2,\varepsilon)}\left( \R/\Z \right) $. This shows that $ C^\infty(\R/Z)\ni f\mapsto\int_{\R/\Z}\left|\gamma'' \right|^2 f\d x $ extends to a continuous linear map of class $ H^{\sigma(t-2,\varepsilon)}(\R/Z) $.

    According to part \textrm{(i)} we find with $ \left|\gamma' \right|=1 $, 
    \begin{align*}\textstyle
      \big|\langle \int_{\R/\Z}\left|\gamma'' \right|^2\gamma'\d x, V_{\gamma}f \rangle  \big|\le \left\|\gamma'' \right\|^2_{L^2}\left|V_{\gamma}f \right|\overset{\textrm{(i)}}{\le}\left\|\gamma'' \right\|^2_{L^2}\left\|f\mapsto V_{\gamma}f \right\|_{L(H^{1-t},\R^n)}\left\|f \right\|_{H^{1-t}}
    \end{align*}
    for all $ f\in C^\infty(\R/\Z) $, so that $ f\mapsto \langle \int_{\R/\Z}\left|\gamma'' \right|^2\gamma'\d x, V_{\gamma}f \rangle $ extends to a continuous linear map of class $ H^{t-1}(\R/\Z)\equiv \left( H^{1-t}(\R/\Z) \right)^* $. 
    In view of the explicit expression \cref{eq:Eb_lowerorderrightinverse} for $ \D(\BendingEnergy)_{\gamma}Y_\gamma $ this proves the lemma, because $ H^{t-1}(\R/\Z)\subset H^{\sigma(t-2,\varepsilon)(\R/\Z)} $ for all $ t\ge 2 $.
\end{proof}

\subsection{Bending energy with self-avoidance term}
\label{sec:PS-elastic-knot}
Now we deal with linear combinations of bending energy $\BendingEnergy$ and some choice
of knot energy $\FL$ in any one of the three knot energy families introduced
in \cref{sec:intro}, i.e., 
$\FL\in\{\TP^{(p,q)},\intM^{(p,q)},E^{\alpha,p}\}$ 
in the respective range \cref{eq:individual-energy-spaces}
for the parameters $p,q,$ and $\alpha,p$. 
For each energy space
 one has $s>\frac1{\rho}$ so that the fractional Sobolev
space $W^{1+s,\rho}(\R/\Z,\R^n)$
embeds compactly into $C^1(\R/\Z,\R^n)$; see
\cite[Theorem A.2]{knappmann-etal_2023}. Since $W^{2,2}$ embeds 
continuously into $W^{1+\sigma,2}$ for any $\sigma\in (0,1)$ (see, e.g.,
\cite[Lemma 2.20]{steenebruegge_2023}), and the latter embeds
compactly into $W^{1+s,\rho}(\R/\Z,\R^n)$ as long as
$s\in (0,1)$, $\rho\in (1,\infty)$ and $\sigma-s >\max\{
\frac{1}{2}-\frac1{\rho},0\}$ by virtue of \cite[Proposition A.3(i)]{matt-etal_2023},
we can check for each energy separately 
under the additional restriction \cref{eq:shrink-range}
of the parameter range that we
have a compact embedding of $W^{2,2}(\R/\Z,\R^n)$ into
the respective energy space. Therefore, to set the stage for the application of \cref{cor:PS-total-energy}, we define $\HL \ceq H^2(\R/\Z,\R^n)$,
$\OL \ceq H^2_\injreg(\R/\Z,\R^n)$, $\BL \ceq C^1(\R/\Z,\R^n)$ as we did
for the pure bending energy $\BendingEnergy$, and add the Banach spaces
$\CL \ceq W^{1+s,\rho}(\R/\Z,\R^n)$ with its open subset $\UL\ceq
W^{1+s,\rho}_\injreg(\R/\Z,\R^n)$ of immersed knots, where we choose
the fractional differentiability $s$ and integrability $\rho$ adapted to
the respective knot energy as in \cref{eq:individual-energy-spaces}, only
with the slightly more restricted parameter ranges \cref{eq:shrink-range}
assumed in \cref{thm:bending-knot-energy}.

\begin{refproof}[Proof of \cref{thm:bending-knot-energy}]
The bending energy $\BendingEnergy$
as well as all knot energies considered are non-negative. Moreover, it follows from
\cite[Proposition 4.8]{okabe-schrader_2023} that $\BendingEnergy\in C^1(\OL)$. 
Continuous
differentiability of the knot energies on the open subset $\UL$ of their 
energy space in the respective full parameter range given in \cref{eq:individual-energy-spaces} was shown in \cite[Theorem 1.4 \& subsequent remark]{blatt-reiter_2015a}
(see also \cite{wings_2018}) for $\TP^{(p,q)}$, in 
\cite[Theorem 3]{blatt-reiter_2015b}
for $\intM^{(p,q)}$, and in \cite[Theorem 5.1]{matt-etal_2023}
for $E^{\alpha,p}$. The slightly
restricted parameter ranges 
lead to a compact injection $\iota:\HL\hookrightarrow \CL$ with $\iota(\OL)\subset
\UL$ for any of the three knot energy families.
The conditions \ref{condE} and \ref{condP} follow from \cref{cor:basic-assumptions} for $s=1$. 
In addition, we have checked conditions \ref{condB*}, \ref{condS}, and \ref{condD} for $\BendingEnergy|_{\AL^2_\bz}$
in \cref{thm:PS_elastic}. 

As to condition \ref{condP*} required in \cref{cor:PS-total-energy},
observe that the explicit formula of the projection $\pr_\g^{\AL^2_\bz}(h)$
in \cref{eq:projection-ALsz} and \cref{eq:projection-ALs} based on the differential
$\D \varSigma $ of the logarithmic strain $\varSigma$ (see \cref{eq:log-strain-bounded-differential}) and its right
inverse $Y_\g$ given by \cref{eq:right-inverse-Y} does \emph{not} require the full regularity
$h\in H^2(\R/\Z,\R^n)$ to make sense. 
In fact, the starting point of the proof of \cref{prop:projection-ALs} uses \cite[Lemma 5.2]{knappmann-etal_2023} which is 
formulated for general fractional Sobolev spaces $W^{1+s,\rho}$, where $s\in (0,1),$
$\rho>1$ satisfy $s-\frac1{\rho} >0$. The energy spaces of all three knot energy families considered
here meet these requirements on the parameters $s$ and $\rho$; see \cref{eq:individual-energy-spaces}. So, we simply set 
$$
        \PR_\g \eta
        \ceq
        (\Id_\CL-Y_\g  \D \varSigma_\g) \, \eta - \pars[\big]{(\Id_\CL-Y_\g \, \D \varSigma_\g) \, \eta}(0)
        \quad\Fo\eta\in\CL,
$$
and copy the proof of \cref{prop:projection-ALs} line by line replacing
$h\in H^2(\R/\Z,\R^n)$ by $\eta\in\CL$ in each step, arriving at the estimate
corresponding to \cref{eq:projection-bound}, which confirms \ref{condP*}.

So, it remains to establish condition
\ref{condL*} for each of the knot energies, i.e., for $\FL\in\{\TP^{(p,q)},
\intM^{(p,q)},E^{\alpha,p}\}$. The key observation is the following: Any sequence
$(\g_k)_k\subset\AL^2_\bz$
with uniformly bounded knot energy, i.e., with $F\ceq\sup_{k\in\N}
\widetilde{\FL}(\g_k)
<\infty$, where $\widetilde{\FL}=\FL\circ\iota|_{\AL^2_\bz}$, satisfies a uniform bilipschitz
estimate depending only on the energy bound $F$ and the fixed  parameters
$s$ and $\rho$ determining the respective energy space as in 
\cref{eq:individual-energy-spaces}. In other words, there is a
constant $c=c(s,\rho,F)>0$ such that 
\begin{equation}\label{eq:uniform-bilip}
|\g_k(u)-\g_k(v)|\ge  c \,|u-v|_{\R/\Z}\quad\Foa u,v\in\R/\Z,\,k\in\N.
\end{equation}

For $\TP^{(p,q)}$ this was proven in \cite[Proposition 2.7]{blatt-reiter_2015a},
for $\intM^{(p,q)}$ in \cite[Proposition 2.1]{blatt-reiter_2015b}, and for
$E^{\alpha,p}$ in \cite[Proof of Theorem 1.1, pp. 21-22]{matt-etal_2023}.
The sequences $(\g_k)_k$
considered  in \ref{condL*} are also assumed to converge weakly in $\HL$
to some limit curve $\g_\infty\in\HL$.
By the compact embedding $\HL\hookrightarrow \BL=C^1(\R/\Z,\R^n)$ and the subsequence principle
one also has $\g_k\to\g_\infty$ in $C^1$ as $k\to\infty$ so that the bilipschitz
estimate \cref{eq:uniform-bilip} transfers to the limit curve $\g_\infty$; 
hence
$\g_\infty$ itself is a knot, satisfying $|\g_\infty'|=\lim_{k\to\infty}|\g_k'|=1$, and
$\g_\infty(0)=\lim_{k\to\infty}\g_k(0)=\bz$. To summarize,
$\g_\infty\in \AL^2_\bz$, so that \ref{condL*} is satisfied for
each $\FL\in\{\TP^{(p,q)},
\intM^{(p,q)},E^{\alpha,p}\}$. Therefore, \cref{cor:PS-total-energy} is applicable which implies the statement of the theorem.
\end{refproof}

\begin{refproof}[Proof of \Cref{thm:bending-moebius-energy}]
All the properties concerning  $\BendingEnergy$ as well as properties
\ref{condE} and \ref{condP} have been shown in the proof of
\Cref{thm:bending-knot-energy} for the same choice of spaces
$\HL \ceq H^2(\R/\Z,\R^n)$, $\BL \ceq C^1(\R/\Z,\R^n)$ and $\OL\ceq
H^2_\injreg(\R/\Z,\R^n)$. Also, \ref{condP*} follows in the same way if
we choose the slightly smaller Hilbert space $\CL \ceq H^{3/2 +\nu}(\R/\Z,\R^n)$
for some $\nu\in (0,1/2)$,
instead of the natural energy space $H^{3/2}(\R/\Z,\R^n)$ for the
M\"obius energy $E^{2,1}$ (see \cite[Theorem 1.1]{blatt_2012a}).
In \cite[Theorem 3.1]{reiter-schumacher_2021} it was shown that 
$E^{2,1}$ is of class $C^{1,1}_\loc$ on the open subset
$\UL \ceq H^{3/2 +\nu}_\injreg(\R/\Z,\R^n)$ of regular knots in $\CL$.
Finally, for the proof of \ref{condL*} we use the uniform
bilipschitz estimate for $E^{2,1}$ proved, e.g., in 
\cite[Lemma 3.1]{gilsbach-vdm_2018} 
together with the compact embedding $\HL\hookrightarrow\CL$ as
in the proof of \Cref{thm:bending-knot-energy}.
\end{refproof}

\begin{refproof}[Proof of \Cref{cor:minimizing-knots-total-energy}]
According to \Cref{thm:bending-knot-energy,thm:bending-moebius-energy} and their  proofs
all total 
energies $\TL$ under consideration are of class $C^1$ and
satisfy the (PS)-condition, and condition \ref{condL*}; hence
\ref{condL**} on the Hilbert submanifold $\AL^2_\bz$. Moreover,
all energies
are non-negative which verifies \ref{lower-bound}. 
Therefore,  \Cref{cor:critical+minimizers} implies
the result.
\end{refproof}

\begin{refproof}[Proof of \Cref{thm:gradient-flow-total-energy}]
In order to apply our abstract result on gradient flows,  
\Cref{thm:gradflow_subconvergence}, we notice 
first that the smooth arclength submanifold $\AL^2_\bz$ of the Hilbert
space $\HL \ceq H^2(\R/\Z,\R^n)$ inherits the Riemannian metric from the inner
product on $\HL$.
Moreover, the required  $C^{1,1}_\loc$-regularity
of all energies considered in \cref{eq:subfamilies}  holds true
 on the larger open subsets of regular knots in $\HL$,
 which  follows
from \cite[Theorem 3.1]{reiter-schumacher_2021} for the M\"obius energy
$E^{2,1}$, from \cite[Theorem 3.5]{steenebruegge_2023} for the tangent-point
energies $\TP^{(p,q)}$ listed in \cref{eq:subfamilies}, from
\cite[Theorem 2.1]{knappmann-etal_2023} for $\intM^{(p,2)}$, and
from \cite[Theorem 1.55]{knappmann_2020} for the classic integral Menger
curvature 
$\ML_p
=2^p\intM^{(p,p)} $
for $p>3$. The additional upper bounds on $p$ encoded
in  \cref{eq:subfamilies} and \cref{eq:shrink-range} are necessary for
the validity of the (PS)-condition established in 
\Cref{thm:bending-knot-energy} and \Cref{thm:bending-moebius-energy}. 
All energies are non-negative and therefore
satisfy condition \ref{lower-bound}, and we verified already condition
\ref{condL*} in the proof of \Cref{thm:bending-knot-energy} which 
implies condition \ref{condL**} on submanifolds of Hilbert spaces. Now,
\Cref{thm:gradflow_subconvergence} implies all statements on the gradient flow, including
the fact that the critical point in the limit as $t\to\infty$ is ambiently
isotopic to the initial knot $\g_0$ by means of 
\Cref{thm:components-knot-classes}. 
\end{refproof}

\subsection{Tangent-point energy}\label{sec:tanpoint}

The most difficult case is when the energy is a pure self-avoidance energy,
which naturally has a non-local structure so that now the leading order
terms are of fractional Sobolev regularity. We 
concentrate here on the tangent-point energies $\TP^{(p,2)}$ for
the range $p\in (4,5)$. 
As shown in \cite[Theorem 1.1]{blatt-reiter_2015a},
the underlying energy spaces are the Hilbert
spaces $H^{1+s}(\Circle,\AmbSpace)$ with $s \ceq (p-3)/2$.

A crucial tool will be the following splitting of the differential of $\TP^{(p,2)}$,
which is due to Blatt and Reiter, see \cite{blatt-reiter_2015a}.
Let $\gamma \in H^{1+s}_\injreg(\Circle,\AmbSpace)$ be parametrized by arclength
and let $h \in H^{1+s}(\Circle,\AmbSpace)$.
Then we have
\begin{equation}
    \D(\TP^{(p,2)})_\gamma h 
    =
    2 \, Q^{\TP}(\gamma,h) + R^{\TP}_\gamma h,
    \label{eq:splitting-TP}
\end{equation}
where $Q^{\TP}$ is a symmetric bilinear form (which will soon turn out to encode an elliptic operator), and $R^{\TP}_{\gamma}$ is a ``remainder term'', which is well-defined for all $ h\in \BL=C^1(\R/\Z,\R^n) $.
Blatt and Reiter also characterize $\QL^{\TP}$ as a Fourier multiplier.
More precisely, they show (see \cite[Proposition~4.1 and its proof]{blatt-reiter_2015a}):
\begin{proposition}\label{lem:fourier-Q}
        There are $\varrho_k \in \C$, $k \in \Z$ such that 
        \[
                Q^{\TP}(f,h)
                =
                \sum_{k \in \Z} \varrho_k \inner{ \hat f_k, \hat h_k }_{\C^\AmbDim}
                \quad 
                \text{for all $f,h \in H^{\frac{p-1}{2}}(\Circle,\AmbSpace)$}
                ,
        \]      
        where $\hat{f}_k,\hat{h}_k \in \C^\AmbDim$ are the Fourier coefficients of $f$ and $h$, respectively.
        Moreover, there is a constant $c_{\TP} > 0$ such that 
        \[
                \varrho_k = c_{\TP} \abs{k}^{p-1} + o(\abs{k}^{p-1})
                \quad 
                \text{as $\left| k \right|\to \infty$,}
        \]
        where $ 0=\rho_0<\rho_{k}=\rho_{-k} $ for all $ k\in \N $, and $ \rho_{k} $ is non-decreasing in $ \left|k \right| $. 
\end{proposition}
The asymptotics of $(\varrho_k)_{k \in \N}$ imply that $f \mapsto Q^{\TP}(f,\cdot)$ is an elliptic operator of order $p-1$ in the following sense:

\begin{lemma}[$ Q^{\TP} $ elliptic]\label{lem:EllipticRegularityQTP}
    Let $ f\in L^2(\R/\Z,\R^n) $. Then the following holds. 
    \begin{enumerate}
        \item If $ f\in H^{t}\left( \R/\Z,\R^n \right) $ for some $ t\in \R $, then $ Q^{\TP}(f,\cdot) $ extends to a linear functional in $ H^{t-(p-1)}\left( \R/\Z,\R^n \right)\equiv(H^{(p-1)-t}\left( \R/\Z,\R^n \right)  )^* $, satisfying 
        \begin{equation}\label{eq:46eq1}
            \left\|Q^{\TP}(f,\cdot) \right\|_{H^{t-(p-1)}}\le C_{Q}\left\|f \right\|_{H^t}
        \end{equation}
        for some constant $ C_Q $ independent of f. 
        \item If $ Q^{\TP}(f,\cdot) \in H^{t-(p-1)}(\R/\Z,\R^n)$ for some $ t\in \R $, then $ f\in H^{t}(\R/\Z,\R^n) $. 
        \item Let $ s:=\frac{p-1}{2}-1 $. Then there are constants $ 0<C_1\le C_2<\infty $ such that 
        \begin{equation}\label{eq:46eq2}
            C_1\left\|f \right\|^2_{H^{1+s}}\le \left\|f \right\|_{L^2}^2+Q^{TP}(f,f)\le C_2\left\|f \right\|_{H^{1+s}}^2\,\Foa f\in H^{1+s}(\R/\Z,\R^n).
        \end{equation}
    \end{enumerate}
\end{lemma}

\begin{proof}
    \cref{lem:fourier-Q} implies that there is some $ k_0\in\N $ such that 
    \begin{equation}\label{eq:Qineq}\textstyle
        \frac{c_{\TP}}{2}\left|k \right|^{p-1}\le \rho_k\le 2 c_{\TP}\left|k \right|^{p-1} \quad \Foa \left|k \right|> k_0.
    \end{equation}
    We estimate for $ f,h\in C^\infty(\R/\Z,\R^n) $ by means of \cref{lem:fourier-Q} and \cref{eq:Qineq}: 
    \begin{equation}\label{eq:estimateQelip}
    \begin{aligned}
        Q^{\TP}(f,h)
         &\textstyle\le \sum_{k\in\Z\setminus\{0\}}\rho_{k}\left|\hat{f}_{k} \right|\left|\hat{h}_{k} \right|
        \le\rho_{k_0}\sum_{0<\left|k \right|\le k_0}\left|\hat{f}_{k} \right|\left|\hat{h}_{k} \right|+2c_{\TP}\sum_{\left|k \right|>k_0}\left|k \right|^{t}\left|\hat{f}_{k} \right|\left|\hat{h}_{k} \right|\left|k \right|^{(p-1)-t}\\
        &\textstyle\le C_Q \sum_{k\in\Z\setminus\{0\}}( 1+\left|k \right|^t\left|k \right|^{p-1-t} )\left|\hat{f}_{k} \right|\left|\hat{h}_{k} \right|\\
        &\textstyle\le C_Q \sum_{k\in\Z\setminus\{0\}}\left( 1+\left|k \right|^t \right)\left|\hat{f}_{k} \right|( 1+\left|k \right|^{p-1-t} )\left|\hat{h}_{k} \right|
        \le C_Q \left\|f \right\|_{H^t}\left\|h \right\|_{H^{p-1-t}},
    \end{aligned}
    \end{equation}
    where we have set $ C_Q\ceq\max\{\rho_{k_0},c_{\TP}\} $.
     By density of $ C^\infty(\R/\Z,\R^n) $ in $ H^{\sigma}(\R/\Z,\R^n) $ for all $ \sigma\in \R $ (see, e.g. \cite[Statement 4 in Section 2.3]{AgranovichSobolevSpaces2015}) this estimate extends to all $ f\in H^{t}(\R/\Z,\R^n) $ and $ h\in H^{p-1-t}(\R/\Z,\R^n) $. 
    Hence, $ Q(f,\cdot)\in H^{t-(p-1)} \left( \R/\Z,\R^n \right)$ if $ f\in H^t(\R/\Z,\R^n) $ 
    satisfying \cref{eq:46eq1}. Choosing $ t=\frac{p-1}{2} $ and $ h\ceq f $ in \cref{eq:estimateQelip} yields the upper bound on $ \left\|f \right\|_{L^2} +Q^{\TP}(f,f)$ in \cref{eq:46eq2}. 
    Similarly, as in \cref{eq:estimateQelip} we use \cref{lem:fourier-Q} to estimate  
   \begin{align*}
     2Q^{\TP}(f,f)&\textstyle\ge \frac{2\rho_1}{k_0^{2+2s}}\sum_{0<\left|k \right|\le k_0}\left| k \right|^{2+2s} \left|\hat{f}_{k} \right|^2+c_{\TP}\sum_{\left|k \right|>k_0}\left|k \right|^{2+2s}\left|\hat{f}_{k} \right|^2
        \ge\rho_{*}\sum_{k\in\Z}\left|k \right|^{2+2s}\left|\hat{f}_{k} \right|^2
   \end{align*}
   for $ \rho_{*}\ceq \min\{\frac{2\rho_{1}}{k_0^{2+2s}},c_{\TP}\} >0$. Therefore, the lower bound in \cref{eq:46eq2} follows from the fact that 
   \begin{align*}\textstyle
    \left\|f \right\|_{H^{1+s}}^2=\sum_{k\in\Z}(1+\left|k \right|^2)^{1+s}\left|\hat{f}_{k} \right|^2 \le 2^{1+s}\left( \sum_{k\in\Z} \left|f_k \right|^2+\sum_{k\in \Z}\left|k \right|^{2+2s}\left|\hat{f}_{k} \right|^2\right),
   \end{align*}
   if we set $ C_1=\frac{1}{2^{2+s}}\min\{\rho_{*},2\} $. 

   It remains to show the second statement, which does \emph{not} come with an a priori estimate for $ \left\|f \right\|_{H^{t}} $ in terms of $ \left\|Q^{TP}(f,\cdot) \right\|_{H^{t-(p-1)}} $. Considering only those $ h\in C^\infty(\R/\Z,\R^n) $ with $ \hat{h}_{k}=0 $ for $ \left|k \right|=0,\dots,k_0 $, we infer from \cref{eq:Qineq} 
   \begin{align*}\textstyle
    \frac{c_{\TP}}{2}\sum_{\left|k \right|>0}\left\langle \hat{f}_{k}\left|k \right|^t,\left|k \right|^{p-1-t}\hat{h}_{k} \right\rangle_{\C^n} \le \left\|Q^{\TP}(f,\cdot) \right\|_{H^{t-(p-1)}}\left\|h \right\|_{H^{(p-1)-t}}
   \end{align*}
   for all such $ h $. By density of smooth functions with the first $ 2k_{0}+1 $ Fourier coefficients vanishing we obtain
   \begin{align*}\textstyle
    &\|( \delta_{k,k_0}\hat{f}_{k}\left|k \right|^t )_{k}\|_{\ell^2}\\
    &\textstyle=\sup\left\{\left|\sum_{\left|k \right|>k_0}\left\langle \hat{f}_{k}\left|k \right|^t,\left|k \right|^{p-1-t}\hat{h}_{k} \right\rangle_{\C^n}  \right|: h\in C^\infty,\hat{h}_{k}=0\text{ for } \left|k \right|\le k_0,\left\|h \right\|_{H^{(p-1)-t}}\le 1\right\}\\
    &\le \left\|Q^{\TP}(f,\cdot) \right\|_{H^{t-(p-1)}},
   \end{align*}
   where $ \delta_{k,k_0}$  equals $1 $ if $ \left|k \right|>k_0 $ and $ 0 $ if $ \left|k \right|\le k_0 $. This implies that $ (\hat{f}_{k}\left|k \right|^{t})_{k}\in \ell^2 $, which together with the assumption that $ f\in L^2(\R/\Z,\R^n) $ shows that $ ((1+k^2)^{\frac{t}{2}}\hat{f}_{k})_{k}\in \ell^2 $, or equivalently $ f\in H^{t}(\R/\Z,\R^n).$
\end{proof}

We are now going to verify the assumption of our
abstract result, \Cref{thm:PS-abstract},
guaranteeing the (PS)-condition for the tangent-point energies.

\begin{refproof}[Proof of \Cref{thm:PS-TP}]
Fix $p\in (4,5)$ and $s \ceq s(p) = (p-3)/ 2$. 
Throughout the proof
we abbreviate the energy by $\TP \ceq \TP^{(p,2)}$ and consider the open 
subset $\OL \ceq H^{1+s}_\injreg(\R/\Z,\R^n)$ of the Hilbert space
$\HL \ceq H^{1+s}(\R/\Z,\R^n)$, which compactly embeds into the
Banach space $\BL \ceq C^1(\R/\Z,\R^n)$. Clearly, $\TP$ is non-negative,
and it is of class $C^1(\OL)$ due to \cite[Theorem 1.4 \& subsequent
remark]{blatt-reiter_2015a} and \cite{wings_2018}. Since $s\in (1/2,1)$
and $\bz\in\R^n$ are fixed, we can apply
\Cref{cor:basic-assumptions} to deduce condition \ref{condE} for the 
arclength submanifold $\AL^{1+s}_\bz$. Moreover, the linear projection
\[
        \smash{\pr_\g^{\AL^{1+s}_\bz} \colon T_\g\OL\to T_\g\AL^{1+s}_\bz}
\]
has property \ref{condP}.

Towards condition \ref{condL} we take a (PS)-sequence $(\g_m)_m\subset\AL^{1+s}_\bz$
for $\TP|_{\AL^{1+s}_\bz}$, converging strongly
in the Hilbert space $\HL$ to some curve $\g_\infty\in
\HL$. Since $T \ceq \sup_{m\in\N}\TP(\g_m)<\infty$ we can apply the uniform
bilipschitz estimate \cite[Proposition 2.7]{blatt-reiter_2015a} to deduce
the existence of some constant $b=b(T,p)>0$ such that
\begin{equation}\label{eq:bilip-TP}
|\g_m(u)-\g_m(v)|\ge b|u-v|_{\R/\Z}\quad\Foa u,v\in\R/\Z,\,m\in\N.
\end{equation}
By the compact embedding $\HL\hookrightarrow\BL$ we find a $C^1$-convergent
subsequence $\g_{m_k}\to\g_\infty$ as $k\to\infty$, so that \cref{eq:bilip-TP}
is conserved in the limit; hence $\g_\infty$ is a knot, satisfying in
addition
$$
%\textstyle
|\g_\infty'(u)|=\lim_{k\to\infty}|\g_{m_k}'(u)|=1\Foa u\in\R/\Z
\quad\AND\quad \g_\infty(0)=\lim_{k\to\infty}\g_{m_k}(0)=\bz.
$$
Therefore, $\g_\infty\in\AL^{1+s}_\bz$, which proves condition \ref{condL}.

The uniform energy bound $T$ implies by means of 
\cite[Proposition 2.5]{blatt-reiter_2015a} the existence of a constant 
$C=C(T,p)$  such that $\sup_{m\in\N}\lfloor\g_m'\rfloor_{s,2}\le C$,
where $\lfloor\cdot\rfloor_{s,2}$ denotes the seminorm
characterizing the fractional
Sobolev space $H^{1+s}$ as subspace of
the standard Sobolev space~$H^1$. 
In addition,
$$
|\g_m(u)|\le |\g_m(u)-\g_m(0)|+|\bz|\le 1+|\bz|\quad\Foa u\in\R/\Z,\,m\in\N,
$$
since $|\g_m'|=1$ on $\R/\Z$, so that $\sup_{m\in\N}\|\g_m\|_{H^{1+s}}\le
2+|\bz| + C,$ which establishes condition \ref{condB}.

The splitting \cref{eq:splitting-TP} of $\D (\TP)_\g h$ leads us to 
define as quadratic term $\QL$ in condition \ref{condS}
\begin{equation}\label{eq:def-QL}
        \QL(\g,h) \ceq 2 \, Q^{\TP}(\g,h)+2\langle \g,h\rangle_{L^2}.
\end{equation}

By definition, $ \QL $ is a continuous bilinear form and satisfies \ref{condS1} with $\mu(\cdot)\equiv 1$. 
Because of \cref{eq:46eq2} in \cref{lem:EllipticRegularityQTP} the bilinear form $ \QL $ satisfies \ref{condS2}
with coercivity constant $c \ceq 2C_1$.

It remains to verify condition \ref{condS3} for the remainder term $ R^{\TP}_{\gamma}h-2 \left\langle \gamma,h \right\rangle_{L^2}  $
to obtain the full
splitting condition \ref{condS}. The $L^2$-inner product is bounded on bounded sets (and hence on  \Econtrolled{\TPp} sets) since
\begin{equation}\label{eq:L2-term}
\langle\g,h\rangle_{L^2}\le\norm{\g}_{L^2} \norm{h}_{L^2} \le \norm{\g}_\HL \norm{h}_\BL \quad \Fo \gamma\in \AL^{1+s}_{\bz},h\in \BL.
\end{equation}
According to \cite[Lemma 4.2]{blatt-reiter_2015a} the term $R^{\TP}_{\gamma}h$ is the finite sum
of expressions, each of which can be bounded from above by
\begin{equation}\label{eq:bounding-R-TP}
\textstyle
C\int_{\R/\Z}\int_\frac{1}{2}^\frac{1}{2}\int_0^1\int_0^1\frac{|\g'(u+s_1w)-\g'(u+s_2w)|^2}{|w|^{p-2}}\d \theta_1\d\theta_2\d w\d u\cdot \|h\|_{C^1},
\end{equation}
where $s_i\in\{0,\theta_i\}$ for $i=1,2$, and
$C=C(\TP(\g)) $ is a constant that depends non-decreasingly on the energy value $ \TP(\gamma) $ through the bilipschitz constant $ b=b(\TP(\gamma),p)>0 $ again by 
\cite[Proposition 2.7]{blatt-reiter_2015a}.
By inserting $\pm\g'(u)$ in the numerator of the multiple integral in 
\cref{eq:bounding-R-TP} if
$s_i\not= 0$ for $i=1,2$, we bound this integral by means of the triangle
inequality from above by two triple integrals of the form treated in 
\Cref{lem:triple-integrals} in \Cref{sec:arclength-reparam}, 
each of which
are bounded from above by the seminorm $\lfloor\g'\rfloor_{s,2}^2$.
So, $ \RL $ is bounded on \Econtrolled{\TPp} sets, and thus satisfies \ref{condS3}. 

Finally, we are left with verifying condition \ref{condD}. Combining \cref{eq:projection-ALs,eq:projection-ALsz,eq:log-strain-differential} we can write 
\begin{align*}
    ( \Id_{\HL}-\pr^{\AL^{1+s}_{\bz}}_{\gamma})h \overset{\eqref{eq:projection-ALs},\eqref{eq:projection-ALsz}}{=}Y_{\gamma}\D\Sigma_{\gamma}h+\pr^{\AL^{1+s}}_{\gamma}h(0)\overset{\eqref{eq:log-strain-differential}}{=}Y_\gamma \left\langle \gamma',h' \right\rangle +\pr^{\AL^{1+s}}_{\gamma}h(0).
\end{align*}
Since $ \D(\TPp)_{\gamma}\tilde{h}\, $ contains\footnote{The explicit formula for the differential of $ \TPp $ is restated in \cref{eq:DTP_formula} in the proof of \cref{lem:RegularityTPRemainderII} below.} only differences and derivatives of a general test function $ \tilde{h} $, {(see \cite[(1.11) in Theorem 1.4]{blatt-reiter_2015a} specified to the present case $ q=2 $)} we only need to bound 
$|\D(\TPp)_{\gamma}Y_\gamma \left\langle \gamma',h' \right\rangle |$
from above by $ C\left\|h \right\|_{\BL} $, where the constant $ C=C(\gamma) $ remains controlled on \Econtrolled{\TPp} subsets of $ \AL_{\bz}^{1+s} $. This follows from \cref{eq:regremainderTPcondD} of \cref{lem:RegularityTPRemainderII} below, since $ H^{3/2+\varepsilon}(\R/\Z,\R^n) $ embeds into $ C^1(\R/Z,\R^n) $ for all $ \varepsilon>0 $. Thus, \ref{condD} holds true for $ S_{\gamma}\ceq\D(\TPp)_{\gamma}Y_\gamma \left\langle \gamma',(\cdot)' \right\rangle  $.
\end{refproof}

\begin{refproof}[Proof of \Cref{cor:minimizing-knots-TP} and \Cref{thm:gradient-flow-TP}.]
It was shown in \cite[Theorem 3.5]{steenebruegge_2023} that the 
tangent-point energy $\TP^{(p,2)}$ for any $p\in (4,5)$ is of class
$C^{1,1}_\loc$ on the open subset $\OL \ceq H^{1+s(p)}_\injreg(\R/\Z,\R^n)$
of the Hilbert space $\HL \ceq H^{1+s(p)}(\R/\Z,\R^n)$
for $s(p)$ as in \cref{eq:s-TP}. In fact $ \TPp $ is even smooth on its natural energy space, see \cite[Theorem 4.6]{doehrer-etal_2025}. The (PS)-condition holds according
to \Cref{thm:PS-TP}. Of course, $\TP^{(p,2)}$ is non-negative; hence
satisfies \ref{lower-bound}. For the verification of
condition \ref{condL}  in the proof of
\Cref{thm:PS-TP}  we only needed the uniform boundedness of
the energy along the sequence, which is assumed to converge in $\HL$,
so that \ref{condL**} follows. So, we may apply 
\Cref{thm:gradflow_subconvergence}
and
\Cref{koro:minimizer_existence}
to conclude.
\end{refproof}

\begin{lemma}\label{lem:RegularityTPRemainderII}
    Let $ \gamma\in \AL^{t} $ for some $ t\ge {(p-1)}/{2} $ and $ p\in(4,5) $. Then the mapping $ h\mapsto\D(\TPp)_{\gamma}Y_\gamma \left\langle \gamma',h' \right\rangle  $ is of class $ H^{t-1-\frac{p}{2}-\varepsilon}(\R/\Z,\R^n) $ for every $ \varepsilon>0$. 
    In particular, $ t={(p-1)}/{2} $, then for each $ \varepsilon>0 $ there is a constant $ C_{\varepsilon} =C_{\varepsilon}(\left\|\gamma \right\|_{H^{{(p-1)}/{2}}},\TPp(\gamma))$ depending non-decreasingly on the $ H^{{(p-1)}/{2}} $-norm and on the energy value of $ \gamma $, such that 
    \begin{equation}\label{eq:regremainderTPcondD}
    \begin{aligned}
        |\D(\TPp)_{\gamma}Y_\gamma \left\langle \gamma',h' \right\rangle  |\le C_{\varepsilon}\left\|h \right\|_{H^{\frac{3}{2}+\varepsilon}}\quad \Foa h\in H^{\frac{3}{2}+\varepsilon}(\R/\Z,\R^n).
    \end{aligned}
    \end{equation}
\end{lemma}
\begin{proof}
    Abbreviate $ \TP\ceq\TPp $ for some fixed $ p\in(4,5) $. It suffices to prove that\\
    $ \D(\TP)_{\gamma}Y_\gamma \left\langle \gamma',h' \right\rangle  $ is a finite sum of expressions of type\footnote{We follow the convention of \cite{blatt-reiter_2015a} that $ \bigoasterisk $ denotes any kind of product.}
    \begin{equation}\label{eq:Rgpsi-structure}
        \textstyle
        \int_{\R/\Z}\int_{-\frac{1}{2}}^{\frac{1}{2}}\idotsint_{[0,1]^K}g^{(p)}(u,w,s_1,\dots,s_{K})\bigoasterisk
        \psi(u+s_Kw)\d \theta_1\cdots\d\theta_K\d w\d u,
    \end{equation}
    where 
    \begin{equation}\label{eq:g-formstructure}\textstyle
        g^{(p)}(u,w,s_1,\dots,s_{K})\ceq\GL\left( \left|\frac{\triangle \gamma}{w} \right| \right)\frac{|\gamma'(u+s_1w)-\gamma'(u+s_2w)|^2}{\left|w \right|^{p-2}}\bigoasterisk\limits_{i=1}^{K-1}\gamma'(u+s_iw),
    \end{equation}
    $ \GL^{(p)} $ is analytic on $ (0,\infty) $, $ s_i\in\{0,\theta_i\} $ for $ i=1,\dots,K $, $ K\ge3 $. 
    Here, either $ \psi=h' $, or $ \psi\equiv V_{\gamma}f $, for $ f\ceq \left\langle \gamma',h' \right\rangle  $, where $ V_\gamma f $ is the vector defining the right-inverse $ Y_\gamma $ in \cref{eq:right-inverse-Y}:
    \begin{equation}\label{eq:right-inverse-formula-again}
    \begin{aligned}
          Y_\g f(x)&=\textstyle
	    \int_0^x  f(u) \, \g'(u) \d u+ \int_0^x \P_{\g'(u)}^\perp V_{\gamma} f \d u\\
        &\qec\textstyle \Upsilon^{1}f(x)+\Upsilon^2f(x)\quad \Foa f\in H^{s}(\R/\Z,\R^n).
    \end{aligned}
    \end{equation}
    Indeed, once $ \D(\TP)_\gamma Y_\gamma \left\langle \gamma',h' \right\rangle  $ is shown to have that specific structure, one can check that the proof of \cite[Proposition 4.3]{blatt-reiter_2015a} can be adopted without any changes to show that each expression of the form \cref{eq:Rgpsi-structure} as a function of $ h $ in case $ \psi=h' $, is of class $ H^{t-1-\frac{p}{2}-\varepsilon}(\R/\Z,\R^n) $ for every $ \varepsilon>0 $. If, however, $ \psi=V_\gamma f $, then recall that $ \gamma $ has bilipschitz constant $b_{\gamma}>0$ depending merely on the energy value $ \TP(\gamma) $, by \cite[Proposition 2.7]{blatt-reiter_2015a} (cf \cref{eq:bilip-TP}), so that we simply estimate by means of $ \left|\gamma' \right|=1 $, 
    \begin{align*}\textstyle
        |g^{(p)}(u,w,s_1,\dots,s_K)|\le C(\TP(\gamma))\frac{|\gamma'(u+s_1w)-\gamma'(u+s_2w)|^2}{\left|w \right|^{p-2}}\,\,\Foa u\in\R/\Z,\, \left|w \right|\le \frac{1}{2},
    \end{align*}
    and for $  s_1,\dots,s_K $.
    Inserting $ \pm \gamma'(u) $ in the numerator on the right-hand side allows us to apply \cref{lem:triple-integrals} to conclude 
\begin{equation}\label{eq:Vf-estimate}
\begin{aligned}\textstyle
    &\textstyle|\int_{\R/\Z}\int_{-\frac{1}{2}}^{\frac{1}{2}}\idotsint_{[0,1]^K}g^{(p)}(u,w,s_1,\dots,s_K)\bigoasterisk V_\gamma f \theta_1\cdots\d\theta_K\d w\d u|
    \le 4C(\TP(\gamma))\lfloor\g'\rfloor_{s,2}^2\left|V_\gamma f \right|.
\end{aligned}
\end{equation}
Apply \cref{lem:DEbYisLowerOrderIIII} \textrm{(i)} for $ f=\left\langle \gamma',h' \right\rangle  $ and \cref{lem:ProductRuleAsDistributionII} for $ \tau=\sigma(t-1,\varepsilon)=t-1 $ to conclude that this is bounded by $ \tilde{C}\left\|h \right\|_{H^{1-t}} $ with a constant $ \tilde{C} $ depending non-decreasingly on $ \left\|\gamma \right\|_{H^t} $ and $ \TP(\gamma) $. Since $ H^{t-1}\subset H^{t-1-\frac{p}{2}-\varepsilon} (\R/\Z,\R^n)$ for all $ \varepsilon>0 $, this proves the claim. In addition, combining \cref{eq:Vf-estimate} for $ t=(p-1)/2$, with an analogous estimate in the case $ \psi=h' $ leads to 
\cref{eq:regremainderTPcondD}. 

So, it remains to establish the structure \cref{eq:Rgpsi-structure} for all summands of $ \D(\TP)_{\gamma}Y_\gamma \left\langle \gamma',h' \right\rangle  $. To that extent we specify the formula for the first variation $ \delta \TPpq[](\gamma,h) $ given in \cite[Theorem 1.4]{blatt-reiter_2015a} to the present case $ q=2,\,p\in(4,5) $ to obtain for general $ h\in H^{1+s}(\R/\Z,\R^n) $, $ s=(p-3)/2 $, (using the abbreviation $ \triangle F=F(u+w)-F(u) $, for $ F:\R/\Z\to\R^n $)
\begin{equation}
\begin{aligned}
    \label{eq:DTP_formula}
        &\textstyle\D(\TP)_{\gamma}h=\int_{\R/\Z}\int_{-\frac{1}{2}}^{\frac{1}{2}}
		2\frac{\langle \P^\perp_{\gamma'(u)}(\triangle\gamma),\triangle h -\left\langle \triangle \gamma,\gamma'(u) \right\rangle h'(u) \rangle }{\left|\triangle \gamma \right|^p}\dd u \dd w\\
		&\textstyle\hspace{12pt}+\int_{\R/\Z}\int_{-\frac{1}{2}}^{\frac{1}{2}}|\P^\perp_{\gamma'(u)}(\triangle \gamma) |^2\big(-p\frac{\left\langle \triangle \gamma,\triangle h \right\rangle}{\left|\triangle \gamma \right|^{p+2}}
		+\frac{(\left\langle \gamma'(u),h'(u) \right\rangle+\left\langle \gamma'(u+w),h'(u+w) \right\rangle  )}{\left|\triangle \gamma \right|^{p}}\big)\dd u \dd w\\
        &\textstyle\hspace{12pt}\qec\widetilde{Q}_{\gamma}h+\widetilde{R}_3h+\widetilde{R}_4h.
\end{aligned}
\end{equation}
Notice that for the last two summands $ \widetilde{R}_3h=R_3^{(p)} $ and $ \widetilde{R}_4h=R_4^{(p)} $, where $ R_i^{(p)} $ are defined in \cite[(4.3)]{blatt-reiter_2015a}. In the present context we have to replace $ h $ by the more specific function $ Y_\gamma f $ given in \cref{eq:right-inverse-formula-again} for $ f=\left\langle \gamma',h' \right\rangle  $, satisfying \cref{eq:derivative-right-inverse}
and therefore \cref{eq:first_derivative_rightInverse}. Consequently, $ \tilde{R}_{4}Y_{\gamma}\left\langle \gamma',h' \right\rangle=\tilde{R}_{4}h =R_4^{(p)} $, and the structure of that term has been shown in \cite[Lemma 4.2]{blatt-reiter_2015a} to be of the form \cref{eq:Rgpsi-structure} for $ \psi=h' $ if one artificially adds the factors $ \gamma'(u+s_1w) $ and $ \gamma'(u+s_2w) $ to the last product in \cref{eq:g-formstructure}. 

The term $ \widetilde{R}_{3}\gamma $ has the {same form \cref{eq:Rgpsi-structure} for $ \psi=h' $,} which translates via \cref{eq:derivative-right-inverse} directly to one summand of the form \cref{eq:Rgpsi-structure} with $ \psi:=f\gamma'=\left\langle \gamma',h' \right\rangle \gamma' $, and a second summand of the form \cref{eq:Rgpsi-structure} with $ \psi\ceq \P^{\perp}_{\gamma'}V_\gamma f=V_\gamma f-\langle V_\gamma f,\gamma' \rangle \gamma'$. In both summands we can absorb the $ \gamma' $-terms of the respective $ \psi $ into the extended product to obtain $ \bigoasterisk_{i=1}^{K+1}\gamma'(u+s_i w) $, with $ s_K=s_{K+1} $, and end up with terms that are again of the form \cref{eq:Rgpsi-structure} with either $ \psi=h' $ or $ \psi=V_{\gamma}f$. As shown, the first term is treated as in \cite[Proposition 4.3]{blatt-reiter_2015a}, and the second by direct estimates as in \cref{eq:Vf-estimate} above. 

Finally, we need to investigate $ \widetilde{Q}_{\gamma}h $ in \cref{eq:DTP_formula} for $ h\ceq Y_{\gamma}f $, $ f=\left\langle \gamma',h' \right\rangle  $. By means of 
\[\textstyle\triangle\gamma=w\int_{0}^{1}\gamma'(t_1)\d \theta_1, \text{ and }\triangle(Y_\gamma f)=w\int_{0}^{1}(Y_\gamma f)'(t_2)\d \theta_2,\]
for $ t_i=u+\theta_i w,\,i=1,2$, we rewrite the numerator of $ \widetilde{Q}_{\gamma}Y_\gamma f $  as 
\begin{equation}\label{eq:ABsplitting}
\begin{aligned}\textstyle
            2\inner[\big]{ 
                \P^\perp_{\gamma'(u)} \Delta\gamma, \Delta (Y f)
                -
                \inner{ \Delta \gamma,\gamma'(u) } \pars{Y f }'(u) 
            }
            =
            2w^2 \! \int_0^1\int_0^1 
                \!\pars{\mathbf{A} + \mathbf{B}}
            \d \theta_1 \d \theta_2,
\end{aligned}
\end{equation}
where - according to the splitting \cref{eq:right-inverse-formula-again} - 
\begin{equation}\label{eq:splitting-A-part}
\begin{aligned}
    \mathbf{A} 
            &\ceq 
            \inner[\big]{
                 \P^\perp_{\gamma'(u)} \gamma'(t_1)
                 , 
                 (\Upsilon^1f)'(t_2)
                 {-
                 \inner{ \gamma'(t_2),\gamma'(u) } (\Upsilon^1f)'(u)}
            }\\
            &\overset{\cref{eq:right-inverse-formula-again}}{=}\inner[\big]{
                 \P^\perp_{\gamma'(u)} \gamma'(t_1)
                 , 
                 \gamma'(t_2)
            }f(t_2),
\end{aligned}
\end{equation}
and, with the abbreviation $ V\ceq V_{\gamma}f $,
\begin{equation}\label{eq:splitting-B-part}
\begin{aligned}
    \mathbf{B} 
            &\ceq 
            \inner[\big]{
                 \P^\perp_{\gamma'(u)} \gamma'(t_1)
                 , 
                 (\Upsilon^2f)'(t_2)
                 -
                 \inner{ \gamma'(t_2),\gamma'(u) } (\Upsilon^2 f)'(u) 
            }
            \\
             &=
            \inner[\big]{ 
                P^\perp_{\gamma'(u)} \gamma'(t_1)
                ,
                V 
                - \gamma'(t_2) \inner{\gamma'(t_2),V} 
                - \inner{\gamma'(t_2),\gamma'(u)} (V-\left\langle V,\gamma'(u) \right\rangle \gamma'(u))
            }\\
            &\overset{\cref{eq:right-inverse-formula-again}}{=}
            \inner[\big]{
                 \P^\perp_{\gamma'(u)} \gamma'(t_1)
                 , V
            }(1-\left\langle \gamma'(t_2),\gamma'(u) \right\rangle )-
            \inner[\big]{
                 \P^\perp_{\gamma'(u)} \gamma'(t_1)
                 , 
                \gamma'(t_2)               
            }
            \left\langle V,\gamma'(t_2) \right\rangle.
\end{aligned}
\end{equation}
We expand the prefactor of $ f(t_2) $ in \cref{eq:splitting-A-part}, which also appears in \cref{eq:splitting-B-part}, by means of the identity $ 1-\left\langle a,b \right\rangle=\frac{1}{2}\left|a-b \right|^2  $ for $ a,b\in \mathbb{S}^{n-1} $, as 
\begin{equation}\label{eq:vdMscan-C}
\begin{aligned}
    &\left\langle \gamma'(t_1)-\left\langle \gamma'(t_1),\gamma'(u) \right\rangle\gamma'(u),\gamma'(t_2)  \right\rangle \\
    =&\left[\left\langle \gamma'(t_1),\gamma'(t_2) \right\rangle -1
    +1-\left\langle \gamma'(t_1),\gamma'(u) \right\rangle 
    +\left\langle \gamma'(t_1),\gamma'(u) \right\rangle (1-\left\langle \gamma'(u),\gamma'(t_2) \right\rangle )
    \right]\\
    =&\textstyle\frac{1}{2}\left[-\left|\gamma'(t_1)-\gamma'(t_2) \right|^2+\left|\gamma'(t_1)-\gamma'(u) \right|^2+\left\langle \gamma'(t_1),\gamma'(u) \right\rangle\left|\gamma'(u)-\gamma'(t_2) \right|^2 \right],
\end{aligned}
\end{equation}
whereas the first summand of $ \mathbf{B} $ equals 
\begin{equation}\label{eq:vdMscan-D}
\begin{aligned}\textstyle
    \frac12 \left\langle \gamma'(t_1)-\left\langle \gamma'(t_1),\gamma'(u) \right\rangle\gamma'(u),V  \right\rangle \left|\gamma'(t_2)-\gamma'(u) \right|^2.
\end{aligned}
\end{equation}
Combining \eqref{eq:ABsplitting}--\eqref{eq:splitting-B-part} with \cref{eq:vdMscan-C,eq:vdMscan-D}, every summand constituting $ \widetilde{Q}_{\gamma}Y_\gamma f $ for $ f=\left\langle \gamma',h' \right\rangle  $ has indeed the form \cref{eq:Rgpsi-structure,eq:g-formstructure}, for $ \psi=h' $ or $ \psi =V_\gamma f $.
\end{proof}  

%% file: Smoothness.tex
\section{Regularity of critical points}\label{sec:Smoothness}
This section is devoted to the proof of \cref{thm:smoothness}, which states that critical points of linear combinations 
\[\TL_{\nu\theta}\ceq\nu\BendingEnergy+\theta\TPp, \quad \nu,\theta\ge 0,\,\,\nu+\theta>0,\]
of bending and tangent-point energies are smooth, including the cases where only one of those energies appears. The case $ \nu>0,\theta=0 $ will be covered by \Cref{thm:elastic-smoothness-new}. For mixed energies, $ \nu,\theta>0 $, the result follows from \Cref{thm:smoothness-T-critical-new}. The pure tangent-point energy is treated in \Cref{thm:smoothness-TP-new}, which covers the case $ \nu=0,\theta>0 $, and thus completes the proof of \Cref{thm:smoothness}. 

Before proceeding to these three scenarios, let us make some general comments. If $ \theta=0 $ we work on the submanifold 
\[\NL^2\subset\OL\ceq H^{2}_{\reg}(\R/\Z,\R^n)\subset H^{2}(\R/\Z,\R^n)\]
of immersed arclength parametrized curves. In all other cases ($ \theta>0 $) we restrict to the submanifolds 
\[\AL^{1+s}\subset \OL\ceq H^{1+s}_{\injreg}(\R/\Z,\R^n)\subset H^{1+s}(\R/\Z,\R^n)\]
of arclength parametrized knots, where $ s=1 $ whenever $ \nu>0 $, and $ s=s(p)=(p-3)/2 $ for $ p\in (4,5) $ if $ \nu=0 $ and $ \theta>0 $. 
For a fixed point $ \bz\in \R^n $ and $ \ML\in \{\NL^{1+s},\AL^{1+s}\} $ we may also consider the submanifolds $ \ML_{\bz}\ceq\{\gamma\in \ML: \gamma(0)=\bz\} $ without changing the regularity statements. Indeed, \cref{cor:critical-points_wo_fixed_point} implies that every $ \EL|_{\ML_{\bz}}  $-critical point for any $ \EL\in \{\TL_{\nu\theta}:\nu,\theta\ge0,\nu+\theta>0\} $, is also $ \EL|_{\ML} $-critical, because $ \EL $ is invariant under translations. 
Finally, recall from \cref{rem:explicit-lagrange-multiplier}, that for arclength-parametrized $ \EL $-critical curves $ \gamma\in\ML $ we have the Euler-Lagrange equation $ \D\EL_{\gamma}-\Lambda\D\Sigma_{\gamma}=0 $ with the explicit Lagrange multiplier 
\[\Lambda=\D\EL_{\gamma}Y_{\gamma}\in \left( H^{s}( \R/\Z) \right)^*\equiv H^{-s}(\R/\Z).\]
Consequently, by \cref{eq:log-strain-differential}, we are going to investigate the criticality equation 
\begin{equation}    \label{eq:CriticalityWithoutLambda}
\begin{aligned}
    \D\EL_{\gamma}h=\D\EL_{\gamma}Y_{\gamma}\left\langle \gamma',h' \right\rangle\quad \Foa h\in H^{1+s}(\R/\Z,\R^n).
\end{aligned}
\end{equation}
\subsection{Bending energy}\label{sec:RegularityBending}

\begin{theorem}
\label{thm:elastic-smoothness-new}
Fix a point $\bz\in\R^n$ and let $\UL$ be an open subset of $\NL^2$ or $\NL^2_\bz$.
Then every critical point of $\BendingEnergy$ on $\UL$ is of class $C^\infty(\Circle,\AmbSpace)$.
\end{theorem}
\begin{proof}\Cref{cor:critical-points_wo_fixed_point} allows us to assume that $ \gamma $ is $ \BendingEnergy\vert_{\UL} $-critical for $ \UL\subset \NL^2$. 
    We exploit that $\gamma$ is parametrized by arclength and the formulas \cref{eq:Diff_elastic_formula}
    for $\D(\BendingEnergy)_\gamma$.
    Substituting this into \cref{eq:CriticalityWithoutLambda} for $ \EL\ceq\BendingEnergy $, we get the following equation
    \begin{equation*}
        \D(\BendingEnergy)_\gamma h
        =
        Q(\gamma,h)
        +
        R_\gamma h
        =
        F_\gamma  h
        \quad 
        \text{for all $h \in H^2(\Circle,\AmbSpace)$}
        ,
    \end{equation*}
    \\[-5ex]
    where
    \begin{equation*}
        Q(\gamma,h)
        \ceq 
        \textstyle \int_{\Circle} \inner{\gamma'',h''} \dd u
        ,
        \;\;
        % \\
        R_\gamma  h
        \ceq
        \textstyle
        - \tfrac{3}{2} \int_{\Circle} \abs{\gamma''}^2 \inner{\gamma',h'} \dd u
        ,
        \;\;
        F_\gamma  h       
        \ceq 
        D(\BendingEnergy)_\gamma \, Y_\gamma \inner{\gamma',h'}
        .
    \end{equation*}
    This is the weak formulation of the elliptic equation 
    \begin{equation}
        \Delta^2 \gamma = F_\gamma - R_\gamma.
        \label{eq:Bending:EllipticEquation}
    \end{equation}
    Let $\varepsilon >0$ and suppose that we know already that $\gamma \in H^t(\Circle,\AmbSpace)$ for some $t \geq 2$.
    Integration by parts (in the distributional sense) and \cref{cor:gamma''gamma''} reveal that 
    \begin{equation}\label{eq:95A_vdM}
    \begin{aligned}
        R_\gamma 
        = 
        \tfrac{3}{2}  \pars[\big]{ \inner{\gamma'',\gamma''} \, \gamma'}' \in H^{\sigma(t-2,\varepsilon)-1}(\Circle,\AmbSpace)\quad\Foa\varepsilon>0
        ,
    \end{aligned}
    \end{equation}
    where $\sigma(\cdot,\cdot)$ is the function defined in \cref{lem:DEbYisLowerOrderIIII}.
    Moreover, the second part of that lemma also shows that 
    \begin{equation}\label{eq:95B_vdM}
    \begin{aligned}
        F_\gamma \in H^{\sigma(t-2,\varepsilon)-1}(\Circle,\AmbSpace)\quad\Foa\varepsilon>0
        .
    \end{aligned}
    \end{equation}
    So, if $\gamma \in H^t(\Circle,\AmbSpace)$, then equation \cref{eq:Bending:EllipticEquation} implies 
    \[
        \Delta^2 \gamma \in 
        H^{\sigma(t-2,\varepsilon)-1}(\Circle,\AmbSpace).
    \]
    Now, elliptic regularity for the bi-Laplacian, see \cite[Theorem 6.33]{folland_1995}, implies
    \[
        \gamma \in 
        H^{4 + \sigma(t-2,\varepsilon)-1}(\Circle,\AmbSpace)
        \subset 
        H^{t+1/2 - \varepsilon}(\Circle,\AmbSpace), \quad\Foa\varepsilon>0.
    \]
    So, if we fix $0 < \varepsilon <1/2$, say $ \varepsilon=\frac{1}{4} $, then we get the same strict increase in regularity for each $t \geq 2$.
    Starting at $t = 2$ and iterating this, we obtain
    \[
        \gamma \in \bigcap_{t \geq 2} H^{t}(\Circle,\R^n) = C^\infty(\Circle,\R^n).
    \]
\end{proof}
 
\begin{corollary}[$ E_b $ critical curves]
    Let $ \bz\in \R^n $ and $ \ML\in\{\NL^2,\NL^2_{\bz}\} $.
    Any critical curve $ \gamma $ of the restricted bending energy $ E_b\big\vert_{\ML} $ is critical for $ E_b-\lambda\LL $ on the open subset $ H^2_{\mathrm{r}}\left( \R/\Z,\R^n \right)\subset H^2\left( \R/\Z,\R^n \right)$ for $ \lambda=-E_b(\gamma)\le-(2\pi)^2 $. 
\end{corollary}
\begin{proof}
    Since $ \BendingEnergy $ is parameter-invariant and positively $ (-1) $-homogeneous, we may apply \cref{thm:critical-points}, which shows that $ \gamma $ is critical for $ \BendingEnergy-\lambda\LL $ for $ \lambda=-\BendingEnergy(\gamma) $. Fenchel's theorem \cite[Theorem 3, p. 405]{CarmoDifferentialGeometryCurves2016}, on total curvature together with Hölder's inequality gives $ \BendingEnergy(\gamma)\ge (2\pi)^2 $.
\end{proof}

\subsection{Bending energy with tangent-point energy}\label{sec:RegularityMixed}

\begin{theorem}\label{thm:smoothness-T-critical-new}
Let $\theta>0,p \in \intervaloo{4,5},\bz\in\R^n$, and let $\UL$ be an open subset of $\AL^2$ or $\AL^2_\bz$.
Then every critical point of $\MixedEnergy=\BendingEnergy+\theta\TPp$ on $\UL$ is of class $C^\infty(\Circle,\AmbSpace)$.
\end{theorem}
\begin{proof}As in the proof of \cref{thm:elastic-smoothness-new} we may focus on $ \MixedEnergy\vert_{\UL} $-critical knots $ \gamma $ for $ \UL\subset\AL^{2} $. 
    The criticality equation is a minor modification of \cref{eq:Bending:EllipticEquation}:
    \begin{equation}
        \Delta^2 \gamma = F_\gamma - R_\gamma - \theta \, \D(\TP^{(p,2)})_\gamma.
        \label{eq:Mixed:EllipticEquation}
    \end{equation}
    For $ t\ge 2>{(p-1)}/{2} $ and $ \gamma\in \AL^{t} $ \cite[Proposition 4.3]{blatt-reiter_2015a} implies $ R^{\TP}_{\gamma}\in H^{t-1-p/2-\varepsilon}(\R/\Z,\R^n) $ for all $ \varepsilon>0 $, where $ R^{\TP}_{\gamma} $ appears as a remainder term in \cref{eq:splitting-TP}.
    Combining this for $ \varepsilon\in(0,(p-4)/2] $ with the first part of \cref{lem:EllipticRegularityQTP} and \cref{eq:splitting-TP} leads to    
    \[
        \D(\TP^{(p,2)})_\gamma  \in H^{t-(p-1)}(\Circle,\AmbSpace).
    \]
    Together with \cref{eq:95A_vdM,eq:95B_vdM}, we obtain
    \[
        \Delta^2 \gamma \in H^{\sigma(t-2,\varepsilon)-1}(\Circle,\AmbSpace) \cap H^{t-(p-1)}(\Circle,\AmbSpace)
        \quad\Foa\varepsilon>0.
    \]
    By elliptic regularity, \cite[Theorem 6.33]{folland_1995}, we arrive at
    \[
        \gamma 
        \in
        H^{\sigma(t-2,\varepsilon)+3}(\Circle,\AmbSpace) \cap H^{t-(p-1)+4}(\Circle,\AmbSpace).
    \]
    Since $p-1 < 4$ and by choice of, say $ \varepsilon=\frac{1}{4} $ as in the proof of \cref{thm:elastic-smoothness-new},
    we gain the fixed amount $ \min\{4-(p-1),\frac{1}{4}\} >0$ in regularity for each $ t\ge2 $, and therefore conclude that $ \gamma $ is smooth.
\end{proof}

The previous result in combination with \cref{thm:critical-points} and \cref{cor:critical-points_wo_fixed_point} yields the following corollary.  
\begin{corollary}[Critical knots under length constraint]
    Let $ \bz\in \R^n $ and $ \ML\in\{\AL^2_\bz,\AL^2\} $.
    For any fixed $ \theta>0$ and $p\in(4,5) $ every critical knot $ \gamma $ of the restricted energy $ \TL\big\vert_{\ML} $ is critical for $ \TL-\lambda \, \LL $ on the open subset $ H^2_{\mathrm{ir}}\left( \R/\Z,\R^n \right)\subset H^2\left( \R/\Z,\R^n \right)$ for some $ \lambda\in \R $.     
\end{corollary}

\begin{remark}[Homogeneous total energies]
\label{rem:homogeneous-eneriges}
One could add a different power of the knot energy (depending on its scaling
behavior) to obtain a positively homogeneous total energy as was done
in \cite[(1.4)]{gilsbach-etal_2023}. Here, for example, one could consider
$\SL \ceq \BendingEnergy + \theta \, (\TP^{(p,2)}\circ\iota)^{4-p}$ as a positively $(-1)$-homogeneous
total energy, which implies by means of \Cref{thm:critical-points}
that $\lambda=-\SL(\g)\le -\BendingEnergy(\g)\le -(2\pi)^2.$ If
$\g$ is non-trivially knotted, one even has $\lambda< -(4\pi)^2$
due to the Far\'y-Milnor theorem \cite[Corollary 4.6]{milnor_1950}.
Very recently, the sharp lower bound 
\begin{equation}\label{eq:TPmin}
\begin{aligned}
    \TP_\textnormal{min} \ceq 
\pi^{p-3}\int_0^\pi \sin^{4-p}z\d z
\end{aligned}
\end{equation}
for $\TP^{(p,2)}$ has been established
in \cite{doehrer-dohmen_2025}, which leads to the improved upper bound
$-(2\pi)^2 - \theta \, \TP_\textnormal{min}^{1/(p-4)}$ on $\lambda$, or
even the strict upper bound $-(4\pi)^2 - \theta \, \TP_\textnormal{min}^{1/(p-4)}$,
if $\g$ is non-trivially knotted.
\end{remark}
\subsection{Tangent-point energy}\label{sec:RegularityTP}

\begin{theorem}\label{thm:smoothness-TP-new}
    Let $p\in (4,5)$, $s=s(p) = (p-3)/2$, and  $\bz\in\R^n$. 
    Let $\UL$ be an open subset of the manifolds $\AL^{1+s}$ or $\AL^{1+s}_\bz$.
    Then every critical knot $\g$ of $\TP^{(p,2)}$ on $\UL$ is of class $C^\infty(\Circle,\AmbSpace)$. 
\end{theorem}
\begin{proof}
    By substituting the splitting \cref{eq:splitting-TP} into \cref{eq:CriticalityWithoutLambda} for $ \EL=\TPp $, we 
    get the following criticality equation:
    \[
        2Q^{\TP}(\gamma,h)
        =
        F_\gamma h
        - 
        R^{\TP}_\gamma h
        \quad
        \text{for all $h \in H^{1+s}(\Circle,\AmbSpace)$}
        ,
    \]
    where 
    \[
        F_\gamma  h \ceq \D(\TP^{(p,2)})_\gamma \, Y_\gamma \inner{\gamma',h'}.
    \]
    Now suppose that we know already that $\gamma \in H^t(\Circle,\R^n)$ for some $t \geq 1 + s$.
    Then \cite[Proposition 4.3]{blatt-reiter_2015a} shows that 
    \[
        R^{\TP}_\gamma \in H^{t-1-p/2-\varepsilon}(\Circle,\R^n)\quad\Foa\varepsilon>0.
    \]
    Moreover, by \cref{lem:RegularityTPRemainderII}, 
    
    \[
        F_\gamma \in
        H^{t-1-p/2-\varepsilon}(\Circle,\R^n) \quad\Foa\varepsilon>0
        .
    \]
    Combining this, we see that 
    \[
        Q^{\TP}(\gamma,\cdot) \in H^{t-1-p/2-\varepsilon}(\Circle,\R^n) \Foa \varepsilon>0.
    \]
    By elliptic regularity (see the second part of \cref{lem:EllipticRegularityQTP}), $\gamma$ must have $p-1$ derivatives more than $\QL^{\TP}(\gamma,\cdot)$, thus
    \[
        \gamma \in 
        H^{t-1-p/2-\varepsilon + (p-1)}(\Circle,\R^n)
        =
        H^{t+(p-4)/2 - \varepsilon}(\Circle,\R^n)\quad\Foa\varepsilon>0
        .
    \]
    Since $p \in \intervaloo{4,5}$, we may fix $\varepsilon$ such that $0 < \varepsilon < (p-4)/2$, e.g., $ \varepsilon=(p-4)/{4}$, and get a strict increase in regularity.
    Starting at $t = 1 + s$ and iterating this, we obtain
    \[
        \gamma \in \bigcap_{t \geq 1+s} H^{t}(\Circle,\R^n) = C^\infty(\Circle,\R^n).
    \]
\end{proof}
Since $ \TPp $ is positively $ (4-p) $-homogeneous we obtain length-constrained critical knots as before, with an explicit bound on the scalar Lagrange multiplier. 
\begin{corollary}[$ \TPp $ critical knots]\label{cor:TP_critical_knots}
    Let $ \bz\in\R^n $ and $ \ML\in\{\AL^{1+s},\AL^{1+s}_\bz\} $. 
    Any critical knot $ \gamma $ of the restricted energy $ \TPp\big\vert_{\ML} $ is critical for $ \TPp-\lambda\,\LL $ on the open subset $ H^{1+s}_{\mathrm{ir}}\left( \R/\Z,\R^n \right)\subset H^{1+s}\left( \R/\Z,\R^n \right)$ for $ \lambda=(4-p)\TPp(\gamma)\le(4-p)\TP_{\min} $, where $ \TP_{\min} $ is defined in \cref{eq:TPmin}.
\end{corollary}

\begin{remark}
    In \cite[Corollary 1.3]{steenebruegge-vorderobermeier_2022} it was shown that critical points of integral Menger curvature $ \Mpz $ for $ p\in( \tfrac{3}{7},\tfrac{8}{7} ) $ under a fixed-length constraint are even real-analytic. The technique of that proof, however, does not seem to carry over to arclength-constrained critical knots, but it relies on a splitting of the differential $ \D\Mpz $ provided in \cite[Section 4]{blatt-reiter_2015b}. This splitting is of the same kind as that of $ \D\TPp $, $ p\in(4,5) $, see \cref{eq:splitting-TP}, that we frequently used in the present paper. Thus, the same method applies to the tangent-point energies. As a consequence, length constrained $ \TPp $-critical knots, and therefore by \cref{cor:TP_critical_knots} \textit{also $ \TPp|_{\ML} $-critical knots, are real-analytic for $ \ML\in \{\AL^{1+s},\AL^{1+s}_{\bz}\} $ for any $ \bz\in \R^n $, where $ s=s(p)={(p-3)}/{2} $ $ p\in(4,5) $. }
\end{remark}

%% file: ReparametrizationToArcLength.tex
\section{Reparametrization to arclength}\label{sec:arclength-reparam}

\begin{theorem}\label[theorem]{thm:arclength}
Let $s\in \intervaloc{1/2,1}$. The (rescaling and) reparametrization operator 
\begin{align*}
\ARC: H^{1+s}_\reg\left( \R/\Z,\R^n \right) \to H^{1+s}_\unitspeed\left( \R/\Z,\R^n \right),
\quad	\gamma &\mapsto \frac{1}{\mathcal{L}(\gamma)} \, \gamma\circ\psi^{-1}_\gamma,
\end{align*} 
from the regular $H^{1+s}$-curves to the arclength
parametrized $H^{1+s}$-curves 
is continuous with respect to the $H^{1+s}$ topology, where 
$\psi_{\gamma}(t) \ceq \frac{1}{\mathcal{L}(\gamma)}\mathcal{L}(\gamma\mid_{[0,t]})$ for $t\in [0,1]$.
	\end{theorem}
Similar results have been shown in different settings; 
see, e.g., \cite[Theorem 2]{anosov_1981}. We follow the 
approach of Reiter in \cite[Theorem 1.13]{reiter_2009} who proved this result in the case $s=1$.
	Note that the result stated also holds for non-injective curves.

We prove \Cref{thm:arclength} by showing that the mapping $\ARC$ 
is the composition of 
continuous functions, dealt with in 
three separate lemmas. For that we define
\begin{align*}
	H^{1+s}_\homeom((0,1)) 
	&\ceq 
	\braces[\big]{
		\psi\in H^{1+s}((0,1)) 
		:  
		\psi'(x)>0 \Foa x\in[0,1],\; \psi(0)=0,\; \psi(1)=1
	}
	,
	\\
	H^{1+s}_\periodic((0,1)) 
	&\ceq 
	\braces[\big]{
		\psi\in H^{1+s}_\homeom((0,1)) 
		:  
		\psi'(0)=\psi'(1)
	}.
	\end{align*}
	The space $H^{1+s}_\homeom((0,1))$ is the set of orientation-preserving 
	$H^{1+s}$-diffeomorphisms of  the interval $[0,1]$, and its
	subspace $H^{1+s}_\periodic((0,1))$ consists of those diffeomorphisms, that 
	have a $1$-periodic derivative.
	The first step is to show continuity of the inversion operator: 
	\begin{lemma}
		\label[lemma]{inversioniscontinuous}
		For any $s\in (\frac{1}{2},1]$ the inversion operator $^{-1}:H^{1+s}_\homeom((0,1))\to H^{1+s}_\homeom((0,1))$ 
		is well-defined and continuous. In addition, we have for any 
		$\psi\in H^{1+s}_\homeom((0,1))$, $(\psi^{-1})'=\frac{1}{\psi'\circ\psi^{-1}}$.
	\end{lemma}
	\begin{proof}
For $s=1$ this is proven in 
\cite[Theorem 1.16]{reiter_2009}, so we may assume that $s\in (\frac{1}{2},1)$.
Since $H^{1+s}\equiv H^{1+s}((0,1))$ 
embeds continuously into $C^1([0,1])$, any 
$\psi\in H^{1+s}_e\equiv H^{1+s}_e((0,1))$ is a diffeomorphism
with a  continuously differentiable inverse $\phi \ceq \psi^{-1}$
satisfying $\phi(0)=0$, $\phi(1)=1$, and $\phi'=\frac{1}{\psi'\circ\phi}\ge 
\frac{1}{\left\|\psi' \right\|_{C^0}}>0$. 
By means of the change of variables $x=\psi(w),$ $y=\psi(v)$ with $
\d x=\psi'(w)\d w$ and $\d y=\psi'(v)\d v$ we can express the seminorm
of $\psi'\circ\phi$ as
\begin{align}
\lfloor\psi' &\circ\phi\rfloor_{s,2}^2
\textstyle
=\iint_{[0,1]^2}\frac{\left|\psi'\circ \phi(x) - \psi'\circ \phi (y) \right|^2}{\left|x-y \right|^{1+2s}}\d x\d y 
= \iint_{[0,1]^2}\frac{\left|\psi'(w) - \psi'(v) \right|^2}{\left|\psi(w)-\psi(v) \right|^{1+2s}}\psi'(w)\psi'(u)\d w\d v\notag\\
&
\textstyle
\le  \left\|\psi' \right\|^2_{C^0}\iint_{[0,1]^2}\frac{\left|\psi'(w) 
- \psi'(v) \right|^2}{\left|\psi(w)-\psi(v) \right|^{1+2s}}\d w\d v\le
c_\psi^{-(1+2s)}\left\|\psi' \right\|^2_{C^0}
\lfloor\psi'\rfloor_{s,2}^2<\infty,\label{eq:erste-chain}
\end{align}
where we used the bilipschitz estimate $|\psi(w)-\psi(v)|\ge c_\psi|w-v|$ 
for
all $v,w\in [0,1]$ for the diffeomorphism $\psi\in H^{1+s}_e$.
Since the map $r:[c,\infty)\to(0,\frac{1}{c}],$ $r(t)=\frac{1}{t}$ 
is Lipschitz for any $c>0$, the derivative $\phi'=r\circ\psi'\circ
\phi$ is of class $H^{s}((0,1))$ and therefore $\phi\in H^{1+s}$, so that
the inversion operator is well-defined. 

Now consider a sequence $(\psi_k)_k\subset H^{1+s}_\homeom$ converging 
to $\psi_\infty$ in $H^{1+s}_\homeom$ as $k\to\infty$. 
By virtue of the compact embedding $H^{1+s}\hookrightarrow C^1$ we
may assume without loss of generality 
that 
there is some constant $c>0$ such that $\psi'_k(x)\ge c$ for all $k\in\N$.
We want to prove that  $\phi_k \ceq \psi^{-1}_{k}$ converge to $\phi_\infty \ceq 
\psi^{-1}_\infty$ in $H^{1+s}$ as $k\to\infty$.
With 
\begin{align*}
	\left\|\phi_k-\phi_\infty \right\|_{C^0}
	&=\left\|\phi_\infty\circ\psi_\infty\circ\phi_k-\phi_\infty
	\circ\psi_k\circ\phi_k \right\|_{C^0}
	=\left\|\phi_\infty\circ\psi_\infty-\phi_\infty\psi_k \right\|_{C^0},
\end{align*}
where we used the fact that $\phi_k$ is a homeomorphism of 
$[0,1]$ onto itself, we obtain $
\left\|\phi_k-\phi_\infty \right\|_{C^0}
\le \left\|\phi_\infty' \right\|_{C^0}
\left\|\psi_k-\psi_\infty \right\|_{C^0}\to 0$ as $k\to\infty$.
		Furthermore,  
		\begin{align*}
\left\|\phi'_k-\phi_\infty' \right\|_{C^0}&\textstyle
\le \frac1{c^2}\left\|\psi'_k\circ\phi_k-\psi_\infty'
\circ \phi_\infty \right\|_{C^0}\\
			&
			\textstyle
\le \frac1{c^2}\left( \left\|\psi'_k\circ\phi_k-\psi_\infty'
\circ \phi_k \right\|_{C^0}+\left\|\psi_\infty'
\circ\phi_k-\psi_\infty'\circ \phi_\infty \right\| \right),
		\end{align*}
where the first summand on the right-hand side
equals $\left\|\psi'_k-\psi_\infty' \right\|_{C^0}\to 0$, 
since $\phi_k$ is a diffeomorphism and $\psi_k\to\psi_\infty$ 
in $C^1([0,1])$ as $k\to\infty$. 
The second summand goes to zero, as $\psi_\infty'$ is uniformly 
continuous on $[0,1]$, 
and we already showed that $\phi_k\to\phi_\infty$ 
in $C^0([0,1])$ as $k\to\infty$. 
It remains to show convergence with respect to the fractional Sobolev norm. 
\begin{align}\label{eq:seminorm-basic-estimate}
&\textstyle
\lfloor\phi'_k-\phi'_\infty\rfloor_{s,2}
\textstyle
=  \big\lfloor\frac{\psi'_k\circ\phi_k-\psi_\infty'\circ 
\phi_\infty}{(\psi'_k\circ\phi_k)\cdot 
(\psi_\infty'\circ\phi_\infty)}\big\rfloor_{s,2}\\
			&\le
			\textstyle
\lfloor\psi'_k\circ\phi_k\!-\!\psi_\infty'
\circ \phi_\infty\rfloor_{s,2}
\big\|\frac{1}{(\psi'_k\circ \phi_k)\cdot (\psi_\infty'\circ\phi_\infty)} 
			\big\|_{C^0}
+ \|\psi'_k\circ\phi_k\!-\!\psi_\infty'\circ \phi_\infty\|_{C^0}
\big\lfloor\frac{1}{(\psi'_k\circ\phi_k)
\cdot (\psi_\infty'\circ\phi_\infty)}\big\rfloor_{s,2},\notag
		\end{align}
where we used \cite[Proposition A.5]{knappmann-etal_2023}.
The second summand on the right-hand side of 
\cref{eq:seminorm-basic-estimate} converges to zero as $k\to\infty$, 
since the $C^0$-norm equals $\|\psi_k'-\psi_\infty'\|_{C^0}\to 0$, and
the seminorm is bounded from above by 
\begin{align}\label{eq:upper-bound-second-seminorm}
\textstyle
\frac1{c^2}\lfloor(\psi_k'\circ\phi_k)
\cdot (\psi_\infty'\circ\phi_\infty)\rfloor_{s,2}&
\textstyle
\le\frac1{c^2}\big(
\|\psi_k'\|_{C^0}\lfloor\psi_\infty'\circ\phi_\infty\rfloor_{s,2}+
\lfloor\psi_k'\circ\phi_k\rfloor_{s,2}\|\psi_\infty'\|_{C^0}\big)\\
&
\textstyle
\overset{\cref{eq:erste-chain}}{\le}
\frac1{c^2}\|\psi_k'\|_{C^0}\|\psi_\infty'\|_{C^0}
\big(c_{\psi_\infty}^{-(1+2s)/2}
\lfloor\psi_\infty'\rfloor_{s,2}+
c_{\psi_k}^{-(1+2s)/2}
\lfloor\psi_k'\rfloor_{s,2}\big),\notag
\end{align}
according
to \cite[Lemma A.6]{knappmann-etal_2023} and \cref{eq:erste-chain}, where the 
bilipschitz constants $c_{\psi_k}$ of the $\psi_k$ remain
uniformly positively bounded from below because of the convergence
$\psi_k\to\psi_\infty$ in $C^1$ as $k\to\infty$, and by
the fact that the
diffeomorphism $\psi_\infty$ has a strictly positive bilipschitz
constant $c_{\psi_\infty}$. Therefore, the supremum of 
 \cref{eq:upper-bound-second-seminorm} over $k\in\N$ is finite.

The $C^0$-norm in the first summand of the right-hand side
of \cref{eq:seminorm-basic-estimate} is bounded from above by 
$\frac{1}{c^2}$,
so
it remains to investigate the seminorm in that first summand:   
\begin{equation}
\label{eq:inverse_fracnorm}
\lfloor\psi_k'\circ\phi_k-\psi_\infty'\circ\phi_\infty\rfloor_{s,2}
\le\lfloor\psi_k'\circ\phi_k-\psi_\infty'\circ\phi_k\rfloor_{s,2}+
\lfloor\psi_\infty'\circ\phi_k-\psi_\infty'\circ\phi_\infty\rfloor_{s,2}.
\end{equation}
The first term in \cref{eq:inverse_fracnorm}
goes to zero as $k\to\infty$ since 
\begin{align*}
&\textstyle
\lfloor\psi'_{k}\circ\phi_k-\psi_\infty'\circ\phi_k\rfloor_{s,2}^2
=\iint_{[0,1]^2}\frac{|\psi'_{k}(\phi_k(x))-\psi_\infty'(\phi_k(x))
-\psi'_{k}(\phi_k(y))+\psi_\infty'(\phi_k(y))|^2}{\left|x-y \right|^{1+2s}}\d x\d y \\
&
\textstyle
= \iint_{[0,1]^2}\frac{|\psi'_{k}(v)-\psi_\infty'(v)-\psi'_{k}(w)
+\psi_\infty'(w)|^2}{\left|\psi_k(w)-\psi_k(v) \right|^{1+2s}} 
\psi'_{k}(w) \psi'_{k}(v)\d w\d v
\le c_{\psi_k}^{-(1+2s)}\|\psi'_k \|^2_{C^0}
\lfloor \psi_k'\!-\!\psi_\infty'\rfloor_{s,2}^2
\to 0,
\end{align*}
where we used the 
change of variables $w=\phi_k(x)$ and 
$v=\phi_k(y)$ and the fact that the bilipschitz constants $c_{\psi_k}$
are uniformly positively bounded from below.

For the square of the
second summand in \cref*{eq:inverse_fracnorm} we estimate 
the integrand of the corresponding double integral as
\begin{align*}
&
\textstyle
\frac{|\psi_\infty'(\phi_k(x))-\psi_\infty'(\phi_\infty(x))-(\psi_\infty'
(\phi_k(y))
-\psi_\infty'(\phi_\infty(y)))|^2}{\left|x-y \right|^{1+2s}}
\textstyle
\le \frac{|\psi_\infty'(\phi_k(x))-\psi_\infty'(\phi_k(y))|^2}{\left|x
-y \right|^{1+2s}}+\frac{|\psi_\infty'(\phi_\infty(x))
-\psi_\infty'(\phi_\infty(y))|^2}{\left|x-y \right|^{1+2s}}.
\end{align*}
The last summand is integrable since it equals the integrand of
the squared seminorm $\lfloor\psi_\infty'\circ\phi_\infty\rfloor_{s,2}$,
which can be estimated as 
 in \cref{eq:erste-chain},
and for the first term we have for $x\not = y$ (by setting $w=\phi_k(x)$ and $v=\phi_k(y)$)
\begin{align*}
\textstyle
\frac{|\psi_\infty'(\phi_k(x))-\psi_\infty'(\phi_k(y))|^2}{\left|x-y \right|^{1+2s}}
&
\textstyle
=\frac{|\psi_\infty'(w)-\psi_\infty'(v)|^2}{\left|\psi_k(w)-\psi_k(v) \right|^{1+2s}}
\le
c_{\psi_k}^{-(1+2s)}\frac{|\psi_\infty'(w)
-\psi_\infty'(v)|^2}{\left|w-v \right|^{1+2s}}
\le C\frac{|\psi_\infty'(w)
-\psi_\infty'(v)|^2}{\left|w-v \right|^{1+2s}},
\end{align*}
where we used again that the bilipschitz constants $c_{\psi_k}$ are 
uniformly positively bounded from below.
The last term is integrable since it is the integrand of
the squared seminorm $\lfloor\psi_\infty'\rfloor_{s,2}$.
So we found an integrable majorant for the integrand almost everywhere 
(as the set $\{(x,y)\in[0,1]^2 \mid x\not=y\}$ is a set of measure zero). 
By the dominated convergence theorem we have 
$\lfloor \psi_\infty'\circ\phi_k-\psi_\infty'\circ\phi_\infty
\rfloor_{s,2}\to 0$, and therefore by \cref{eq:seminorm-basic-estimate}
also $\lfloor\phi_k'-\phi_\infty'\rfloor_{s,2}\to 0$ as $k\to\infty$.
\end{proof}
Analogously to \cite[Corollary 1.17]{reiter_2009} we have the following corollary:
\begin{corollary}\label{cor:restriction-inverse-operator}
The restriction of the inversion operator $^{-1}$ to $H^{1+s}_\periodic((0,1))$ 
is a continuous map onto $H^{1+s}_\periodic((0,1))$.
\end{corollary}
Next we show that the composition operator is well-defined and continuous. 
\begin{lemma}
\label[lemma]{composition_continuous}
Let $s\in (\frac{1}{2},1]$. The composition operator 
\begin{align}\label{eq:composition-operator}
H_p^{1+s}((0,1))\times H^{1+s}(\R/\Z,\R^n)\to H^{1+s}(\R/\Z,\R^n),\quad
(\psi,h)\mapsto h\circ\psi,
\end{align}
is well-defined and continuous. 
\end{lemma}
\begin{proof}
For $s=1$ one finds a proof in \cite[Theorem 1.18]{reiter_2009}, so
as before we focus on $s\in (\frac{1}{2},1)$. The chain rule \cite[Lemma A.4]{knappmann-etal_2023} yields $h\circ\psi\in H^{1+s}(\R/\Z,\R^n)$, so that the composition operator \cref{eq:composition-operator}
is well-defined.

As to continuity, 
consider a sequence $(\psi_k,h_k)_k\subset H_p^{1+s}((0,1))
\times H^{1+s}(\R/\Z,\R^n)$ converging in $H^{1+s}\times H^{1+s}$ to a limiting 
pair $(\psi_\infty,h_\infty)\in H_p^{1+s}((0,1))\times H^{1+s}(\R/\Z,\R^n)$.
Then we have 
\begin{align*}
\textstyle
\|h_k\circ\psi_k-h_\infty\circ\psi_\infty\|_{C^0}
&\textstyle
\le\|h_k\circ\psi_k-h_\infty\circ\psi_k\|_{C^0}
+\|h_\infty\circ\psi_k-h_\infty\circ\psi_\infty\|_{C^0}\\
&\le \|h_k-h_\infty\|_{C^0}+\|h_\infty'\|_{C^0}
\|\psi_k-\psi_\infty\|_{C^0}\to 0,
\end{align*}
where we used, that $\psi_k$ is a diffeomorphism for each $k\in\N$. 
Similarly, for the  derivative,
\begin{align}
\|(h_k\circ\psi_k)'- (h_\infty\circ\psi_\infty)'\|_{C^0}&
\le\|\big(\|h_k'-h_\infty'\|_{C^0}
+\|h_\infty'\circ\psi_k-h_\infty'\circ\psi_\infty\|_{C^0} \big)\|\psi_k'\|_{C^0}
\notag\\
&+\|h_\infty'\circ\psi_\infty\|_{C^0}\|\psi_k'-\psi_\infty'\|_{C^0}\to 0\quad\As
k\to\infty,\label{eq:zwischen}
\end{align}
since
$h_\infty'$ is uniformly continuous on $\R/\Z$ and 
$\psi_k\to\psi_\infty$ in $C^1$. 
Lastly, 
the seminorm $\lfloor(h_k\circ\psi_k)'-(h_\infty\circ\psi_\infty)'\rfloor_{s,2}
$ is bounded from above by
\begin{align*}
&\lfloor(h_k'\circ\psi_k)\psi_k'-(h_\infty'\circ \psi_\infty)\psi_k'
\rfloor_{s,2}
+\lfloor(h_\infty'\circ\psi_\infty)\psi_k'-(h_\infty'\circ \psi_\infty)
\psi_\infty\rfloor_{s,2}\\
&\le
\big(\|h_k'-h_\infty'\|_{C^0}
-\|h_\infty'\circ\psi_k-h_\infty'\circ\psi_\infty\|_{C^0}\big)
\lfloor\psi_k'\rfloor_{s,2}\\
& + \big(
\lfloor h_k'\circ\psi_k-h_\infty'\circ\psi_k\rfloor_{s,2}+
\lfloor h_\infty'\circ\psi_k-h_\infty'\circ\psi_\infty\rfloor_{s,2}\big)
\|\psi_k'\|_{C^0}\\
& + \|h_\infty'\circ\psi_\infty\|_{C^0}\lfloor\psi_k'-\psi_\infty'\rfloor_{s,2}
+\lfloor h_\infty'\circ\psi_\infty\rfloor_{s,2}\|\psi_k'-\psi_\infty'\|_{C^0},
\end{align*}
by means of the product rule \cite[Proposition A.5]{knappmann-etal_2023}.
The $k$-dependent $C^0$-norms in the first and third line of the right-hand side
converge to zero as $k\to\infty$
as seen in 
\cref{eq:zwischen}, and the seminorms $\lfloor\psi_k'\rfloor_{s,2}$
 are uniformly bounded in $k$, as well as
$\lfloor h_\infty'\circ\psi_\infty\rfloor_{s,2}$ and
the norm $\|h_\infty'\circ\psi_\infty\|_{C^0}$.
Since also the norms $\|\psi_k'\|_{C^0}$ are uniformly bounded, it suffices to
show that the $k$-dependent seminorms in the middle line 
converge to zero. The first
seminorm can be treated as the first summand on the right-hand side
of \cref{eq:inverse_fracnorm} with the only difference that now the change
of variables reads as $w=\psi_k(x)$ and $v=\psi_k(y)$, so that one
has to use the bilipschitz constants $c_{\phi_k}=\|\psi_k'\|_{C^0}^{-1}$
of the inverse functions $\phi_k=\psi_k^{-1}$ satisfying
$c_{\phi_k}\ge
\frac{1}{2}\|\psi_\infty'\|_{C^0}>0$ for $k\gg 1$, because $\psi_k\to\psi_\infty$ in $C^1$. 
The second seminorm in that middle line resembles the second
summand of \cref{eq:inverse_fracnorm} and can be treated analogously 
using dominated
convergence, only now with the substitution $w=\psi_k(x)$ and $v=\psi_k(w)$ which
again requires the bilipschitz constants $c_{\phi_k}$ which are uniformly
positively bounded from below as $k\to\infty$. 
	\end{proof}
\begin{lemma}
\label[lemma]{lem:parametrization_continuous}
The operator 
\begin{equation}\label{eq:arclength-operator}
\textstyle
\varPsi \colon H_\reg^{1+s}(\R/\Z,\R^n)\to H^{1+s}_\periodic((0,1)),\quad
\gamma\mapsto \big( t\mapsto \frac{1}{\mathcal{L}(\gamma)} \int_{0}^{t}\left|\gamma'(\tau) \right|\d \tau \big),
\end{equation}
is well-defined and continuous.
	\end{lemma}
	\begin{proof}
It is clear that $\varPsi \, \g \in H^{1+s}_p((0,1))$ for any 
$\g\in H^{1+s}_\reg(\R/\Z,\R^n)$ so that $\varPsi$ is well-defined.
For a sequence of curves $(\gamma_k)_{k\in\N}\subset H_\reg^{1+s}(\R/\Z,\R^n)$ 
with $\gamma_k\to\gamma_\infty$ in $H_\reg^{1+s}(\R/\Z,\R^n)$ 
as $k\to\infty$, 
there exists a constant
 $c>0$, such that we have
 $\min\{|\gamma_k'(x)|,|\g_\infty'(x)|\}\ge c$ 
 for all $x\in\R/\Z$ and $k\in\N$, which implies $\LL(\g_k)\ge c$ for
 all $k\in\N$.
 In addition, $\psi_k \ceq \varPsi \, \gamma_k$ is of class  
 $C^1$ with $\psi_k'=\frac{1}{\mathcal{L}(\gamma_k)}\left|\gamma_k' \right|$
 for each $k\in\N$. 
Furthermore,  
$\mathcal{L}(\gamma_k)\to\mathcal{L}(\gamma_\infty)$ as $k\to\infty$
by the $C^1$-convergence $\g_k\to\g_\infty$; hence
\begin{align*}
&
\textstyle
\|\psi_k(\cdot)-\psi_\infty\|_{C^0}
= \norm[\Big]{ 
	\frac{\mathcal{L}(\gamma_\infty)\int_{0}^{\cdot}|\gamma_k'|\d\tau
	-\mathcal{L}(\gamma_k)\int_{0}^{\cdot}|\gamma_\infty'|\d\tau}{\mathcal{L}(
	\gamma_k)\mathcal{L}(\gamma_\infty)} 
}_{C^0}\\
&
\textstyle
\le\frac{1}{c}\int_{0}^{1}|\gamma_k'(\tau)-\gamma_\infty'(\tau)|\d\tau 
+ \frac{1}{c^2}| \mathcal{L}(\gamma_k)-\mathcal{L}(\gamma_\infty)|
\int_{0}^{1}|\gamma_\infty'(\tau)|\d\tau\to 0\quad\As k\to\infty.
\end{align*}
For the $C^0$-norm of the derivative we get
\begin{align*}
\textstyle
\|\psi_k'-\psi_\infty'\|_{C^0}&
\textstyle
=\big\|\frac{1}{\mathcal{L}(\gamma_k)}
|\gamma_k'| -\frac{1}{\mathcal{L}(\gamma_\infty)}|\gamma_\infty'|\big\|_{C^0}
\textstyle
\le\frac{1}{c}\|\gamma_k'-\gamma_\infty'\|_{C^0}
+\frac{1}{c^2}|\mathcal{L}(\gamma_k)-\mathcal{L}(\gamma_\infty)|
\|\gamma_\infty'\|_{C^0},
\end{align*}
which converges to zero as $k\to\infty$. Analogously,
\begin{align*}
\textstyle
\lfloor\psi_k'-\psi_\infty'\rfloor_{s,2}\le\frac{1}{c} 
\lfloor\gamma_k'-\gamma_\infty'\rfloor_{s,2}+\frac{1}{c^2}
|\mathcal{L}(\gamma_k)-\mathcal{L}(\gamma_\infty)|\lfloor\gamma_\infty'
\rfloor_{s,2}\to 0\quad\As k\to\infty. 
\end{align*}
	\end{proof}
\begin{refproof}[Proof of \Cref{thm:arclength}]
By 
\Cref{inversioniscontinuous,composition_continuous,lem:parametrization_continuous} 
the operator is continuous as the composition of continuous operators. 
Furthermore, for $\gamma\in H^{1+s}_\reg(\R/\Z,\R^n)$, 
\begin{align*}
\textstyle
|\ARC(\gamma)'(x)|&
\textstyle
=\big|\frac{1}{\mathcal{L}(\gamma)}
\left( \gamma\circ(\varPsi \, \gamma)^{-1}\right)'(x) \big|
=\frac{1}{\mathcal{L}(\gamma)}\left|\gamma'\left( \left( \varPsi \, \gamma \right)^{-1} (x)\right)(\left( \varPsi \, \gamma \right)^{-1})'(x) \right|\\
&
\textstyle
=\frac{1}{\mathcal{L}(\gamma)}\big|\frac{\gamma'\left( (\varPsi \, \gamma)^{-1}(x) \right)}{(\varPsi \, \gamma)'((\varPsi \, \gamma)^{-1}(x))} \big|=1.
\end{align*}
\end{refproof}

\begin{lemma}[Triple integrals]
\label{lem:triple-integrals}
Let $s\in (0,1)$, $\rho\in [1,\infty)$, and $\g\in W^{1+s,\rho}(\R/\Z,\R^n)$.
Then
\begin{equation}\label{eq:triple-integrals}
\int_{\R/\Z}\int_{-\frac{1}{2}}^\frac{1}{2}\int_0^1\frac{|\g'(u+\theta w)-\g'(u)|^\rho}{
|w|^{1+s\rho}}\d\theta\d w\d u\le \lfloor\g'\rfloor_{s,\rho}^\rho.
\end{equation}
\end{lemma}
\begin{proof}
Using Fubini and the change of variable $\tilde{w} \ceq \theta w$ with $\d \tilde{w}=
\theta\d w$ and $\tilde{w}(\pm\frac{1}{2})=\pm\frac{\theta}2$ one can bound 
the
triple integral in \cref{eq:triple-integrals} by means of $\theta\in [0,1]$
from above by
\begin{equation}\label{eq:triple-integral-upper-bound}\textstyle
\int_0^1\!\!\int_{\R/\Z}\int_{-\frac{1}{2}}^\frac{1}{2}\frac{|\g'(u+\tilde{w})-\g'(u)|^\rho}{|\tilde{w}|^{1+s\rho}}\d \tilde{w}\d u\d\theta\le 
\lfloor\g'\rfloor_{s,\rho}^\rho.
\end{equation}
\end{proof}

%% file: ProductRule_BehzadanHolst.tex
\section{Product rules}\label{sec:ProductRule}

Here we gather a few product rules, i.e., statements of well-definedness and existence of bilinear compositions of functions in certain Sobolev classes. We use these results throughout the main body of the work, in particular, for the regularity theory in \cref{sec:Smoothness}.

We start with the following result, which can be found as Theorem~7.4 in \cite{zbMATH07447116} but is included here for the reader's convenience.
Notice that this theorem for bounded Lipschitz domains can be generalized to compact smooth manifolds such as $\Circle$ by the usual tricks of the trade such as choosing a smooth atlas and subordinate smooth partition of unity.
\begin{theorem}\label{theo:ProductRuleBehzadanHolst}
    Let $\varOmega \subset \DomSpace$ be a bounded Lipschitz domain.
    Let $ \sigma,\sigma_1,\sigma_2\in\R$ with $ \sigma\ge0 $, and $p,p_1,p_2 \in \intervalco{1,\infty}$ satisfy 
    one of the following two conditions for each $i \in \braces{1,2}$:
    \begin{itemize}
        \item
        $p_i \leq p$,
        $\sigma_i \geq \sigma$,
        $\pars[\big]{\sigma_i - \tfrac{\DomDim}{p_i}} \geq \pars[\big]{ \sigma - \tfrac{\DomDim}{p}}$,
        and
        $\pars[\big]{\sigma_1 - \tfrac{\DomDim}{p_1}} + \pars[\big]{\sigma_2 - \tfrac{\DomDim}{p_2}}
        >
        \pars[\big]{\sigma - \tfrac{\DomDim}{p}}
        $.
        \item          
        $\max \braces{p_1,p_2} > p$,
        $\sigma_i > \sigma$,
        $\pars[\big]{\sigma_i - \tfrac{\DomDim}{p_i}} 
        > 
        \pars[\big]{ \sigma - \tfrac{\DomDim}{p}}$,
        and
        $\sigma_1 + \sigma_2 - \sigma
        >
        \tfrac{\DomDim}{\min \braces{p_1,p_2}}$.
    \end{itemize}
    Then the product map
    \[
        W^{\sigma_1,p_1}(\Domain) \times W^{\sigma_2,p_2}(\Domain)
        \to 
        W^{\sigma,p}(\Domain),
        \;\;
        (u_1,u_2) \mapsto u_1 \cdot u_2
    \]
    is well-defined and continuous.
\end{theorem}

\begin{lemma}\label{lem:ProductRuleAsDistribution}
    Let $t \geq 0$. 
    Then for all $\varepsilon > 0$ the map
    \[
        H^{t}(\Circle) \times H^{t}(\Circle)
        \to 
        H^{\sigma(t,\varepsilon)}(\Circle),
        \;\;
        (u_1,u_2) \mapsto u_1 \cdot u_2
    \]
    is well-defined and continuous, where (as in \cref{lem:DEbYisLowerOrderIIII})
    \[
        \sigma(\tau,\varepsilon) \ceq 
        \begin{cases}
            2\,\tau - \frac{1}{2} - \varepsilon ,&\text{for $ \tau \in [0,\frac{1}{2}]$},\\
            \tau                                ,&\text{for $\tau \in(\frac{1}{2},\infty)$}.
        \end{cases}
    \]
    Here, we regard $ u_1\cdot u_2 $ as a distribution whenever $ \sigma(t,\varepsilon)<0 $.
\end{lemma}

\begin{proof}
    If $ t>\frac{1}{2} $, $ H^{t}(\R/\Z) $ is a Banach algebra, since then $ H^{t}\hookrightarrow L^\infty $, and there is nothing to prove. So, we may focus on $ t\in[0,\frac{1}{2}] $ and distinguish three cases. 

    \textit{{Case I:} $t \in (\frac{1}{4},\frac12]$.}
    In this case, $ \sigma(t,\varepsilon)=2t-\frac{1}{2}-\varepsilon\le \frac{1}{2}-\varepsilon $ for all $ \varepsilon>0 $, and we may choose $ \varepsilon_{0}(t)\ceq\frac{1}{2}(2t-\frac{1}{2})>0 $, such that $ \sigma(t,\varepsilon)>0 $ for all $ \varepsilon\in (0,\varepsilon_{0}(t)] $. We may therefor apply \cref{theo:ProductRuleBehzadanHolst} for $ \sigma_{1}=\sigma_2\ceq t $, and $ \sigma\ceq \sigma(t,\varepsilon)=t+(t-\frac{1}{2})-\varepsilon \le t-\varepsilon<t$ for $  \varepsilon\in (0,\varepsilon_{0}(t)] $, and $ p_1=p_2=p\ceq 2 $, to obtain the result first for all $  \varepsilon\in (0,\varepsilon_{0}(t)] $. For $ \varepsilon>\varepsilon_{0}(t) $ simply observe that $ \sigma(t,\varepsilon) <\sigma(t,\varepsilon_{0}(t))$ so that $ H^{\sigma(t,\varepsilon_{0}(t))}(\R/\Z)\subset H^{\sigma(t,\varepsilon)}(\R/\Z) $. 

    \textit{{Case II:} $ t \in  (0,1/4]$.}
    In this case $ \sigma\equiv\sigma(t,\varepsilon)\le -\varepsilon <0$, and $ u_1\cdot u_2 $ acts on smooth functions as 
    \[\textstyle\int_{\R/\Z}u_1\cdot u_2\cdot v\dd x\quad \Fo v\in C^\infty(\R/\Z).\]
    By the Sobolev embedding theorem (see, e.g., \cite[Theorem 2.3.4]{AgranovichSobolevSpaces2015})
    $ H^{-\sigma}(\R/\Z)\hookrightarrow L^{q}(\Circle) $ for all $ q\in [1,2/(1+2\sigma)]= [1,1/(2t-\varepsilon)]$, whereas $ H^{t}(\R/\Z)\hookrightarrow L^{2/(1-2t)} (\Circle)$. 

    In view of Hölder's inequality with three factors we compute for $ p=\frac{2}{1-2t} $, $ \frac{1}{p}+\frac{1}{p}=1-2t $, to find $ q=\frac{1}{2t}\in [1,1/(2t-\varepsilon)]$, so that 
    \begin{align*}\textstyle
        \left|\int_{\R/\Z}u_1\cdot u_2\cdot v\d x \right|\le \left\|u_1 \right\|_{L^p}\left\|u_2 \right\|_{L^p}\left\|v \right\|_{L^q}\le C \left\|u_1 \right\|_{H^t}\left\|u_2 \right\|_{H^t}\left\|v \right\|_{H^{-\sigma}}\Foa v\in C^\infty(\R/\Z),
    \end{align*}
    which by density of $ C^\infty(\R/\Z) $ in $ H^{-\sigma}(\R/\Z) $ (see, e.g. \cite[Statement 4 in Section 2.3]{AgranovichSobolevSpaces2015}) yields the desired estimate for all $ v\in H^{-\sigma}(\R/\Z) $.

    \textit{{Case III:} $t=0$.}
    Now, $ \sigma(t,\varepsilon)=-\frac{1}{2}-\varepsilon <-\frac{1}{2}<0$, so we can use the known embedding $ H^{-\sigma}(\R/\Z)\hookrightarrow L^\infty(\R/\Z) $ to establish 
    \begin{align*}\textstyle
        \left|\int_{\R/\Z}u_1\cdot u_2\cdot v\d x \right|\le \left\|u_1 \right\|_{L^2}\left\|u_2 \right\|_{L^2}\left\|v \right\|_{L^\infty}\le C \left\|u_1 \right\|_{H^t}\left\|u_2 \right\|_{H^t}\left\|v \right\|_{H^{-\sigma}}\Foa v\in H^{-\sigma}(\R/\Z),
    \end{align*}
    which concludes the proof. 
\end{proof}
\begin{corollary}\label{cor:gamma''gamma''}
    Let $t \geq 2$. 
    Then for all $\varepsilon > 0$ the maps 
    \begin{alignat*}{2}
        H^{t}(\Circle,\R^n)
        &\to 
        H^{\sigma(t-2,\varepsilon)}(\Circle),
        \quad
        &\gamma &\mapsto \inner{\gamma'',\gamma''},
        \\
        H^{t}(\Circle,\R^n)
        &\to 
        H^{\sigma(t-2,\varepsilon)}(\Circle,\AmbSpace),
        \quad
        &\gamma &\mapsto \inner{\gamma'',\gamma''} \, \gamma',
        \\
        H^{t}(\Circle,\R^n)
        &\to 
        H^{\sigma(t-2,\varepsilon)-1}(\Circle,\R^n),
        \quad
        &\gamma &\mapsto \pars[\big]{ \inner{ \gamma'',\gamma''} \, \gamma'}',
    \end{alignat*}
    are well-defined and continuous, where $\sigma(\cdot,\cdot)$ is the function from \cref{lem:ProductRuleAsDistribution}.
\end{corollary}
\begin{proof}
    The statement for the first map is an immediate consequence of \cref{lem:ProductRuleAsDistribution}. For the second map we note, that it suffices to prove the statement for all $ \varepsilon\in (0,\frac{1}{2}) $ by the ordering $ H^{\sigma(t-2,\varepsilon_1)}(\R/\Z)\subset H^{\sigma(t-2,\varepsilon_2)}(\R/\Z)$ for all $ 0<\varepsilon_1\le \varepsilon_2 $ and all $ t\ge 2 $. \\
    We distinguish 3 cases. 

    \textit{{Case I:} $t\in (2+\frac{1}{2},\infty)$.}
    In this case, $ \sigma(t-2,\varepsilon)=t-2 \,$ by definition of $ \sigma(\cdot,\cdot) $, and $ H^{t-2}(\R/\Z) $ is a Banach algebra. Since $ \gamma'\in H^{t-1}(\R/\Z,\R^n)\subset H^{t-2}(\R/\Z,\R^n) $ we conclude. 

    \textit{{Case II:} $t\in [2+\frac{1}{4}+\frac{\varepsilon}{2},2+\frac{1}{2}]$.}
    Now $ \sigma(t-2,\varepsilon)\ge0 $ for all $ \varepsilon\in(0,1/2) $, and we can use \cref{theo:ProductRuleBehzadanHolst} with the choice $ \sigma\ceq \sigma(t-2,\varepsilon) $, $ \sigma_1\ceq \sigma, \sigma_2\ceq t-1>1>\frac{1}{2}-\varepsilon\ge \sigma $, and $ p_1=p_2=p \ceq 2$ to settle this case. 

    \textit{{Case III:}\, $t\in [2,2+\frac{1}{4}+\frac{\varepsilon}{2}]$.} Here, $ \sigma(t-2,\varepsilon)<0 $, and we distinguish two subcases depending on $ \varepsilon\in (0,\frac{1}{2}) $.

    \textit{{Case IIIa:} $\varepsilon\in (2(t-2),\frac{1}{2})$, if this interval is non-empty.}
    Then $ \sigma(t-2,\varepsilon) <-\frac{1}{2}$ such that $ -\sigma\ceq-\sigma(t-2,\varepsilon)>\frac{1}{2} $, which implies by Morrey's embedding $ H^{-\sigma}(\R/\Z,\R^n)\hookrightarrow C^0(\R/\Z,\R^n) $ and thus,
    \begin{align*}\textstyle
        \left|\int_{\R/\Z}\left|\gamma'' \right|^2 \left\langle \gamma',v \right\rangle \d x \right|&\le \left\|\gamma'' \right\|^2_{L^2}\left\|\gamma' \right\|_{C^0}\left\|v \right\|_{L^\infty}
        \le \left\|\gamma \right\|^3_{H^t}\left\|v \right\|_{H^{-\sigma}}.
    \end{align*}

    \textit{{Case IIIb:} $\varepsilon\in (0,2(t-2))$, if this interval is non-empty.}
    Now $ \sigma(t-2,\varepsilon)\in[-\frac{1}{2},0) $, and we can apply \cref{theo:ProductRuleBehzadanHolst} to the product $ \left\langle \gamma',v \right\rangle  $ for $ v\in H^{-\sigma(t-2,\varepsilon)}(\R/\Z,\R^n) $ with the choice $ \sigma_1\ceq t-1\ge 1>\sigma_2\equiv\sigma\ceq-\sigma(t-2,\varepsilon) \in (0,\frac{1}{2}]$ to find $ \left\langle \gamma',v \right\rangle \in H^{-\sigma(t-2,\varepsilon)}(\R/\Z,\R^n) $ with the estimate 
   \begin{equation}\label{eq:C4_abschaetzung}
   \begin{aligned}
    \left\|\left\langle \gamma',v \right\rangle  \right\|_{H^{-\sigma(t-2,\varepsilon)}}\le C\left\|\gamma' \right\|_{H^{t-1}}\left\|v \right\|_{H^{-\sigma(t-2,\varepsilon)}}
        \le C\left\|\gamma \right\|_{H^{t}}\left\|v \right\|_{H^{-\sigma(t-2,\varepsilon)}}.
   \end{aligned}
   \end{equation}
    We already know that $ \left\langle \gamma'',\gamma'' \right\rangle \in H^{\sigma(t-2,\varepsilon)} $ with $ \left\|\left\langle \gamma'',\gamma'' \right\rangle  \right\|_{H^{\sigma(t-2,\varepsilon)}}\le \tilde{C}\left\|\gamma \right\|_{H^t}^2 $, which implies 
    \begin{align*}\textstyle
        \left|\int_{\R/\Z}\left|\gamma'' \right|^2f\d x \right|\le \tilde{C}\left\|\gamma \right\|^2_{H^t}\left\|f \right\|_{H^{-\sigma(t-2,\varepsilon)}}\Foa f\in H^{\sigma(t-2,\varepsilon)}(\R/\Z,\R^n);
    \end{align*}
    hence by \cref{eq:C4_abschaetzung} 
    \begin{align*}\textstyle
        \left|\int_{\R/\Z}\left|\gamma'' \right|^2 \left\langle \gamma',v \right\rangle \d x \right|\le \tilde{C}C\left\|\gamma\right\|^3_{H^{t}}\left\|v \right\|_{H^{-\sigma(t-2,\varepsilon)}} \Foa v\in H^{-\sigma(t-2,\varepsilon)}(\R/\Z,\R^n). 
    \end{align*}
\end{proof}

\begin{lemma}\label{lem:ProductRuleAsDistributionII}
    Let $t \geq 1$.
    Then for all $\varepsilon >0$ and all $\tau \leq \sigma(t-1,\varepsilon)$
    the following map
    \[
        H^t(\Circle,\AmbSpace)
        \times
        H^{1-\tau}(\Circle,\AmbSpace) \to H^{1-t}(\Circle),
        \quad 
        (\gamma,h) \mapsto \inner{\gamma',h'}
    \]
    is well-defined and continuous,
    where $\sigma(\cdot,\cdot)$ is the function defined in \cref{lem:ProductRuleAsDistribution}.
\end{lemma}
\begin{proof}
    From \cref{lem:ProductRuleAsDistribution} we know that the following map is well-defined and continuous for all $ \tau\le \sigma(t-1,\varepsilon),\varepsilon>0 $:
    \begin{align*}
        H^{t}(\R/\Z,\R^n)\times H^{t-1}(\R/\Z)\to H^{\tau}(\R/\Z,\R^n), (\gamma,f)\mapsto f\gamma'.
    \end{align*}
    Hence,  there is a constant $ C\ge 0 $ such that 
    \begin{equation}\label{eq:C5_abschaetzung}
    \begin{aligned}
        \left\|f\gamma' \right\|_{H^{\tau}}\le C\left\|f \right\|_{H^{t-1}}\left\|\gamma \right\|_{H^t},\Foa \gamma\in H^t(\R/\Z,\R^n),\,f\in H^{t-1}(\R/\Z),
    \end{aligned}
    \end{equation}
    for all $ \tau\le \sigma(t-1,\varepsilon)$ and $\varepsilon>0$.
    Since $ \tau $ may be negative, we interpret $ f\gamma' $ for a given $ f\in C^\infty(\R/\Z) \subset H^{t-1}(\R/\Z)$ as a distribution and obtain for $ h\in C^\infty(\R/\Z,\R^n) $ that 
    \begin{align*}\textstyle
        \int_{\R/\Z}\left\langle f\gamma',h' \right\rangle \d x =\int_{\R/\Z}\left\langle \gamma',h' \right\rangle f\d x
    \end{align*}
    may be interpreted as $ \left\langle \gamma',h' \right\rangle  $ acting on $ f $ as a distribution, satisfying the following estimate by means of \cref{eq:C5_abschaetzung},
    \begin{align*}\textstyle
        |\int_{\R/\Z}\left\langle \gamma',h' \right\rangle f\d x |\le \left\|f\gamma' \right\|_{H^{\tau}}\left\|h' \right\|_{H^{-\tau}}&\le C \left\|f \right\|_{H^{t-1}}\left\|\gamma \right\|_{H^{t}}\left\|h' \right\|_{H^{-\tau}}
        \le C \left\|f \right\|_{H^{t-1}}\left\|\gamma \right\|_{H^{t}}\left\|h \right\|_{H^{1-\tau}}
    \end{align*}
    for all $ f\in C^\infty(\R/\Z) $ and $ h\in C^\infty(\R/\Z,\R^n)$. By density \cite[Statement 4, Section 2.3]{AgranovichSobolevSpaces2015} this holds for all $ f\in H^{t-1}(\R/\Z) $ and $ h\in H^{1-\tau}(\R/\Z,\R^n) $, which proves the claim. 
\end{proof}

%% file: Acknowledgments.tex
\section*{Acknowledgments}
N.~F.{}, H.~S.{}, and D.~S.{} were  
partially funded by  {the DFG}-Graduiertenkolleg \emph{Energy, Entropy, and Dissipative Dynamics (EDDy)},  project no. 320021702/GRK2326.
The fourth author's work is partially funded by 
the Excellence Initiative of the German federal and state governments.